,

# Yoneda lemma

**Gintaras Valiukevičius**
*Veja*
*Ukmergė, Lithuania*
*e-mail: arklys@vivaldi.net*

**Abstract.**
We are studying properties of the name appointment in various categories of enriched graphs. The Yoneda lemma is generalized for continuous transforms between transports of enriched original graphs.

We are going to define the joint pair of superfunctors $X \otimes T$ and $Y^T$ in the category $Cat_V$ of all $V$-enriched categories. For a values category $V$ we demand existence of *truly joint pair* of functors. One adjoint is defined with the biproduct $r \otimes s$ and another coadjoint is an exponential functor $t^s$. As usually it provides the bijective *name appointment*

$$V(r \times s; t) \to V(r; t^s)$$

with commuting diagrams for arrows $u : r \to r'$ and $v : t \to t'$

$$\begin{array}{ccc} V(r \otimes s; t) & \to & V(r; t^S) \\ (u \otimes 1_s)^* \uparrow & & \downarrow (v^s)_* \\ V(r' \otimes s; t) & \to & V(r'; (t')^s) \end{array}$$

For concreteness we can take some examples of values categories $V$. At this moment we don't need of large sets. All sets will be small sets, i. e. the members of some *universe*.

The first example is a category of small sets $Set$. We take the *Carte tensor product* $X \times Y$ of two sets $X$ and $Y$ as an example of *biproduct*. It will be called a *Carte biproduct*. For some *assistant set* $T$ it provides a functor $X \times T$ within the category of sets $Set$. This functor has the *coadjoint functor* called *exponential functor* $Y^T$. It is defined with the *functional space* $Y^T = Set(T; Y)$ compounded of all mappings $\phi : T \to Y$ as points. The arrows will be defined by *changing of target space*. For a mapping $f : X \to Y$ exponential functor will appoint the mapping $f^T : X^T \to Y^T$ defined by *composition of mappings*. For a mapping $\phi : T \to X$ it appoints the composed mapping $\phi \circ f : T \to Y$.

For such joint pair of functors the *unit transform* is defined by collection of *sections mappings*

$$\lambda_X : X \to (X \times T)^T$$

which for a point $x \in X$ appoints the *section* $\lambda_X(x) : T \to X \times T$. It is a mapping which for a point $s \in T$ appoints the *couple* $\langle x, s \rangle \in X \times T$, we shall write

$$\lambda_X(x)_t = \langle x, t \rangle \in X \times T \ .$$

The *counit transform* is defined by collection of *evaluation mappings*

$$\mathbf{ev}_Y : Y^T \times T \to Y$$

which for a couple $\langle \phi, s \rangle$ of mapping $\phi \in Y^T$ and point $s \in T$ appoints as *value* the point $\phi(s) \in Y$.

For such *joint pair of functors* we have bijective *name appointment*

$$Set(X \times T; Y) =\to Set(X; Y^T) \ .$$

It for a mapping over Carte product $f : X \times T \to Y$, which will be called a *form*, appoints another mapping to the functional space $g : X \to Y^T$. It will be called a *name form*. Its values will be called *partial forms*. .

Another example will be the category of ordered spaces $Mn$. Its objects will be *ordered spaces* $X$ and arrows will be *monotone mappings* between ordered spaces $f : X \Longrightarrow Y$.

An ordered space can be considered as a *convergence space*. In this space for *onepoint sequence* $x \in X$ we appoint a *limit set* compounded of all bigger points

$$\mathtt{lim}\, x = \{y : x \leq y\} \subset X \ .$$

The *monotone mappings* $f : X \to Y$ then coincide with continuous mappings defined by inclusion of limit's image

$$f(\mathtt{lim}\, x) \subset \mathtt{lim}\, f(x) \ .$$

Such convention allows us to use the same *nomenclature* for different areas of mathematics.

We take the *Carte tensor product* of ordered spaces $X \times Y$ defined as Carte tensor product of sets $X \times Y$ with *initial order relation* for the projections $p_1 : X \times Y \to X$ and $p_2 : X \times Y \to Y$. It will be again Carte tensor product of corresponding convergence space.

For some *assistant ordered space* $T$ defined functor $X \times T$ has a coadjoint exponential functor defined with a new functional space $Y^T = Mn(T; Y)$ compounded of all monotone mappings $\phi : T \Longrightarrow Y$. The functional space in the category of sets $Set$ for a small set $T$ coincides with infinite Carte tensor product

$$Set(T; Y) = \prod_{t \in T} Y$$

i. e. each mapping $\phi : T \to Y$ can be considered as collection of values

$$\langle \phi(t) \uparrow t \in T \rangle \in \prod_{t \in T} Y \ .$$

For ordered space $Y$ we get an order relation in such product defined again as the order relation *initial for these projections* $p_t : \prod_{t \in T} Y \to Y$. The monotone mappings compound a *partial space* in such ordered space

$$Mn(T; Y) \subset \prod_{t \in T} Y \ .$$

Such situation can be red also for the convergence spaces.

The unit transform will be defined by the collection of the same sections mappings $\lambda_X : X \to (X \times T)^T$. We need only to check that it will be monotone or continuous for corresponding convergence spaces

$$\lambda_X : X \to (X \times T)^T \subset \prod_{t \in T}(X \times T) \ .$$

First we shall check that every section $\lambda_X(x) : T \to X \times T$ is continuous. We need to check the continuity only for each projection.

The composition with the first projection is get as composition with the unique continuous mapping to the final *onepoint space*

$$\lambda_X(x) \circ p_1 : T \to X \times T \to \{x\} \subset X \ ,$$

so it again coincides with the unique arrow to onepoint set.

The composition with the second projection coincides with the identic mapping over assistant space

$$\lambda_X(x) \circ p_2 : T \to X \times T \to T = \texttt{Id}_T : T \to T \ .$$

So both projections provide continuous mappings, therefore sections are continuous.

The continuity of sections mapping $\lambda_X : X \to (X \times T)^T$ must be checked again for every projection $p_t : (X \times T)^T \to X \times T$ and further for two projections over Carte tensor product $X \times T$. We get the identity mapping $\texttt{Id}_X : X \to X$ for the first projection

$$\lambda_X \circ p_t \circ p_1 = \texttt{Id}_X : X \to (X \times T)^T \to X \times T \to X \ ,$$

and we get unique arrow to the onepoint set $! : X \to \{t\}$ for the second projection

$$\lambda_X \circ p_t \circ p_2 = ! : X \to (X \times T)^T \to X \times T \to \{t\} \ .$$

So we get the continuous projection $\lambda_X \circ p_t : X \to X \times T$, consequently also the sections mapping $\lambda_X : X \to (X \times T)^T$ itself is continuous.

The counit transform in the category *Set* was defined by collection of evaluation mappings

$$\mathbf{ev}_Y : (\prod_{t \in T} Y) \times T \to Y \ .$$

Now we want to restrict such mapping over smaller set of continuous mappings

$$\mathbf{ev}_Y : Y^T \times T \to Y \ .$$

In functional space $Y^T$ we have defined an initial convergence, which usually is named as a weak convergence. There are only special cases when we can get continuous evaluation for weak convergence in functional space. But for the ordered spaces evaluation mapping becomes monotone. Let we have inequality for monotone mappings
$$\phi(t) \leq \psi(t) \uparrow t \in T \ .$$
For pair of arguments $s \leq t$ we can check inequality between acquired values

$$\phi(s) \leq \psi(t) \ .$$

Indeed. For transitive relation we get

$$\phi(s) \leq \phi(t) \leq \psi(t) \ .$$

For taken joint pair of functors we get again the bijective name appointment
$$Mn(X \times T; Y) =\to Mn(X; Y^T)$$
which for arbitrary form as monotone mapping over the product of ordered spaces $f : X \times T \to Y$ appoints the monotone mapping of *partial forms* over the source space $X$.

The ordered space is an instance of a set $X$ with binary relation $R \subset X \times X$, or more generally an instance of a graph $X$. We can understand the inequality $a \geq b$ as an edge $f : a \to b$ between the vertexes. Otherwise sometimes is useful with the sign of inequality to denote the arrows from some another category. Such simplified notation can be applied when we want to work with different sorts of arrows as we see for the **2**-*categories*. Usually somebody prefer to speak about *maps* and *cells*, but that language can be applied only for special class of our considered *superfunctors*.

The graph $X$ is understood as a set of *edges* between the *vertexes*. The set of vertexes will be denoted $X_0$. We shall work with *directed graphs*. Every edge has its direction. A *graphs transport* between two graphs $f : X \Longrightarrow Y$ is defined by the *vertexes mapping* $f_0 : X_0 \to Y_0$ and the family of *edges mappings* $f_{x,y} : X(x;y) \to Y(f_0(x); f_0(y))$. They are composed with usual rules
$$(f \circ g)_0 = f_0 \circ g_0 \ , \quad (f \circ g)_{x,y} = f_{x,y} \circ g_{f_0(x), f_0(y)} \ ,$$
and we have an *identity transport* defined with identity mappings

$$\text{Id}_{X_0} : X_0 \to X_0 \ , \quad \text{Id}_{X(x,y)} : X(x;y) \to X(x;y) \ .$$

Identity transports are neutral for composition of graphs transports. They will be denoted as identity graphs transport $\text{Id}_X : X \Longrightarrow X$.

The graphs transports will be arrows in the category of graphs $Grph$. We define a *Carte biproduct* $X \times Y$. It has the Carte tensor product of vertexes sets $X_0 \times Y_0$ and the Carte tensor product of edges set

$$X \times Y(\langle x, y \rangle; \langle x', y' \rangle) = X(x; x') \times Y(y; y') \uparrow \langle x, y \rangle \in X_0 \times Y_0, \langle x', y' \rangle \in X_0 \times Y_0 \ .$$

For chosen *assistant graph* $T$ we shall define a joint pair of superfunctors. The adjoint superfunctor is defined with the Carte biproduct $X \times T$. For a graph $X$ it appoints the graphs biproduct with the assistant graph $X \times T$ and for a graphs transport $f : X \Longrightarrow X'$ it appoints the graphs transport between Carte biproducts defined by the biproduct of taken graphs transport with the identity transport of the assistant graph

$$f \times \text{Id}_T : X \times T \Longrightarrow X' \times T \ .$$

The coadjoint superfunctor $Y^T$ is defined with the graph of all transforms between graphs transports over the assistant graph $\phi : T \Longrightarrow Y$. These graphs transports become vertexes, and edges will be transforms $\alpha : \phi \to \psi$ between them defined with some collection of edges

$$\alpha_t \in Y(\phi_0(t); \psi_0(t)) \ .$$

For an graphs transport $g : Y \Longrightarrow Y'$ it appoints the new graphs transport defined by changing of target space with taken graphs transport. For a graphs transport $\phi : T \Longrightarrow Y$ it appoints the new graphs transport $\phi \circ g : T \Longrightarrow Y'$ and for a transform between graphs transports $\alpha : \phi \to \psi$ it appoints the new transform $\alpha * T : \phi \circ f \to \psi \circ f$ defined with the new collection of arrows

$$g_{\phi_0(t),\psi_0(t)}(\alpha_t) \in Y'(g_0(\phi_0(t)); g_0(\psi_0(t))) \ .$$

We can see that another changing of target space with another graphs transport $h : Y' \Longrightarrow Y''$ provides the same result as changing of target space with composition of graphs transports $g \circ h : Y \Longrightarrow Y''$. Also the changing of target space with identity transport $\text{Id}_Y : Y \Longrightarrow Y$ provides the identity graphs transport again. So we have got indeed the superfunctor $Y^T$.

We get a *partial unit transform* for such joint pair of superfunctors defined with the collection of sections transports $\lambda_X : X \Longrightarrow (X \times T)^T$ for the assistant graph $T$ and the source graph $X$ with chosen unit arrows $1_s \in T(s; s)$ and $1_a \in X(a; a)$. For vertexes we take ordinary sections mappings

$$\lambda_{X_0} : X_0 \to X_0^{T_0} \ .$$

For the edge $\alpha \in X(x, y)$ we appoint the transform between two sections $\lambda_X(x) : T \Longrightarrow X \times T$ and $\lambda_X(y) : T \Longrightarrow X \times T$. These sections are graphs transports with the mappings for vertexes $\lambda_{X_0}(x) : T_0 \to X_0 \times T_0$ defined

$$\lambda_{X_0}(x)_t = \langle x, t \rangle \in X_0 \times T_0$$

and the edges mappings are provided

$$\lambda_X(\alpha) = \langle \alpha, 1_s \rangle \in X(x; y) \times T(s; t) \uparrow \alpha \in X(x; y) \ .$$

The unit transform allows us to define a *name appointment* for the forms over the Carte biproduct

$$Grph(X \times T; Y) \to Grph(X; Y^T) \ .$$

For a graphs transport $f : X \times Y \Longrightarrow Y$ it appoints the composition of sections mapping with the graphs transport provided by changing of target space with taken graphs transport

$$g = \lambda_X \circ f^T : X \Longrightarrow Y^T .$$

For a vertex $x \in X_0$ it appoints the section which for a vertex $t \in T_0$ appoints the couple $\langle x, t \rangle \in X_0 \times T_0$, and further the changing of target space will appoint the graphs transport $T \Longrightarrow Y$ which for a vertex $t \in T$ appoints the vertex in the target space $f(x, t) \in Y_0$ and for an edge $u \in T(s; t)$ it will appoint the image of couple $\langle 1_x, u \rangle \in X \times T(\langle x, s \rangle; \langle x, t \rangle)$ by taken graphs transport

$$f(1_x, u) \in Y(f_0(x, s); f_0(x, t)) .$$

For an edge $\alpha \in X(a; b)$ it appoints the transform between two sections defined with collection of edges $\langle \alpha, 1_s \rangle \in X \times T(\langle a, s \rangle; \langle b, s \rangle)$ and further the transform between graphs transports to the target graph $Y$ defined with collection of edges $f(\alpha, 1_s) \in Y(f_0(a, s); f_0(b, s))$.

The partial counit transform is defined with collection of evaluation transports $\mathbf{ev}_Y : Y^T \times T \Longrightarrow Y$ in the target graph $Y$ with mappings of contact arrows

$$Y_{c,d,e} : Y(c; d) \times Y(d; e) \to Y(c; d) \circ Y(d; e) \subset Y(c; e) .$$

Evaluation transport can be defined in two ways.

For the couple $\langle \alpha, u \rangle$ of a transform $\alpha : \phi \to \psi$ between to graphs transports $\phi, \psi : T \Longrightarrow Y$ and an edge between two vertexes in source graph $u \in T(s; t)$ we can define the *precontact arrow*

$$\mathbf{ev}_Y^+(\langle \alpha, u \rangle) = \alpha_s \circ \psi_{s,t}(u) ,$$

and another *postcontact arrow*

$$\mathbf{ev}_Y^-(\langle \alpha, u \rangle) = \phi_{s,t}(u) \circ \alpha_t .$$

The counit transport allows us to define the *realization appointment* for partial forms transport over the source graph $X$

$$Grph(X; Y^T) \to Grph(X \times T; Y) .$$

With precontact arrows we get the *prerealization*, and with postcontact arrows we define the *postrealization*.

For a graphs transport $g : X \Longrightarrow Y^T$ it appoints the biproduct of taken graphs transport with the identity transport of the assistant graph $g \times \mathtt{Id}_T : X \times T \Longrightarrow Y^T \times T$ composed with the evaluation transport $\mathbf{ev}_Y : Y^T \times T \Longrightarrow Y$

$$f = (g \times \mathtt{Id}_Y) \circ \mathbf{ev}_Y : X \times T \Longrightarrow Y .$$

For a couple $\langle x, t \rangle$ of $x \in X_0$ and $t \in T_0$ it appoints the couple $\langle g_0(x), t \rangle$ of the graphs transport $g_0(x) \in (Y^T)_0$ and the same vertex in the assistant graph $t \in T_0$. Further we get the value of this graphs transport over taken vertex $g_0(x)_t \in Y_0$. For a couple $\langle \alpha, u \rangle$ of an edge in the source graph $\alpha \in X(a;b)$ and an edge in the assistant graph $u \in T(s;t)$ we get the couple $\langle g(\alpha), u \rangle$ of transform $g(\alpha) : g_0(a) \to g_0(b)$ between graphs transports $g_0(a), g_0(b) : T \Longrightarrow Y$ defined with collection of edges $g(\alpha)_s \in Y(g_0(a)_s; g_0(b)_s)$ and the same edge $u \in T(s;t)$. Further the evaluation transport appoints the precontact arrow

$$g(\alpha)_s \circ g_0(b)_u \in Y(g_0(a)_s; g_0(b)_t)$$

or the postcontact arrow

$$g_0(a)_u \circ g(\alpha)_t \in Y(g_0(a)_s; g_0(b)_t) \ .$$

Earlier in **G. Valiukevičius 2022** we have checked when realization appointment *reverses* name of graphs transport $f \in Grph(X \times T; Y)$. For the preevaluation we have used the *predecomposable property*

$$f(\alpha, 1_s) \circ f(1_b, u) = f(\alpha, u) \ ,$$

and for the postevaluation we have used the *postdecomposable property*

$$f(1_a, u) \circ f(\alpha, 1_t) = f(\alpha, u) \ .$$

In these cases we get that the realization appointment is *surjective* and the name appointment is *injective*.

Also we have checked when the name appointment reverses the realization of graphs transports $g : X \Longrightarrow Y^T$. For the preevaluation we have used the *preneutral properties*

$$g(1_a)_s \circ g_0(a)_u = g_0(a)_u) \ , \quad g(\alpha)_s \circ g_0(b)_{1_s} = g(\alpha)_s \ ,$$

and for the postevaluation we have used the *postneutral properties*

$$g_0(a)_u \circ g(1_a)_t = g_0(a)_u \ , \quad g_0(a)_{1_s} \circ g(\alpha)_s = g(\alpha)_s \ .$$

In these cases we get that the name appointment is surjective and the realization appointment is injective.

A transform $\alpha : \phi \to \psi$ is called *natural* if we get commuting quadrats from equality of precontact and postcontact arrows

$$\alpha_s \circ \psi_0(u) = \phi_0(u) \circ \alpha_t \ .$$

$$\begin{array}{ccc} \phi_0(s) & \xrightarrow{\phi_{s,t}(u)} & \phi_0(t) \\ \alpha_s \downarrow & \searrow & \downarrow \alpha_t \\ \psi_0(s) & \xrightarrow[\psi_{s,t}(u)]{} & \psi_0(t) \end{array}$$

For two graphs transports $\phi, \psi : T \Longrightarrow Y$ we shall denote the set of natural transforms $\alpha : \phi \to \psi$ with a new sign $Y^{(T)}(\phi; \psi)$. These will be edges sets for a new exponential graph $Y^{(T)}$. At this moment we are not interested can such new exponential graphs to define a superfunctor.

Having predecomposable and postdecomposable properties we get injective name appointment to the set $Grph(X; Y^T)$ of all transports to the exponential graph. We shall say that taken graphs transport $f : X \times T \Longrightarrow Y$ is *decomposable*. We have seen that such graphs transport can be *corestrained* in the new exponential graph of natural transforms $Y^{(T)}$.

Both preevaluation and postevaluation transports coincide over the exponential graph of natural transforms. This is easy to see from equality

$$f(\alpha, 1_s) \circ f(1_b, u) = f(1_a, u) \circ f(\alpha, 1_s) \ .$$

In this case both realization mappings also coincide and their image contains all decomposable graphs transports.

Having graphs transport $g : X \Longrightarrow Y^T$ with preneutral and postneutral properties we get that both prerealization and postrealization mappings reverses the name of taken graphs transport. We shall say that taken graphs transport $g : X \Longrightarrow Y^T$ has *neutral properties*. Such properties can be visualized by *splitting diagrams*

$$\begin{array}{ccccccc}
 & g_0(a)_{1_s} & & & & g_0(a)_u & \\
g_0(a)_s & \longrightarrow & g_0(a)_s & g_0(a) & \longrightarrow & & g_0(a)_t \\
g(\alpha)_s \downarrow & \searrow & \downarrow g(\alpha)_s & g(1_a)_s \downarrow & \searrow & & \downarrow g(1_a)_t \\
g_0(b)_{1_s} & \longrightarrow & g_0(b)_s & g_0(a)_s & \longrightarrow & & g_0(a)_t \\
 & g_0(b)_{1_s} & & & & g_0(a)_u & \\
\end{array}$$

Such diagrams were applied by **P. Freyd, A. Scedrov 1990** §1.28 to define $g_0(a)_{1_s}$ and $g(1_a)_s$ as splitting idempotents.

In this case the name mapping is surjective and both realization are injective. For the decomposable graphs transport $f : X \times T \Longrightarrow Y$ having name $g : X \Longrightarrow Y^T$ with neutral properties we get bijective name appointment

$$Grph(X \times T; Y) =\to Grph(X; Y^{(T)}) \ .$$

We shall restrict the general category of graphs as we want to get superfunctor $Y^{(T)}$ for the exponential graphs of natural transforms. First we shall define the *original graphs*. In **J. Freyd, A. Scedrov 1990** §1.263 we can see a name *pointed sets* for the sets with chosen origin point. However in 1.4(11)4 is defined the point in arbitrary object from some category, so we don't want to use the name of pointed set for objects in arbitrary category.

In original graph $X$ we choose unit arrows $1_a \in X(a; a)$. The *original graphs transports* $f : X \Longrightarrow Y$ must maintain chosen unit arrows

$$f(1_a) = 1_{f_0(a)} \ .$$

They will be called also as *crude functors*. Such category of original graphs we shall denote $OGrph$. The letter $O$ must denote the graphs with origin point.

We can take yet a smaller category compounded of all original graphs transports between original graphs with contact arrows mappings

$$\alpha \circ \beta \in X(a; c) \uparrow \alpha \in X(a; b), \beta \in X(b; c) \ .$$

They define a graphs transport over *pullbacked product* $X \times_p X \subset X \times X$.

$$\circ : X \times_p X \Longrightarrow X .$$

In original graph $X$ we demand a neutral properties for chosen unit arrows

$$1_a \circ \alpha = \alpha = \alpha \circ 1_b \uparrow \alpha \in X(a;b) .$$

For the category $OGrph^=$ we take original graphs transports with additional property to maintain the contact arrows mappings

$$f(\alpha \circ \beta) = f(\alpha) \circ f(\beta) \uparrow \alpha \in X(a;b), \beta \in X(b;c) .$$

Such graphs transports can be called as *usual functors*, with more concrete language we can call them as *natural functors*. The choosing of unit arrows or contact arrows mappings can be announced as some algebraic operations and special graphs transforms must maintain such operations.

We have the Carte biproduct $X \times Y$ of two original graphs defined with initial structure of operations. For a couple of vertexes $\langle x, y \rangle \in X_0 \times Y_0$ we appoint the couple of unit arrows

$$\langle 1_x, 1_y \rangle \in X(x;x) \times Y(y;y) .$$

For two couples of edges

$$\langle \alpha, \alpha' \rangle \in X(a;b) \times Y(a';b'), \langle \beta', \beta' \rangle \in X(b,c) \times Y(c';d')$$

we can define a contact arrow

$$\langle \alpha, \alpha' \rangle \circ \langle \beta, \beta' \rangle = \langle \alpha \circ \beta, \alpha' \circ \beta' \rangle \in X(a;c) \times Y(a';c') .$$

So for arbitrary assistant original graph $T$ we define superfunctor $X \times T$ in the category of original graphs $OGrph$ or in the smaller category $OGrph^=$ of *natural functors* between original graphs.

For an original graph $X$ it appoints the Carte biproduct $X \times T$ and for an original graphs transports $f : X \Longrightarrow X'$ it appoints the Carte biproduct with the identity transport of the assistant graph $\text{Id}_T : T \Longrightarrow T$

$$f \times \text{Id}_T : X \times T \Longrightarrow X' \times T .$$

Another superfunctor will be defined with the graph of all transforms $Y^T$. The unit arrows and contact arrows mappings will be taken as initial structure in the product $\prod_{t \in T_0} Y$. The unit transform $1_\phi : \phi \to \phi$ is defined with the collection of unit arrows

$$1_{\phi_0(s)} \in Y(\phi_0(s); \phi_0(s)) \uparrow s \in T_0 .$$

The contact arrows mapping also is defined by collection of contact edges

$$(\alpha \circ \beta)_s = \alpha_s \circ \beta_s \uparrow s \in T_0 .$$

The unit arrow $1_\phi : \phi \to \phi$ will be a natural transform, but the composition of natural transforms doesn't need to remain a natural transform. However for the natural original graphs functor over the Carte biproduct $f : X \times T \implies Y$ the name form $g : X \implies Y^T$ will be corestricted in the exponential space a natural transform.

Finally we restrict ourselves with the *categories*. The original graph $X$ becomes a category if its contact arrows are get with associative composition of arrows, i. e. we have the identity provided with associativity isomorphism

$$\alpha \circ (\beta \circ \gamma) = (\alpha \circ \beta) \circ \gamma \uparrow \alpha \in X(a;b), \beta \in X(b;c), \gamma \in X(c;d) \ .$$

Vertexes of category will be called *objects* $x \in X_0$ and edges will be called *arrows* $\alpha \in X(a;b)$. The functors between categories $f : X \implies Y$ will compound the category $Cat$.

For the category of small categories $Cat$ we can define an exponential superfunctor $Y^{(T)}$ of natural transforms. First we need to check that the composition of two natural transforms remains a natural transform.

Let we have three functors $\phi, \psi, \xi : T \implies Y$ and two natural transforms $\alpha : \phi \to \psi$ and $\beta : \psi \to \xi$. Having *commuting quadrats* for arbitrary edge in assistant category $u \in T(s;t)$

$$\phi(u) \circ \alpha_t = \alpha_s \circ \psi(u) \ ,$$

$$\psi(u) \circ \beta_t = \beta_s \circ \xi(u) \ ,$$

we get a new commuting quadrat

$$\phi(u) \circ (\alpha_t \circ \beta_t) = (\alpha_s \circ \beta_s) \circ \xi(u) \ .$$

$$\begin{array}{ccccc}
\phi_0(s) & \xrightarrow{\alpha_s} & \psi_0(s) & \xrightarrow{\beta_s} & \xi_0(s) \\
\phi(u) \downarrow & & \beta(u) \downarrow & & \downarrow \xi(u) \\
\phi_0(t) & \xrightarrow[\alpha_y]{} & \psi_0(t) & \xrightarrow[\beta_t]{} & \xi_0(t)
\end{array}$$

The proving uses three ways from a vertex $\phi_0(s)$ to the vertex $\psi_0(t)$

$$(\alpha_s \circ \beta_s) \circ \xi(u) = \alpha_s \circ (\beta_s \circ \xi(u)) = \alpha_s \circ (\psi(u) \circ \beta_t) =$$

$$(\alpha_s \circ \psi(u)) \circ \beta_t = (\phi(u) \circ \alpha_t) \circ \beta_t = \phi(u) \circ (\alpha_t \circ \beta_t) \ .$$

We also can prove that the space of natural transforms $Y^{(T)}$ defines a superfunctor for the graphs transports $f : X \implies Y$ maintaining contact arrows mappings

$$f(u \circ v) = f(u) \circ f(v) \uparrow u \in X(a;b), v \in X(b;c) \ .$$

The changing of target space $f^T : Y^T \implies (Y')^T$ will maintain the commuting diagrams for natural transform $\alpha : \phi \to \psi$, i. e. we get commuting diagram again in the new target space $Y'$

$$\begin{array}{ccc}
f_0(\phi_0(s)) & \xrightarrow{f(\alpha_s)} & f_0(\psi_0(s)) \\
f(\phi(u)) \downarrow & & \downarrow f(\psi(u)) \\
f_0(\phi_0(t)) & \xrightarrow[f(\alpha_t)]{} & f_0(\psi_0(t))
\end{array}$$

For the original graphs transport $f : X \Longrightarrow Y$ we get the changing of target space again maintaining unit transform defined by by collection of chosen unit arrows $1_{\phi_0(x)} \in X(\phi_0(x); \phi_0(x))$. Also for the natural functor between original graphs we get again the changing of target space maintaining the composition of transforms

$$f(\alpha \circ \beta) = \prod_{s \in T_0} f(\alpha_s \circ \beta_s) = \prod_{s \in T_0} (f(\alpha_s) \circ f(\beta_s)) = f(\alpha) \circ f(\beta) \ .$$

We shall prefer for some functor $F : X \Longrightarrow Y$ to use the usual notation for points
$$F_0(x) \in Y_0 \uparrow x \in X_0$$
and the exponential notation for arrows
$$f^F \in Y(F_0(x); F_0(y)) \uparrow f \in X(x; y) \ .$$

We extend the proving of Yoneda lemma for more general original graphs with contact arrows mappings. Let we have a graphs transport to the sets category $F : X^{op} \Longrightarrow Set$. It will be called a *character* over taken original graph $X$. The representable character is defined by some point $x \in X_0$

$$X(-; x) : X^{op} \Longrightarrow Set \ .$$

It can be understood as a value for the Yoneda mapping

$$(Y_X)_0 : X_0 \to (Set^{X^{op}})_0 \ .$$

This mapping is extended as original graphs transport

$$Y_X : X \Longrightarrow Set^{X^{op}}$$

and this extension will be *fullfaithful imbedding*, i. e. we have bijection between the edges set $X(a, b)$ and the set of transforms between appointed representable characters $Set^{X^{op}}(X(-; a); X(-; b))$. Later we shall describe such transforms as natural or continuous ones.

The choosing of unit arrow $i_x \in X(x; x)$ is a universal arrow from onepoint set $* \in Set_0$ to the representable functor $X(-; x) : X^{op} \Longrightarrow Set$, clr. **S. MacLane 1998** §3.1 Universal arrows. The arbitrary point in the edges set $F : * \to X(a; x)$ is get from universal arrow $i_x : * \to X(x; x)$ with source changing $f^* : X(x; x) \to X(a; x)$

$$\begin{array}{ccccc} * & \xrightarrow{i_x} & X(x;x) & & x \\ & \searrow & \downarrow f^* & & \uparrow f \\ & & X(a;x) & & a \end{array}$$

For arbitrary graphs transport $F : X^{op} \Longrightarrow Set$ we shall get a bijection between the set of points $m \in F_0(x)$ and the set of transforms between representable character $X(-; x)$ and taken graphs transport $F$.

**Proposition 1.** *Over an original graph $X$ with contact arrows mappings*

$$X(a;b) \otimes X(b;c) \to X(a;b) \circ X(b;c) \subset X(a;c)$$

*for arbitrary graphs transport $F : X^{op} \Rightarrow \mathrm{Set}$ we have a bijection between the value set $F_0(x)$ and the set of transforms between representable character $X(-;x)$ and taken graphs transport $F$.*

**Proof:** A transform between two graphs transports $\alpha : X(-;x) \to F$ is defined by collection of mappings between value sets

$$\alpha_a : X(a;x) \to F_0(a) \ .$$

We take the image of unit arrow $i_x \in X(x;x)$

$$m = \alpha_x(i_x) \in F_0(x)$$

as a point correspondent to taken transform $\alpha : X(-;x) \to F$.

Otherwise for arbitrary point $m \in F_0(x)$ we get a mapping

$$\alpha_a : X(a;x) \to F_0(a) \ .$$

For an edge $f \in X(a;x)$ the mapping $f^F : F_0(x) \to F_0(x)$ provides the value of taken point $f^F(m) \in F_0(a)$.

For a unit arrow $i_x \in X(x;x)$ we get the identity mapping as source changing $i_x^* := \mathrm{Id}_{X(x;x)} : X(x;x) \to X(x;x)$. So we get the same point $m = i_x^F(m)$ correspondent for with this point $m \in F_0(x)$ constructed transform $\alpha : X(-;x) \to F$. Otherwise for some transform $\alpha : X(-;x) \to F$ appointed point $\alpha_x(i_x) \in F_0(x)$ generates the same transform $\alpha : X(-;x) \to F$.

$\square$

For an original graph $X$ with associative contact arrows mappings the representable characters are natural functors, and all transforms between two characters are natural ones. For non associative contact arrows mappings a situation becomes more interesting.

In general the source changing isn't natural functor. But for the unit arrow $i_x \in X(x;x)$ we get equality

$$h \circ (f \circ i_x) = (h \circ f) \circ i_x \ .$$

We shall say that unit arrow $i_x \in X(x;x)$ is a *social point* for source changing functor

$$h^*(f^*(i_x)) = (h \circ f)^*(i_x) \ .$$

Every natural transform $\alpha : X(-;x) \Longrightarrow F$ is defined by collection of mappings $\alpha_a : X(a;x) \to F(a)$. For a natural transform we have commuting diagrams for any edge $h \in X(b;a)$

$$\begin{array}{ccc} X(a;x) & \xrightarrow{\alpha_a} & F_0(a) \\ h^* \downarrow & & \downarrow h^F \\ X(b;x) & \xrightarrow{\alpha_b} & K(b) \end{array}$$

Therefore we get again equality of social point $m \in F_0(x)$

$$h^F(f^F(m)) = (h \circ f)^F(m) \ .$$

**Proposition 2.** *For an original graph $X$ with contact arrows mappings*

$$\circ : X(a;b) \times X(b;c) \to X(a;c)$$

*the choosing of social point $m \in F_0(x)$ provides a natural transform*

$$\alpha : X(-;x) \to K$$

*between representable character $X(-;x)$ and taken character $F$.*

**Proof:** The natural transform $\alpha : X(-;x) \to F$ is defined by collection of mappings between appointed sets

$$\alpha_a : X(a;x) \to K(a)$$

with commuting diagrams for every edge $h \in X(b;a)$

$$\begin{array}{ccc} X(a;x) & \xrightarrow{\alpha_a} & F_0(a) \\ h^* \downarrow & & \downarrow h^F \\ X(b;x) & \xrightarrow[\alpha_b]{} & F_0(b) \end{array}$$

For the unit arrow $1_x \in X(x;x)$ we have a social point $\alpha_x(1_x) \in F_0(x)$. Otherwise each social point $m = \alpha_x(1_x) \in F_0(x)$ can be get from some uniquely defined natural transform. For generated transform $\alpha : X(-;x) \to F$ we need to check commuting diagram for arbitrary edge $h \in X(b;a)$

$$\begin{array}{ccc} X(a;x) & \xrightarrow{\alpha_a} & F_0(a) \\ h^* \downarrow & & \downarrow h^F \\ X(b;x) & \xrightarrow[\alpha_b]{} & F_0(b) \end{array}$$

The equality of sets mappings suffice to check for each point in source set. Therefore we can apply the larger commuting diagram

$$\begin{array}{c} m \in F_0(x) \\ F(f) \downarrow \quad \searrow \\ \nearrow \quad F_0(a) \xrightarrow[F(h)]{} F_0(b) \\ i_x \in X(x;x) \\ f^* \downarrow \quad \searrow \\ \nearrow \quad \nearrow \\ X(a;x) \xrightarrow[h^*]{} X(b;x) \end{array}$$

□

We have got that the set of natural transforms $\alpha : X(-;x) \to F$ can be identified with the set of social points in values set $F(x)$.

Taking another representable character $F = X(-;y)$ we get bijection between natural transforms $\alpha : X(-;x) \to X(-;y)$ and social edges $\theta \in X(x;y)$. For such natural transform we appoint an image $\alpha_x(1_x) \in X(x;y)$. Conversely for each edge $\theta \in X(x;y)$ we appoint the transform defined with collection of mappings $\alpha_a(f) = f \circ \theta$.

$$\begin{array}{ccc} 1_x \in X(x;x) & \xrightarrow{\alpha_x} & \theta \in X(x;y) \\ f^* \downarrow & & \downarrow f^* \\ f \in X(a;x) & \xrightarrow{\alpha_a} & f \circ \theta \in X(a;y) \end{array}$$

By the former proposition such transform is natural.

In general For the contact arrow $\theta \circ \kappa \in X(x;z) \uparrow \kappa \in X(y;z)$ we don't get vertical composition of correspondent natural transforms. But for composition of natural transforms we get composed contact arrows mappings. We need to check the composition only over unit arrows.

Let we have composable natural transforms $\alpha : X(-;x) \to X(-;y)$ and $\beta : X(-;y) \to X(-;z)$ with corespondent edges $\alpha_x(1_x) = \theta \in X(x;y)$ and $\theta_y(1_y) = \kappa \in X(y;z)$. Then for the composition $\alpha \circ \beta : X(-;x) \to X(-;z)$ we get contact arrow $\beta_y(\alpha_x(1_x)) = \theta \circ \kappa$

$$\begin{array}{ccccc} & & 1_y \in X(y;y) & \xrightarrow{\beta} & \kappa \in X(y;z) \\ & & \theta^* \downarrow & & \theta^* \\ 1_x \in X(x;x) & \xrightarrow{\alpha_x} & X(x;y) & \xrightarrow{\beta_x} & \theta \circ \kappa \in X(x;z) \end{array}$$

As a consequence we deduce important property.

For the two vertexes $x \in X_0$ and $y \in X_0$ having isomorphic natural transform between representable characters

$$\alpha : X(-;x) =\to X(-;y)$$

defined with edge $\alpha_x(1_x) = \theta \in X(x;y)$ we get the inverse natural transform $\alpha^{-1} : X(-;y) \to X(-;x)$ defined with the inverse edge $\alpha^{-1}(1_y) = \theta^{-1} \in X(y;x)$, i. e. we have equalities of contact arrows

$$\theta \circ \theta^{-1} = 1_x , \quad \theta^{-1} \circ \theta = 1_y .$$

For an original graph with associative contact arrows mappings every representable character is natural functor, as for composable edges $h \in X(c;b)$ and $g \in X(b;a)$ the representable character $X(-;x) : X^{op} \Longrightarrow V$ provides composable source changing

$$X(a;x) \xrightarrow{g^*} X(b;x) \xrightarrow{h^*} X(c;x) ,$$

with an equality

$$h \circ (g \circ f) = (h \circ g) \circ f \uparrow f \in X(a;x) ,$$

i. e. $h^*(g^*(f)) = (h \circ g)^*(f)$.

Also for associative contact arrow mappings we get composable natural transforms defined by arbitrary composable edges.

$$\begin{array}{ccccccc}
& & 1_y \in X(y;y) & \xrightarrow{\beta_y} & \kappa \in X(y;z) \\
& & \theta^* \downarrow & & \theta^* \\
1_x \in X(x;x) & \xrightarrow{\alpha_x} & X(x;y) & \xrightarrow{\beta_x} & \theta \circ \kappa \in X(x;z) \\
f^* \downarrow & & f^* \downarrow & & \downarrow f^* \\
f \in X(a;x) & \xrightarrow[f \circ \theta]{} & f \circ \theta \in X(a;y) & \xrightarrow[f \circ (\theta \circ \kappa)]{} & X(a;z)
\end{array}$$

$$(\alpha \circ \beta)_a(f) = f \circ (\theta \circ \kappa) = (f \circ \theta) \circ \kappa \; .$$

We shall notice that Yoneda lemma can be applied for more large stock of transforms between two characters. For the graph $Y$ having edges sets ordered with some order relation we get order relation for parallel mappings of edges set. Let we have two mappings $u, v : X(a;b) \to Y(\phi_0(a); \psi_0(b))$. Then we can take initial order relation between such mappings, i. e. we have an inequality between two mappings if such inequality is true for each pairs of values. We shall denote such inequality as an arrow

$$u \geq v \iff u \to v \uparrow u, v \in Set(X(a;b); Y(\phi_0(a); \psi_0(b))) \; .$$

We could define the set of *continuous transforms* by diagram

$$\begin{array}{ccc}
F_0(x) & \xrightarrow{\alpha_x} & G_0(x) \\
f^F \downarrow & & \downarrow f^G \\
F_0(y) & \xrightarrow[\alpha_y]{} & G_0(y) \\
& & \nearrow
\end{array}$$

It will be the *deformed quadrats* $(-13)$ for *respecting reform* $\alpha$ in **G. Valiukevičius 2009**. The arrows $f^F$ and $f^G$ will be interpreted as *generalized convergences*.

With another direction of inequality we shall define *cocontinuous transforms*

$$\begin{array}{ccc}
F_0(x) & \xrightarrow{\alpha_x} & G_0(x) \\
F(f) \downarrow & & \downarrow G(f) \\
F_0(y) & \xrightarrow[\alpha_y]{} & G_0(y) \\
& & \swarrow
\end{array}$$

It will correspond the *deformed quadrats* $(+13)$ of *wasting reform* $\alpha$.

The deformed quadrats provides instances of *diagrams*. Their edges may be interpreted as functors between different categories. But later we shall define also some coherent transform of continuous quadrats, and we shall use another sort of edges defined with arrows of natural transforms.

The Yoneda lemma for continuous or cocontinuous transforms needs to define for the graphs transport $F : X^{op} \implies Set$ lower social points or upper social points.

For the *lower social point* $m \in F_0(x)$ we demand continuous mapping $h^F : F_0(a) \to F_0(b)$ for every edge $h \in X^{op}(a;b)$ i. e. we should have deformed triangular diagram

$$\begin{array}{c} m \\ \downarrow \searrow \\ F_0(a) \xrightarrow[h^F]{} F_0(b) \nearrow \end{array}$$

For the upper social point $m \in F_0(x)$ we demand that a mapping $h^F : F_0(a) \to F_0(b)$ would be cocontinuous, i. e. we would have a triangular diagram with opposite order relation

$$\begin{array}{c} m \\ \downarrow \searrow \\ F_0(a) \xrightarrow[h^F]{} F_0(b) \swarrow \end{array}$$

We can check that the lower social points generate continuous transforms, and the upper social points generate cocontinuous transforms between representable character $X(-;x)$ and some graphs transport $F$.

**Proposition 3.** *The lower social point $m \in F_0(x)$ generates continuous transform $\alpha : X(-;x) \to F$. And vice versa continuous transforms are generated by some lower social point.*

**Proof:** The proving can be checked by chasing diagram

$$\begin{array}{ccccc} & & m \in F_0(x) & & \\ & & F(f) \downarrow & \searrow & \\ & \nearrow & & & \\ i_x \in X(x;x) & & F_0(a) & \xrightarrow[F(h)]{} & F_0(b) \\ f^* \downarrow & \searrow & & & \nearrow \\ & \nearrow & & \nearrow & \\ X(a;x) & \xrightarrow[h^*]{} & X(b;x) & & \end{array}$$ □

We can say that contact arrows mappings are *lower associative*, if we have associative arrow as inequality between two $V$-graphs transports

$$(h \circ f) \circ \theta \to h \circ (f \circ \theta) \uparrow \theta \in X(c;x), g \in X(b;c), f \in X(a;b) .$$

This will provide all arrows $\theta \in X(c;x)$ being lower social for representable functor $X(-;x)$.

We can repeat the proposition about name appointment for graphs with ordered edges sets. Let we have graphs $X$ with ordered sets of edges $X(a;b)$. We demand that for a graphs transport $F : X \Longrightarrow Y$ the mappings between edges sets would be monotone

$$F_{a,b} : X(a;b) \to X(F_0(a); F_0(b)) .$$

In the Carte biproduct $X \times Y$ we take the initial order relation for edges sets
$$X \times Y(\langle a, a' \rangle, \langle b, b' \rangle) = X(a; b) \times Y(a'; b') \ .$$

In the exponential space $Y^T$ also we can take an initial order relation inside the set of transforms
$$Y^T(\phi; \psi) = \prod_{s \in T_0} Y(\phi_0(s); \psi_0(s)) \ .$$

We get again a joint pair of monotone superfunctors within the categories of graphs with ordered sets of edges.

The choosing of unit arrows in the source graph $1_a \in X(a; a)$ and in the assistant graph $1_t \in T(t; t)$ allows us to define the unit transform with the collection of sections transports $\lambda_X : X \Longrightarrow (X \times T)^T$. The edges mappings
$$\lambda_X(\alpha) : \lambda_X(a) \to \lambda_X(b)$$

are defined by the collection of edges
$$\langle u, 1_s \rangle : \langle a, s \rangle \to \langle b, s \rangle \uparrow u \in X(a; b) \ ,$$

$$\lambda_X(u) = \prod_{s \in T_0} \langle u, 1_s \rangle : \prod_{s \in T_0} \langle a, 1_s \rangle \to \prod_{s \in T_0} \langle b, ; 1_s \rangle \ .$$

and obviously maintain the order relation.

The choosing of monotone contact arrows in the target graph $Y$
$$Y_{c,d,e} : Y(c; d) \times Y(d; c) \to Y(c; e)$$

allows us to define a counit transform $\mathbf{ev}_Y : Y^T \times T \to Y$ with collection of evaluation transports. The contact arrows mappings must be monotone for every its argument. We can take the preevaluation or postevaluation transport. At this moment we are chosen the preevaluation transport defined with precontact arrows $\mathbf{ev}_V^+ : Y^T \times T \Longrightarrow Y$. For a transform $\alpha : \phi \to \psi$ and an edge $u \in T(s; t)$ we appoint
$$\mathbf{ev}_Y^+(\alpha, u) = \alpha_s \circ \psi(u) \in Y(\phi_0(s); \psi_0(t)) \ .$$

Such appointment will be monotone by the properties of chosen monotone contact arrows mappings.

The unit transform defined with sections transports $\lambda_X : X \Longrightarrow (X \times T)^T$ provides the monotonous name appointment for forms
$$g : Grph(X \times T; Y) \to Grph(X; Y^T) \ .$$

With an inequality between two graphs transports $f \geq f'$ we get the inequality between appointed name forms $f^\# \geq (f')^\#$ provided with composition of monotone graphs transports. Such are the section transport and changing of target space also are monotone for chosen inequalities in exponential spaces

$$g(f) = f^{\#} : X \xrightarrow{\lambda_X} (X \times T)^T \xrightarrow{f^T} Y^T .$$

A counit transform defined with preevaluation transports $\mathbf{ev}_Y^+ : Y^T \times T \Longrightarrow Y$ also is monotone appointment for name forms

$$f_+ : Grph(X; Y^T) \to Grph(X \times T; Y) .$$

With an inequality between two graphs transports $g \geq g'$ we get the inequality between their realizations $f_+(g) \geq f_+(g')$ provided with composition of monotone graphs transports. Such is the biproduct with identity transport $\mathtt{Id}_T : T \Longrightarrow T$ and then we have monotone evaluation transport

$$f_+(g) = g^{\flat} : X \times T \xrightarrow{g \times 1_T} Y^T \times T \xrightarrow{\mathbf{ev}_Y^+} Y .$$

For graphs transport $f : X \times T \Longrightarrow Y$ we can ask *lower predecomposition* property

$$f(\alpha, u) \geq f(\alpha, 1_s) \circ f(1_b; u) \uparrow \alpha \in X(a; b), u \in T(s; t) .$$

Then we get inequalities for prerealization of name form

$$f \geq f_+(g(f)) .$$

It can be considered as unit inequality for joint pair of monotonous mappings provided with the form's name and the name form's realization

$$Grph(X \times T; Y) \Leftrightarrow Grph(X; Y^T) .$$

For a graphs transport $g : X \Longrightarrow Y^T$ we ask *upper preneutral* properties

$$g(1_a)_s \circ g_0(a)_u \geq g_0(a)_u ,$$

$$g(\alpha)_s \circ g_0(b)_{1_s} \geq g_0(\alpha)_s .$$

Then we get a counit inequalities for the joint pair of name and realization appointments

$$g \leq g(f_+(g)) .$$

The categorical proving is a same as for the graphs without any order relation in the sets of edges, cl. **G. Valiukevičius 2022**.

Now we shall try to study $V$-enriched spaces. At first we define the category $Grph_V$ of $V$-enriched graphs. We shall use a shorter name $V$-*graphs*. In the set of arrows $V(r; s)$ between two points $r, s \in V_0$ we choose some order relation. We must note that the equality also can be understood as an order relation. The arrows between points of values category $r, s \in V_0$ will be called as *mappings*.

A $V$-graph $X$ will be a set of vertexes $X$ with the sets of edges being objects in the values category $X(x;y) \in V_0$.

A $V$-*graphs transport* $f : X \Longrightarrow Y$ will be defined with an appointment of vertexes $f_0 : X_0 \to Y_0$ and collection of mappings for the edges sets

$$f_{x,y} : X(x;y) \to Y(f_0(x); f_0(y)) .$$

Such transports compound the category with usual composition. The identity transports will be *unit arrows* for such composition. For two $V$-graphs transports $f : X \Longrightarrow Y$ and $g : Y \Longrightarrow Z$ we have composition of vertexes appointments
$$(f \circ g)_0(x) = g_0(f_0(x)) \in Z \uparrow x \in X$$

and composition of edges set mappings

$$(f \circ g)_{x,y} = f_{x,y} \circ g_{f_0(x), f_0(y)} .$$

An identity transport $\text{Id}_X : X \Longrightarrow X$ is defined by identity appointment for vertexes
$$\text{Id}_{X_0} : X_0 \to X_0$$
and identity mappings for the edges sets

$$\text{Id}_{X(x,y)} : X(x;y) \to X(x;y) .$$

This category will be denoted $Grph_V$.

We shall say that two $V$-graphs transports are parallel if they have the same vertexes appointments

$$f_0(x) = g_0(s) \uparrow s \in X_0 .$$

The order relation in the sets of edges $V(a;b)$ will define the order relation between *parallel $V$-graphs transports*

$$f, g : X \Longrightarrow Y .$$

The inequality $f \geq g$ is defined by order relation for the edges set mappings

$$f_{x,y} \geq g_{x;y} : X(a;y) \to Y(f_0(x); f_0(y)) = Y(g_0(x); g_0(y)) .$$

More strong ordering could be defined by factorization property

$$g = f \circ h .$$

But it is waiting for later investigations and applications.

We shall define some *superfunctors* over the category $Grph_V$. At first we have the *primary projection*

$$Grph_V \Longrightarrow \mathbf{SET} .$$

For each $V$-graph $X$ it appoints the set of vertexes $X_0$ and for a $V$-graphs transport $f : X \Longrightarrow Y$ it appoints the appointment $f_0 : X_0 \to Y_0$ of vertexes.

We have more interesting superfunctors. For two $V$-graphs $X$ and $Y$ we appoint a biproduct $X \otimes Y$ defined with Carte product $X_0 \times Y_0$ for the sets of vertexes and biproduct from values category $V$ for the sets of edges

$$X \otimes Y(\langle x, y \rangle; \langle x', y' \rangle) := X(x; x') \otimes Y(y; y') \ .$$

For two $V$-graphs transports $f : X \Longrightarrow X'$ and $g : Y \Longrightarrow Y'$ we take the biproduct of transports

$$f \otimes g : X \otimes Y \Longrightarrow X' \otimes Y'$$

defined with Carte tensor product of vertexes appointments

$$f_0 \times g_0 : X_0 \times Y_0 \to X'_0 \times Y'_0$$

and biproduct of edges set mappings

$$(f \otimes g)_{\langle x,y \rangle, \langle x',y' \rangle} := f_{x,x'} \otimes g_{y;y'} \ .$$

So we have defined a superfunctor over the Carte product of categories

$$Grph_V \times Grph_V \Longrightarrow Grph_V \ .$$

We assume that the biproduct $r \otimes s$ is monotone in values category $V$ for both arguments, so the tensor product of $V$-graphs transports $f \otimes g$ will be monotone also.

Let we have pairs of parallel graphs transports $f \geq f' : X \Longrightarrow X'$ and $g \geq g' : Y \Longrightarrow Y'$. These means that we have the inequalities for edges set mappings

$$f_{x,x'} : X(x;x') \to X'(f_0(x); f_0(x')) \ , \quad f'_{x,x'} : X(x;x') \to X'(f'_0(x); f'_0(x')) \ .$$

These mappings also are parallel as $f_0 = f'_0$, therefore

$$X'(f_0(x); f_0(x')) = X'(f'_0(x); f'_0(x')) \ .$$

Also we have the inequalities for another edges set mappings

$$g_{x,x'} : Y(x;x') \to Y'(g_0(x); g_0(x')) \ , \quad g'_{x,x'} : Y(x;x') \to Y'(g'_0(x); g'_0(x')) \ .$$

These mappings are parallel as $g_0 = g'_0$, therefore

$$Y'(g_0(x); g_0(x')) = Y'(g'_0(x); g'_0(x')) \ .$$

The inequalities
$$f_{x,x'} \geq f'_{x,x'} \ , \quad g_{y,y'} \geq g'_{y,y'}$$

provide the inequality for biproducts

$$f_{x,x'} \otimes g_{x,x'} \geq f'_{x,x'} \otimes g'_{y,y'} \ .$$

This means that we have inequality of parallel graphs transports

$$f \otimes g \geq f' \otimes g' : X \otimes Y \Longrightarrow X' \otimes Y' \ .$$

We shall deal with superfunctor defined by some *assistant V-graph T*. It for a $V$-graph $X$ appoints the biproduct $X \otimes T$ and for a $V$-graphs transport $f : X \Longrightarrow Y$ appoints the biproduct with the identity transport of the assistant graph

$$f \otimes \text{Id}_T : X \otimes T \Longrightarrow Y \otimes T \ .$$

Such appointment will maintain the inequality between parallel $V$-graphs transports $f, g : X \Longrightarrow Y$

$$f \geq g \Longrightarrow f \otimes \text{Id}_T \geq g \otimes \text{Id}_T \ .$$

We also shall define an *exponential superfunctor* $Y^T$. For a $V$-graph $Y$ and a small assistant $V$-graph $T$ we take the set $Y^T$ of all $V$-graphs transports $\phi : T \Longrightarrow Y$. They are defined with the appointments of vertexes $\phi_0 : T_0 \to X_0$ and the mappings of edges sets

$$\phi_{s,t} : T(s;t) \to Y(\phi_0(s); \phi_0(t))) \ .$$

The edges sets $Y^T(\phi; \psi)$ are defined as sets of transforms produced by product of edges sets in values category $V$

$$Y^T(\phi; \psi) = \prod_{s \in T_0} Y(\phi_0(s); \psi_0(s)) \ .$$

We can check that we have got indeed a superfunctor.

For the $V$-graphs transport $g : X \Longrightarrow Y$ we appoint the changing of target space. For a set of transforms between two $V$-graphs transports $\phi, \psi : T \Longrightarrow X$ we appoint the new set of transforms between changed transports

$$\phi \circ g, \psi \circ g : T \Longrightarrow Y \ .$$

The edges set mapping will be provided by the Carte product of mappings

$$\prod_{s \in T_0} X(\phi_0(s); \psi_0(s)) \to \prod_{s \in T_0} g_{\phi_0(s), \psi_0(s)}(X(\phi_0(s); \psi_0(s))) \subset$$

$$\prod_{s \in T_0} Y(g_0(\phi_0(s)); g_0(\psi_0(s))) \ .$$

This superfunctor also is monotone for the order relations between *secondary parallel V-graphs transports* of exponential spaces.

$$g, g' : Y^T(\phi; \psi) \Longrightarrow Y^T(\phi \circ g; \psi \circ g) \ .$$

For the inequality between parallel $V$-graphs transports $g, g' : X \Longrightarrow Y$ with an equality of projection

$$g_0(x) = g'_0(x) \uparrow x \in X_0$$

and inequalities for edges set mappings

$$g_{x,y} \geq g'_{x,y} \uparrow x, y \in X_0$$

we get the $V$-graphs transports $g^T, (g')^T : X^T \Longrightarrow Y^T$ with the equality for the *secondary projections*

$$(g_0)^{T_0} = (g'_0)^{T_0} : (X_0)^{T_0} \to (Y_0)^{T_0}$$

and inequalities for the $V$-graphs transports defining the appointment of the first projection. For a $V$-graphs transport $\phi \in (X^T)_0$ it appoint the new $V$-graphs transport $\phi \circ g \in (Y^T)_0$

$$\phi \circ g \geq \phi \circ g' \uparrow \phi \in X^T$$

with inequalities for edges set mappings

$$T(s;t) \to g_{\phi_0(s),\phi_0(t)}(\phi_{s,t}(T(s;t))) \geq g'_{\phi_0(s),\phi_0(t)}(\phi_{s,t}(T(s;t))) \ .$$

Also we get the inequalities for edges set mappings of $V$-graphs transports provided with changing of target space

$$X(x;y) \to \prod_{s \in T_0} g_{\phi_0(s),\psi_0(s)}(X(\phi_0(s);\psi_0(s))) \geq$$

$$\prod_{s \in T_0} g'_{\phi_0(s),\psi_0(s)}(X(\phi_0(s);\psi_0(s))) \subset \prod_{s \in T_0} Y(g_0(\phi_0(s));g_0(\psi_0(s))) \ .$$

We have got a joint pair of monotone superfunctors $X \otimes T$ and $Y^T$ within the category of $V$-graphs $Grph_V$. A *unit partial transform* for such joint pair of superfunctors is defined by collection of *sections transports* for $V$-graphs $X$ with chosen *unit arrow mappings* $* \to 1_x(*) \subset X(x;x)$. Also we need unit arrow mappings for the assistant $V$-graph $* \to 1_t(*) \subset T(t;t)$.

Sections transport $\lambda_X : X \Longrightarrow (X \otimes T)^T$ for a vertex $x \in X_0$ appoints the section $\lambda_X(x) : T \Longrightarrow X \otimes T$ which for a vertex $t \in T_0$ appoints the couple of vertexes $\langle x, t \rangle \in X_0 \times T_0$ and have edges set mappings

$$\lambda_X(a)_{s,t} : T(s;t) \to 1_a(*) \otimes T(s;t) \subset X(a;a) \otimes T(s;t)$$

defined as a biproduct of chosen unit arrow mappings $* \to 1_a(*) \subset X(a;a)$ with identity mapping $\mathtt{Id}_{T(s;t)} : T(s,t) \to T(s;t)$. For the sections transport an edges set mapping will be a product of biproducts for identity mappings $\mathtt{Id}_{X(a,b)} : X(a,b) \to X(a;b)$ with chosen unit arrow $* \to 1_s(*) \subset T(s;s)$

$$X(a;b) \to (\lambda_X)_{a,b}(X(a;b)) = \prod_{s \in T_0} (X(a;b) \otimes 1_s(*)) \subset (X \otimes T)^T(\lambda_X(a); \lambda_X(b)) \ .$$

There is no sense to speak about monotonous property for such section transports when we work with $V$-graphs. But name appointment of forms rests monotone
$$f^\# : X \xrightarrow{\lambda_X} (X \otimes T)^T \xrightarrow{f^T} Y^T .$$
It is enough to know that changing of target space is monotone for inequality between parallel $V$-graphs transports
$$f \geq g : X \otimes T \Longrightarrow Y .$$

A *counit partial transform* is defined by collection of *evaluation transports* for $V$-graphs $Y$ with *contact arrows mappings*
$$Y(x;y) \otimes Y(y;z) \to Y(x;y) \circ Y(y;z) \subset Y(x;z) .$$

For a couple of $V$-graphs transport $\phi : T \Longrightarrow Y$ and vertex $s \in T_0$ the evaluation transport appoints the value $\phi_0(s) \in Y_0$, and for a set of transform between two $V$-graphs transports $\phi, \psi : T \Longrightarrow Y$
$$\prod_{s \in T_0} Y(\phi_0(s); \psi_0(s)) \subset Y^T(\phi; \psi)$$
and an edges set $T(s;t)$ evaluation transport will be defined by mapping of *precontact arrows*
$$Y^T(\phi; \psi) \otimes T(s;t) \to Y(\phi_0(s); \psi_0(s)) \circ \psi_{s,t}(T(s;t)) \subset Y(\phi_0(s); \psi_0(t)) .$$

We shall demand that changing of target space with contact arrows mappings will be monotone in values category $V$. For arbitrary object $r \in V_0$ we ask monotone mappings
$$V(r; Y(x;y)) \to V(r; Y(x;y) \circ Y(y;z)) \subset V(r; Y(x;z)) ,$$
$$V(r; Y(y;z)) \to V(r; Y(x;y) \circ Y(y;z)) \subset V(r; Y(x;z)) .$$
That provides monotone realization for name forms $g : X \Longrightarrow Y^T$
$$g^\flat : X \otimes T \xrightarrow{g \otimes 1_T} Y^T \otimes T \xrightarrow{\mathbf{ev}_Y} Y .$$

Let we have an inequality between the secondary parallel name forms $g \geq g' : X \Longrightarrow Y^T$. Then we can check the inequality between realized forms. We shall check inequalities for edges set mappings
$$X(x;y) \otimes T(s;t) \to g_{x,y}(X(x;y)) \otimes T(s;t) \to$$
$$Y(g_0(x)_s; g_0(y)_s) \circ (g_0(y))_{s,t}(T(s;t)) \geq$$
$$g'_{x;y}(X(x;y)) \circ (g_0(y))_{s,t}(T(s;t)) \geq$$
$$g'_{x;y}(X(x;y)) \circ (g'_0(y))_{s,t}(T(s;t)) .$$

We are ready to prove the properties of name appointment or realization appointment between the forms into the target $V$-graph $Y$. For chosen unit arrows mappings in the source $V$-graph $* \to 1_a(*) \subset X(a;a)$ and the assistant $V$-graph $* \to 1_s(*) \subset T(s;s)$, and contact arrows mappings in the target $V$-graph
$$Y(x;y) \otimes Y(y;z) \to Y(x;y) \circ Y(y;z) \subset Y(x;z)$$
the $V$-graphs form $f : X \otimes Y \Longrightarrow Z$ will be called *lower decomposable* if we have inequality between mappings

$$X(a;b) \otimes T(s;t) \to f_{\langle a,s\rangle,\langle b,s\rangle}(X(a;b) \otimes 1_s(*)) \circ f_{\langle y,s\rangle,\langle y,t\rangle}(1_b(*) \otimes T(s;t)) \leq$$

$$f_{\langle x,s\rangle,\langle b,t\rangle}(X(a;b) \otimes T(s;t)) \subset Y(f(a,s); f(b,t))$$

Otherwise we shall say that a $V$-enriched graphs transport $g : X \Longrightarrow Y^T$ is *upper neutral* for chosen unit arrows if we have other inequalities between mappings

$$X(a;b) \to g_{a,b}(X(a;b))_s \circ g(b)_{s,s}(1_s(*)) \geq$$
$$g_{a,b}(X(a;b))_s \subset Y(g(a)_s; g(b)_s) ,$$
$$T(s;t) \to g_{a,a}(1_a(*)) \circ g(a)_{s,t}(T(s;t)) \geq$$
$$g(a)_{s,t}(T(s;t)) \subset Y(g(a)_s; g(a)_t) .$$

As name or realization appointments aren't totally defined, we shall use the notion of *composable partial appointments*. One partial appointment $\Phi : A \rightarrowtail B$ is composable with another partial appointment $\Psi : B \rightarrowtail C$ if the first appointment has values in the domain of definition for the second appointment.

**Proposition 4.** *The choosing of unit arrow mappings in the source $V$-graph $* \to 1_a(*) \subset X(a;a)$ and in the assistant $V$-graph $* \to 1_t(*) \subset T(s;s)$ allows us to define the name appointment between the sets of $V$-graphs transports*

$$Grph_V(X \otimes T; Y) \to Grph_V(X; Y^T) .$$

*The contact arrows mappings in the target $V$-graph $Y$*

$$Y(c;d) \otimes Y(d;e) \to Y(c;d) \circ Y(d;e) \subset Y(c;e)$$

*provide evaluation transport with precontact arrows mappings and allow us to define the realization appointment*

$$Grph_V(X; Y^T) \to Grph_V(X \otimes T; Y) .$$

*If contact arrows mappings are monotone, then the composable restriction of these appointments provides joint pair of monotone mappings from $V$-graphs transports over the biproduct $f : X \otimes T \Longrightarrow Y$ which are lower decomposable for chosen mapping of contact arrows in the target $V$-graph $Y$*

$$X(a;b) \otimes T(s;t) \to f_{\langle a,s\rangle,\langle b,s\rangle}(X(a;b) \otimes 1_s(*)) \circ f_{\langle y,s\rangle,\langle y,t\rangle}(1_b(*) \otimes T(s;t)) \leq$$

$$f_{\langle x,s \rangle, \langle b,t \rangle}(X(a;b) \otimes T(s;t)) \subset Y(f(a,s); f(b,t))$$

to $V$-graphs transports $g : X \Longrightarrow Y^T$ which are upper neutral for chosen unit arrows in the source $V$-graph $X$ and the assistant $V$-graph $T$

$$X(a;b) \to g_{a,b}(X(a;b))_s \circ g(b)_{s,s}(1_s(*)) \geq$$

$$g_{a,b}(X(a;b))_s \subset Y(g(a)_s; g(b)_s) ,$$

$$T(s;t) \to g_{a,a}(1_a(*)) \circ g(a)_{s,t}(T(s;t)) \geq$$

$$g(a)_{s,t}(T(s;t)) \subset Y(g(a)_s; g(a)_t) ,$$

i. e. we have inequalities for unit and counit transforms

$$f \geq f_+(g(f)) , \quad g(f_+(g)) \geq g .$$

**Proof:** We shall adapt the earlier categorical proving for the category of graphs **GRPH** in **G. Valiukevičius 2022**. We can see that the name appointment is monotone. The inequality $f \geq f'$ between $V$-graphs transports $f, f' : X \times T \Longrightarrow Y$ is maintained by exponential superfunctor

$$f^T \geq (f')^T$$

and by composition with sections transport

$$\lambda_X : X \Longrightarrow (X \otimes T)^T .$$

Also the realization appointment is monotone. The inequality $g \geq g'$ between $V$-graphs transports $g, g' : X \Longrightarrow Y^T$ is maintained by superfunctor of tensor product

$$f \otimes \text{Id}_T \geq g \otimes \text{Id}_T$$

and by composition with evaluation transport

$$\mathbf{ev}_Y : Y^T \otimes T \Longrightarrow Y .$$

The evaluation transport is monotone as a consequence of monotone contact arrows mappings in the target $V$-graph $Y$.

First we shall check that a realization of name is decreasing. For it we construct a continuous diagram

$$\begin{array}{ccccc} X \otimes T & \stackrel{\lambda_X \otimes 1_T}{\longrightarrow} & (X \otimes T)^T \otimes T & \stackrel{f^T \otimes 1_T}{\longrightarrow} & Y^T \otimes T \\ & \searrow & \downarrow \mathbf{ev}_{(X \otimes T)} & & \downarrow \mathbf{ev}_Y \\ & & X \otimes T & \stackrel{}{\underset{f}{\longrightarrow}} & Y \\ & & & & \nearrow \end{array}$$

The evaluation transport $\mathbf{ev}_{(X \otimes T)} : (X \otimes T)^T \otimes T \Longrightarrow X \otimes T$ is defined by choosing mappings of contact arrows

$$X(a;b) \otimes T(s;t) \to (X(a;b) \otimes *) \otimes (* \otimes T(s;t)) \to$$

$$(X(a;b) \otimes 1_s(*)) \circ (1_b(*) \otimes T(s;t)) = X(a;b) \otimes T(s;t) \ .$$

This is natural choice, but we need to check that with such choice we have got a suitable diagram. At first we shall check the first triangular equality for arbitrary $V$-graphs $X$ and $T$ with chosen unit arrow mappings $* \to 1_a(*) \subset X(a;a)$ and $* \to 1_t(*) \subset T(t;t)$. We shall prove the identic transport for composition of $V$-graphs transports

$$\mathtt{Id}_{(X \otimes T)} : X \otimes T \stackrel{\lambda_X \otimes 1_T}{\longrightarrow} (X \otimes T)^T \otimes T \stackrel{\mathbf{ev}_{(X \otimes T)}}{\longrightarrow} X \otimes T \ .$$

We apply the projection of categories $Grph_V \Longrightarrow \mathbf{SET}$ which for a $V$-graph $X$ appoints its vertexes set $X_0$, and for a $V$-graphs transport $f : X \Longrightarrow Y$ appoints the vertexes appointment $f_0 : X_0 \to Y_0$.

The constructed biproduct of $V$-graphs $X \otimes Y$ is projected to the Carte tensor product of sets $X_0 \times Y_0$.

For an exponential $V$-graph $Y^T$ such projection appoints the set of all $V$-graphs transports $\phi : T \Longrightarrow Y$. It will be called a *primary projection*. The *secondary projection* is defined with deformation of each $V$-graphs transport $\phi : T \Longrightarrow Y$ to its vertexes appointment $\phi_0 : T_0 \to Y_0$. Later we shall need to calculate these projections only for some $V$-graphs transports.

For the secondary projection we get the first triangular equality in the category of sets $\mathbf{SET}$. The evaluation appointment for a $V$-graphs transport $\phi : T \Longrightarrow Y$ and vertex $s \in T_0$ appoints the value $\phi_0(s) \in Y_0$ which depends only on projection $\phi_0 : T_0 \to Y_0$, so we get the first triangular equality also for the primary projection.

For a biproduct of edges sets $X(a;b) \otimes T(s;t)$ the first $V$-graphs transport

$$\lambda_X \otimes \mathtt{Id}_T : X \otimes T \Longrightarrow (X \otimes T)^T \otimes T$$

will appoint the biproduct of transforms set

$$\prod_{s \in T} (X(a;b) \otimes 1_s(*)) \subset (X \otimes T)^T (\lambda_X(a); \lambda_X(b))$$

and the same edges set in the assistant $V$-graph $T(s;t)$. The second $V$-graphs transport will be the evaluation transport on $V$-graphs biproduct $X \otimes T$. It for the former biproduct will appoint the contact arrows mapping

$$X(a;b) \otimes T(s;t) \to (X(a;b) \otimes 1_s(*)) \circ (1_b(*) \otimes T(s;t)) = X(a;b) \otimes T(s;t)) \ .$$

So we can end the checking of the first triangular equality in the category of $V$-graphs $Grph_V$.

Now we shall check for what $V$-graphs transport $f : X \otimes T \Longrightarrow Y$ we get continuous quadrat

$$\begin{array}{ccc} (X \otimes T)^T \otimes T & \stackrel{f^T \otimes 1_T}{\longrightarrow} & Y^T \otimes T \\ \mathbf{ev}_{(X \otimes T)} \downarrow & & \downarrow \mathbf{ev}_Y \\ X \otimes T & \stackrel{}{\underset{f}{\longrightarrow}} & Y \\ & & \nearrow \end{array}$$

We have the commuting quadrat in the category of sets **SET** for the secondary projection, also we have commuting quadrat for the primary projection. It rests to check the continuous diagram for the edges set mappings.

We get the continuous quadrat only for the diagram begun with sections mappings and with lower decomposable $V$-graphs transport $f: X \otimes T \Longrightarrow Y$.

For a set of transforms between two sections

$$\prod_{s \in T_0} (X(a;b) \otimes 1_s(*)) \subset (X \otimes T)^T(\lambda_X(a); \lambda_X(b))$$

the first evaluation $\mathbf{ev}_{(X \otimes T)}$ will be defined with the contact arrows mapping

$$(X(a;b) \otimes 1_s(*)) \circ (1_b(*) \otimes T(s;t)) = X(a;b) \otimes T(s;t) \subset (X \otimes T)(\langle a,s \rangle; \langle b,t \rangle)$$

and taken $V$-graphs transport $f: X \times T \Longrightarrow Y$ will appoint the edges set

$$f_{\langle a,s \rangle; \langle b,t \rangle}(X(a;b) \otimes T(s;t)) \subset Y(f(a,s); f(b,t)) .$$

The changing of the target space with $V$-graphs transport $f: X \otimes T \Longrightarrow Y$ for a set of transforms between sections

$$\prod_{s \in T_0} (X(a;b) \otimes 1_s(*)) \subset (X \otimes T)^T(T\lambda_X(a); \lambda_X(b))$$

will appoint the set of transforms between changed sections

$$\prod_{s \in T_0} f_{\langle a,s \rangle, \langle b,t \rangle}(X(a;b) \otimes 1_s(*)) \subset Y^T(\lambda_X(a) \circ f; \lambda_X(b) \circ f) ,$$

so for the $V$-graphs transport $f$ with lower decomposition the evaluation transport for a biproduct of such set of transforms with the same edges set in assistant graph $T(s;t)$ will appoint the smaller contact arrows mapping

$$f_{\langle a,s \rangle, \langle b,s \rangle}(X(a;b) \otimes 1_s(*)) \circ f_{\langle b,s \rangle, \langle b,t \rangle}(1_b(*) \otimes T(s,;t)) \leq$$

$$f_{\langle a,s \rangle, \langle b,t \rangle}(X(a;b) \otimes T(s;t)) \subset Y(f(a,s); f(b,t)) .$$

Next we shall check that a name of realization is also decreasing $g(f_+(g)) \geq g$. For it we construct another cocontinuous diagram

$$\swarrow Y^T \xleftarrow{(\mathbf{ev}_Y)^T} (Y^T \times T)^T \xleftarrow{(g \otimes 1_T)^T} (X \otimes T)^T$$
$$\nwarrow \quad \uparrow \lambda_{Y^T} \quad \uparrow \lambda_X$$
$$Y^T \xleftarrow{g} X$$

In the $V$-graph $Y$ we have chosen contact arrow mappings

$$Y(c;d) \otimes Y(d;e) \to Y(c;d) \circ Y(d;e) \subset Y(c;e) .$$

but for the sections transport $\lambda_{Y^T} : Y^T \implies (Y^T \otimes T)^T \otimes T$ we have to choose the unit arrows in exponential graph $Y^T$. It will be defined with chosen unit arrows mappings in $V$-graph $X$ composed with arrows sets mappings of taken $V$-graphs transport $g : X \implies Y^T$

$$g_{a,a}(1_a(*)) \subset Y^T(g_0(a); g_0(a)) \ .$$

We shall check the second triangular equality

$$Y^T \xrightarrow{\lambda_{Y^T}} (Y^T \otimes T)^T \xrightarrow{(\mathbf{ev}_Y)^T} Y^T \geq \mathtt{Id}_{Y^T} \ .$$

begun with taken $V$-graphs transport $g : X \implies Y^T$ which is upper neutral to chosen unit arrows mappings in source $V$-graph $X$

$$X(a; b) \to g_{a,b}(X(a; b))_s \circ g_0(b)_{s,s}(1_s(*)) \geq$$
$$g_{a,b}(X(a; b))_s \subset Y(g_0(a)_s; g_0(b)_s) \ ,$$
$$T(s; t) \to g_{a,a}(1_a(*)) \circ g_0(a)_{s,t}(T(s; t)) \geq$$
$$g_0(a)_{s,t}(T(s; t)) \subset Y(g_0(a)_s; g_0(a)_t) \ .$$

At first we shall check the second triangular equality in the category of sets **SET** for the set of all $V$-graphs transports $(Y^T)_0 = \{\phi : T \implies Y\}$

$$(Y^T)_0 \xrightarrow{(\lambda_{(Y^T)})_0} ((Y^T \times T)^T)_0 \xrightarrow{(\mathbf{ev}^T_{Y^T})_0} (Y^T)_0$$

begun with the projection of taken graphs transport $g_0 : X_0 \to (Y^T)_0$. It for a $V$-graphs transport $\phi : T \implies Y$ appoint the section

$$\lambda_{Y^T}(\phi) : T \implies Y^T \times T \ .$$

We can apply the secondary projection to get a new triangular equality for the same $V$-graphs transport $\phi : T \implies Y$, but with a new sections appointment in the category of sets **SET**

$$(Y^T)_0 \xrightarrow{\lambda_{(Y^T)_0}} ((Y^T)_0 \times T_0)^{T_0} \xrightarrow{\mathbf{ev}^{T_0}_{(Y^T)_0}} (Y^T)_0 \ .$$

It for a $V$-graphs transport $\phi : T \implies Y$ will appoint the section in the category of sets

$$\lambda_{(Y^T)_0}(\phi) : T_0 \to (Y^T)_0 \times T_0$$

which for a vertex $t \in T_0$ appoints the couple $\langle \phi, t \rangle \in (Y^T)_0 \times T_0$. This new triangular equality will provide the triangular equality also for the primary projection, as evaluation allows such secondary projection.

Indeed for a set of all mappings $\phi : T_0 \to Y_0$ we have already checked the second triangular equality

$$(Y_0)^{T_0} \xrightarrow{\lambda_{(Y_0)^{T_0}}} ((Y_0)^{T_0} \times T_0)^{T_0} \xrightarrow{\mathbf{ev}^{T_0}_{(Y_0)}} (Y_0)^{T_0} \ .$$

The first appointment will be sections appointment with arbitrary appointment $\phi : T_0 \to Y_0$
$$\lambda_{(Y_0)^{T_0}} : T_0 \to (Y_0)^{T_0} \times T_0 \ .$$

So it will be valid also for imbedded smaller set $(Y^T)_0 \subset (Y_0)^{T_0}$ of all appointments get from $V$-graphs transforms $\phi : T \Longrightarrow Y$.

For a $V$-graphs transport $g_0(a) : T \Longrightarrow Y$ we need else to check the inequality for the edges set mappings. An edges set $T(s;t)$ will be mapped by such graphs transports to the edges set

$$g_0(a)_{s,t}(T(s;t)) \subset Y(g_0(a)_s; g_0(a)_t) \ .$$

The sections transport for such graphs transport will appoint a section

$$\lambda_{Y^T}(g_0(a)) : T \Longrightarrow Y^T \otimes T$$

which for an edges set $T(s;t)$ appoints the biproduct

$$g_{a,a}(1_a(*)) \otimes T(s;t) \subset Y^T \otimes T(\langle g_0(a), s\rangle; \langle g_0(a), t\rangle) \ .$$

The changing of target space with evaluation transport is defined with pre-contact arrows mapping which for such biproduct will appoint the composition

$$g_{a,a}(1_a(*))_s \circ g_0(a)_{s,t}(T(s;t)) \geq g_{s,t}(T(s;t)) \subset Y(g_0(a)_s; g_0(a)_t) \ .$$

So for an edges set $T(s;t)$ we get decreasing mapping.

Finally we need to check the second triangular equality for the edges set mappings. Also we must begin with edges set mappings for taken $V$-graphs transport $g : X \Longrightarrow Y^T$. They for an edges set $X(a;b)$ appoints the set of transforms between two graphs transports

$$g_{a,b}(X(a;b)) \subset Y^T(g_0(a); g_0(b)) \ .$$

For such set of transforms we get the set of transforms between sections as a product

$$\prod_{s \in T_0} (g_{a,b}(X(a;b)) \otimes 1_s(*)) \subset (Y^T \otimes T)^T(\lambda_{Y^T}(g_0(a)); \lambda_{Y^T}(g_0(b))) \ .$$

The changing of target space with evaluation transport $\mathbf{ev}_Y : Y^T \otimes T \Longrightarrow Y$ will appoint the set of precontact arrows

$$g_{a,b}(X(a;b))_s \circ g_0(b)_{s,s}(1_s(*)) \geq g_{a,b}(X(a;b)) \subset Y(g_0(a)_s; g_0(b)_s) \ .$$

So we get decreasing mapping for the edges set $X(a;b)$ again.

We have shown the second triangular inequality begun with taken $V$-graphs transport $g : X \Longrightarrow Y^T$.

Now we shall check the commuting quadrat between defined sections transports

$$\begin{array}{ccc} (Y^T \otimes T)^T & \stackrel{(g \otimes 1_T)^T}{\longleftarrow} & (X \otimes T)^T \\ \lambda_{Y^T} \uparrow & & \uparrow \lambda_X \\ Y^T & \stackrel{}{\underset{g}{\longleftarrow}} & X \end{array}$$

This quadrat is commuting for the projection in the category **SET**. For a vertex $a \in X$ we appoint the $V$-graphs transport $g_0(a) \in Y^T$ which provides the section $\lambda_{Y^T}(g_0(a)) : T \Longrightarrow Y^T \otimes T$ which for a vertex $s \in T_0$ appoints the couple $\langle g_a(a), s \rangle$ with edges set mappings

$$\lambda_{Y^T}(g_0(a))_{s,t} = g_{a,a}(1_a(*)) \otimes \mathtt{Id}_{T(s;t)} : T(s;t) \to g_{a,a}(1_a(*)) \otimes T(s;t) \subset$$

$$Y^T(g_0(a); g_0(a)) \otimes T(s;t) .$$

The same section is get from the section $\lambda_X(a) : T \Longrightarrow X \otimes T$ by changing of target space with $V$-graphs transport $g \otimes \mathtt{Id}_T : X \otimes T \Longrightarrow Y^T \otimes T$. It for a couple of vertexes $\langle a, s \rangle \in X_0 \times T_0$ appoints the couple $\langle g_0(a), s \rangle \in (Y^T)_0 \times T_0$, and the edges set mappings

$$g_{a,a}(1_a(*)) \otimes \mathtt{Id}_{T(s;t)} : T(s;t) \to g_{a,a}(1_a(*)) \otimes T(s,t) \subset$$

$$Y^T(g_0(a); g_0(a)) \otimes T(s;t) .$$

For an edges set $X(a;b)$ taken $V$-graphs transport $g : X \Longrightarrow Y^T$ appoints the set of transforms between two $V$-graphs transports

$$g_{a,b}(X(a;b)) \subset Y^T(g_0(a); g_0(b))$$

which is mapped by sections transport $\lambda_{Y^T} : Y^T \Longrightarrow (T^T \otimes T)^T$ to the set of transforms between two sections

$$\prod_{s \in T_0} g_{a,b}((X(a;b)) \otimes 1_s(*) \subset (Y^T \otimes T)^T(\lambda_{Y^T}(g_0(a)); \lambda_{Y^T}(g_0(b))) .$$

The same set of transform between two sections is provided from the set of transforms

$$\prod_{s \in T_0} X(a;b) \otimes 1_s(*) \subset X^T(\lambda_X(a); \lambda_X(b))$$

by changing of target space with $V$-graphs transport $g \otimes \mathtt{Id}_T : X \times T \Longrightarrow Y^T \otimes T$. For an edges set

$$X(a;b) \otimes 1_s(*) \subset X(a;b) \otimes T(s;s)$$

it appoints the edges set

$$g_{a,b}(X(a;b)) \otimes 1_s(*) \subset Y^T(g_0(a); g_0(b)) \otimes T(s;s)$$

for every vertex $s \in T$.

We have got that the $V$-graphs transports $f \in Grph(X \times T; Y)$ with lower decomposition for chosen contact arrows mappings in the target graph $Y$ are *lower reversing* its name and the $V$-graphs transports $g \in Grph(X; Y^T)$ with upper neutral property for chosen unity arrow mappings in the source $V$-graph $X$ and the assistant $V$-graph $T$ are *upper reversing* its realization. For couple of composable graphs transforms $\langle f, g \rangle$ with both properties together we have get a truly joint pair of monotone name appointment and monotone realization appointment. □

We shall notice that in the case of decomposition equality

$$X(a;b) \otimes T(s;t) \to f_{\langle a,s\rangle,\langle b,s\rangle}(X(a;b) \otimes 1_s(*)) \circ f_{\langle y,s\rangle,\langle y,t\rangle}(1_b(*) \otimes T(s;t)) =$$

$$f_{\langle x,s\rangle,\langle b,t\rangle}(X(a;b) \otimes T(s;t)) \subset Y(f(a,s); f(b,t))$$

we shall get reversing the realization mapping

$$f = f_+(g(f)) \ .$$

In the case of neutral equalities

$$X(a;b) \to g_{a,b}(X(a;b))_s \circ g(b)_{s,s}(1_s(*)) =$$

$$g_{a,b}(X(a;b))_s \subset Y(g(a)_s; g(b)_s) \ ,$$

$$T(s;t) \to g_{a,a}(1_a(*)) \circ g(a)_{s,t}(T(s;t)) =$$

$$g(a)_{s,t}(T(s;t)) \subset Y(g(a)_s; g(a)_t)$$

we shall get reversing property for the name appointment

$$g(f_+(g)) = g \ .$$

We also want for $V$-enriched graphs to define the notion of the set of natural transforms. However for the values category $V$ with ordered edges sets more convenient to work with the set of continuous transforms.

In the exponential space of all $V$-graphs transports $Y^T$ we take the inverse image of order relation $(\geq) \subset V(r;s) \times V(r;s)$ for the names of two mappings provided by preevaluation and postevaluation transports

$$Y^T(\phi;\psi) \otimes T(s;t) \to Y(\phi_0(s); \psi_0(s)) \circ \psi_{s,t}(T(s;t)) \subset Y(\phi_0(s); \psi_0(t)) \ ,$$

$$Y^T(\phi;\psi) \otimes T(s;t) \to \phi_{s,t}(T(s;t)) \circ Y(\phi_0(s); \psi_0(s)) \subset Y(\phi_0(s); \psi_0(t)) \ .$$

So we get the inverse image

$$\to Y_+^{(T)}(\phi;\psi) \to Y^T(\phi;\psi) \Rightarrow \prod_{s,t \in T_0} Y(\phi_0(s); \psi_0(t))^{T(s;t)} \ .$$

For the first precontact evaluation we get the set of cocontinuous transforms $Y_+^{(T)}(\phi; \psi)$

$$\begin{array}{ccc} \phi_0(a) & \xrightarrow{\alpha_a} & \psi_0(a) \\ \downarrow & (+13) & \downarrow \\ \phi_0(b) & \xrightarrow[\alpha_b]{} & \psi_0(b) \\ & & \swarrow \end{array}$$

The sign $(+)$ must denote the property of deformed quadrat $(+13)$ for *wasting reforms*

Otherwise for the first postcontact evaluation we shall get the set of continuous transforms $Y_-^{(T)}(\phi; \psi)$

$$\begin{array}{ccc} \phi_0(a) & \xrightarrow{\alpha_a} & \psi_0(a) \\ \downarrow & (+13) & \downarrow \\ \phi_0(b) & \xrightarrow[\alpha_b]{} & \psi_0(b) \\ & & \nearrow \end{array}$$

The sign $(-)$ must denote the property of deformed quadrat $(-13)$ for *respecting reforms*, At this moment we can't to talk about new superfunctors of exponential spaces $Y_-^{(T)}$ or $Y_+^{(T)}$.

We repeat the proposition when the name of $V$-graphs transport can be corestricted in the new exponential space of continuous transforms.

**Proposition 5.** *The choosing of unit arrow mappings in the source $V$-graph $1_a(*) \subset X(a; a)$ and in the assistant $V$-graph $1_s(*) \subset T(s; s)$ allows us to define the name appointment between the sets of $V$-graphs transports*

$$Grph_V(X \otimes T; Y) \to Grph_V(X; Y^T) \ .$$

*For a contact arrows mapping in the target $V$-graph $Y$*

$$Y(c; d) \otimes Y(d; e) \to Y(c; d) \circ Y(d; c) \subset Y(c; e)$$

*the choosing of evaluation transport defined with precontact or postcontact arrows mappings allows us to define the prerealization or postrealization appointments*

$$Grph_V(X; Y^T) \to Grph_V(X \otimes T; Y) \ .$$

*The name of $V$-graphs transport $f : X \otimes T \Longrightarrow Y$ which has upper predecomposition and lower postdecomposition for the contact arrows mapping in the target $V$-graph $Y$*

$$f_{\langle a,s \rangle, \langle b,s \rangle}(X(a; b) \otimes 1_s((*))) \circ f_{\langle b,s \rangle, \langle b,t \rangle}(1_b(*) \otimes T(s; t)) \geq$$

$$f_{\langle a,s \rangle, \langle b,t \rangle}(X(a; b) \otimes T(s; t)) \geq$$

$$f_{\langle a,s \rangle, \langle a,t \rangle}(1_a(*) \otimes T(s; t)) \circ f_{\langle a,t \rangle, \langle b,t \rangle}(X(a; b) \otimes 1_t(*))$$

is corestricted in exponential $V$-graphs of continuous transforms $Y_{-}^{(T)}$.

**Proof:**

We want to check for the name appointment the inclusion in inverse image of order relation

$$g_{a,b}X(a;b) \to Y^T(g_0(a); g_0(b)) \Rightarrow \prod_{s,t \in T_0} Y(g_0(a)_s; g_0(b)_t)^{T(s;t)} \ .$$

provided by preevaluation and postevaluation transports

$$X \otimes T \xrightarrow{\lambda_X \otimes 1_T} (X \otimes T)^T \otimes T \xrightarrow{f^T \otimes 1_T} Y^T \otimes T \xrightarrow{\mathbf{ev}_Y^\pm} Y$$

Using the cocontinuous diagram for the preevaluation transport

$$X \otimes T \xrightarrow{\lambda_X \otimes 1_T} (X \otimes T)^T \otimes T \xrightarrow{f^T \otimes 1_T} Y^T \otimes T$$
$$\searrow \quad \downarrow \mathbf{ev}_{(X \otimes T)}^+ \quad \downarrow \mathbf{ev}_Y^+$$
$$X \otimes T \xrightarrow{f} Y$$

and the continuous diagram for the postevaluation transport

$$X \otimes T \xrightarrow{\lambda_X \otimes 1_T} (X \otimes T)^T \otimes T \xrightarrow{f^T \otimes 1_T} Y^T \otimes T$$
$$\searrow \quad \downarrow \mathbf{ev}_{(X \otimes T)}^- \quad \downarrow \mathbf{ev}_Y^-$$
$$X \otimes T \xrightarrow{f} Y$$

allows us to check the *corestriction* of name

$$g : X \xrightarrow{\lambda_X} (X \otimes T)^T \xrightarrow{f^T} Y^T$$

in exponential space of continuous transforms $Y_{-}^{(T)}$.

This can be proved by contact arrows mappings

$$\begin{array}{ccc} f_0(a,s) & \to & f_0(a,t) \\ \downarrow & \searrow & \downarrow \\ f_0(b,s) & \to & f_0(b,t) \nearrow \\ & & \nearrow \end{array}$$

Such contact arrow is provided by

$$f_{\langle a,s \rangle, \langle b,t \rangle}(X(a;b) \otimes T(s;t)) : f_0(a,s) \to f_0(b,t)$$

from continuous diagrams drown for both evaluation transports $\mathbf{ev}_Y^\pm$

$$f_{\langle a,s \rangle, \langle b,s \rangle}(X(a;b) \otimes 1_s((*))) \circ f_{\langle b,s \rangle, \langle b,t \rangle}(1_b(*) \otimes T(s;t)) \geq$$

$$f_{\langle a,s \rangle, \langle b,t \rangle}(X(a;b) \otimes T(s;t)) \geq$$

$$f_{\langle a,s\rangle,\langle a,t\rangle}(1_a(*) \otimes T(s;t)) \circ f_{\langle a,t\rangle,\langle b,t\rangle}(X(a;b) \otimes 1_t(*)) \; ,$$

i. e. we have inequalities for the precontact and the postcontact arrows

□

We can also notice that in the case of equalities

$$f_{\langle a,s\rangle,\langle b,s\rangle}(X(a;b) \otimes 1_s((*)) \circ f_{\langle b,s\rangle,\langle b,t\rangle}(1_b(*) \otimes T(s;t)) =$$

$$f_{\langle a,s\rangle,\langle b,t\rangle}(X(a;b) \otimes T(s;t)) =$$

$$f_{\langle a,s\rangle,\langle a,t\rangle}(1_a(*) \otimes T(s;t)) \circ f_{\langle a,t\rangle,\langle b,t\rangle}(X(a;b) \otimes 1_t(*))$$

we get a name appointment corestricted in the exponential space of natural transforms $Y^{(T)}$.

We again shall define the *category of original V-graphs* $OGrph_V$. For the original V-graph we demand the unit arrow mappings $e_a(*) \subset X(a;a)$. The *original V-graphs transport* $f : X \Longrightarrow Y$ must maintain chosen unit arrow mappings

$$f_{a,a}(e_a(*)) = e_{f_0(a)}(*) \; .$$

We shall call them also as *crude functors*. Such category will be denoted $OGrph_V$.

For contact arrows mappings

$$X(a;b) \otimes X(b;c) \to X(a;b) \circ X(b;c) \subset X(a;c)$$

we demand neutral equalities

$$e_a(*) \circ X(a;b) = X(a;b) = X(a;b) \circ 1_b(*) \; .$$

The original V-graphs transport $f : X \Longrightarrow Y$ will be called continuous if we additionally have continuous diagram for chosen contact arrows mappings

$$\begin{array}{ccc} X(a;b) \otimes X(b;c) & \stackrel{f_{a,b} \otimes f_{b,c}}{\longrightarrow} & Y(f_0(a); f_0(b)) \otimes Y(f_0(b); f_0(c)) \\ \circ \downarrow & & \downarrow \circ \\ X(a;c) & \stackrel{}{\underset{f_{a,c}}{\longrightarrow}} & Y(f_0(a); f_0(c)) \\ & & \nearrow \end{array}$$

The *continuous original V-graphs transports* will be also called *continuous functors*. The category of continuous original V-graphs transports will be denoted $OGrph_V^-$.

Correspondingly we define *cocontinuous functors* and the category of cocontinuous original V-graphs transports will be denoted $OGrph_V^+$. *natural functors* are defined as transports of original V-graphs transports maintaining contact arrows mappings. Their category will be denoted $OGrph_V^=$. They are usually considered as functors in V-enriched categories.

A biproduct of original $V$-graphs $X \otimes Y$ provides the superfunctor also in the category of original $V$-graphs $OGrph_V$. The unit arrow mapping is defined by the biproduct of unit arrow mappings

$$1_{\langle x,t \rangle}(*) = 1_x(*) \otimes 1_s(*) \subset X \otimes T(\langle x,s \rangle; \langle x,s \rangle) \ .$$

For the choosing of contact arrows mappings we need additionally to demand a mapping in the values category $V$

$$m_{\langle p,s \rangle, \langle r,t \rangle} : (p \otimes s) \otimes (r \otimes t) \to (p \otimes r) \otimes (s \otimes t) \ .$$

We can say that such biproduct $p \otimes r$ is a *relator to itself*. Usually such arrow is provided for associative and commuting tensor product. We shall call it a *relator's involution*. Such arrows are inverse to their selves, therefore they define isomorphic transform.

The choosing of contact arrows mappings in the biproduct graph $X \otimes Y$ is defined with the biproduct of contact arrows mappings. Let we have the contact arrows mappings

$$X_{a,b,c} : X(a;b) \otimes X(b;c) \to X(a;b) \circ X(b;c) \subset X(a;c) \ ,$$

$$Y_{a',b',c'} : Y(a';b') \otimes Y(b';c') \to Y(a';b') \circ Y(b';c') \subset Y(a';c') \ .$$

Then in the biproduct $X \otimes Y$ we can define the contact arrows mappings

$$(X \otimes Y(\langle a,a' \rangle; \langle b,b' \rangle)) \otimes (X \otimes Y(\langle b,b' \rangle; \langle c,c' \rangle)) =$$

$$(X(a;b) \otimes Y(a';b')) \otimes (X(b;c) \otimes Y(b';c')) \to$$

$$(X(a;b) \otimes X(b;c)) \otimes (Y(a';b') \otimes Y(b';c')) \to$$

$$(X(a;b) \circ X(b;c)) \otimes (Y(a';b') \circ Y(b';c')) \subset$$

$$X(a;c) \otimes Y(a';c') = X \otimes Y(\langle a,a' \rangle; \langle c,c' \rangle) \ .$$

We need to check that chosen unit arrows mappings

$$1_a(*) \otimes 1_{a'}(*) \subset X(a;a) \otimes Y(a';a')$$

are neutral for such defined arrows mappings. First we have an equality

$$(1_a(*) \otimes 1_{a'}(*)) \otimes (X(a;b) \otimes Y(a';b')) \to$$

$$(1_a(*) \otimes X(a;b)) \otimes (1_{a'}(*) \otimes Y(a';b') \to$$

$$(1_a(*) \circ X(a;b)) \otimes (1_{a'}(*) \circ Y(a';b')) = X(a;b) \otimes Y(a';') \ .$$

In the same way we can check another equality

$$(X(a;b) \otimes Y(a';b')) \circ (1_a(*) \otimes 1_{a'}(*)) = X(a;b) \otimes Y(a';b') \ .$$

The biproduct with an assistant original $V$-graph $T$ defines superfunctor in various categories of original $V$-graphs. For an original $V$-graph $X$ it appoints the biproduct $X \otimes T$, and for an original $V$-graphs transport $g : X \Longrightarrow Y$ it appoints the biproduct with the identity transport of assistant $V$-graph

$$g \otimes \mathtt{Id}_T : X \otimes T \Longrightarrow Y \otimes T \ .$$

For a general original $V$-graphs transport we need to check only maintenance of unit arrow mappings. For the original $V$-graphs transport $g : X \Longrightarrow Y$ we shall have maintained the unit arrow mappings

$$g_{a,a}(1_a(*)) = 1_{g_0(a)} \subset Y(g_0(a); g_0(a)) \ .$$

So the unit arrow mappings in the biproduct

$$1_a(*) \otimes 1_s(*) \subset X(a;a) \otimes T(s;s)$$

will be mapped to the unit arrow mapping in another biproduct

$$(g \otimes \mathtt{Id}_T)_{\langle a,s \rangle, \langle a,s \rangle}(1_a \otimes 1_s) = g_{a,a}(1_a(*)) \otimes 1_s(*) = 1_{g_0(a)}(*) \otimes 1_s(*) \subset$$

$$Y(g_0(a); g_0(a)) \otimes T(s;s) \ .$$

For the continuous original $V$-graphs transforms $g : X \Longrightarrow Y$

$$\begin{array}{ccc} X(a;b) \otimes X(b;c) & \stackrel{g_{a,b} \otimes g_{b,c}}{\longrightarrow} & Y(g_0(a); g_0(b)) \otimes Y(g_0(b); g_0(c)) \\ \circ \downarrow & & \downarrow \circ \\ X(a;c) & \stackrel{}{\underset{g_{a,c}}{\longrightarrow}} & Y(g_0(a); g_0(c)) \\ & & \nearrow \end{array}$$

we need to ask that chosen arrows $m_{\langle p,s \rangle, \langle r,t \rangle}$ would define cocontinuous transform for every argument changing $p \to p'$, $r \to r'$, $s \to s'$, $t \to t'$, for example we have cocontinuous quadrats

$$\begin{array}{ccc} (p \otimes s) \otimes (r \otimes t) & \stackrel{m_{\langle p,s \rangle, \langle r,t \rangle}}{\longrightarrow} & (p \otimes r) \otimes (s \otimes t) \\ \downarrow & & \downarrow \\ (p' \otimes s) \otimes (r \otimes t) & \stackrel{}{\underset{m_{\langle p',s \rangle, \langle r,t \rangle}}{\longrightarrow}} & (p' \otimes r) \otimes (s \otimes t) \\ & & \swarrow \end{array}$$

Then the appointed biproduct $g \otimes \mathtt{Id}_T : X \otimes T \Longrightarrow Y \otimes T$ rests continuous for defined contact arrows mappings. The proving is produced gradually by each argument. First we shall check for the composition in the first space. We have the cocontinuous diagrams in the values category $V$

$$\begin{array}{cc} (X(a;b) \otimes T(s;t)) \otimes (X(b;c) \otimes T(t;v)) & \stackrel{m_{X,T}}{\longrightarrow} \\ (g_{a,b} \otimes 1_{T(s;t)}) \otimes (g_{b,c} \otimes 1_{T(t;v)}) \downarrow & \\ (Y(g_0(a); g_0(b)) \otimes (T(s;t)) \otimes Y(g_0(b); g_0(c)) \otimes T(t;v)) & \stackrel{}{\underset{m_{Y',T}}{\longrightarrow}} \end{array}$$

$$
\begin{array}{rl}
\to & (X(a;b) \otimes X(b;c)) \otimes (T(s;t) \otimes T(t;v)) \xrightarrow{\circ} \\
& (g_{a,b}\otimes g_{b,c})\otimes(1_{T(s;t)}\otimes 1_{T(t;v)}) \downarrow \\
\to & (Y(g_0(a);g_0(b)) \otimes Y(g_0(b);g_0(c))) \otimes (T(s;t) \otimes T(t;v)) \xrightarrow{\circ} \\
& \swarrow
\end{array}
$$

$$
\begin{array}{rlll}
\to & X(a;c) \otimes (T(s;t) \otimes T(t;v)) & \xrightarrow{\circ} & X(a;c) \otimes T(s;v) \\
& 1_{X(a;c)}\otimes(1_{T(s;t)}\otimes 1_{T(t;v)}) \downarrow & & \downarrow (1_{X(a;c)}\otimes 1_{T(s;v)}) \\
\to & X(a;c) \otimes (T(s;t) \otimes T(t;v)) & \xrightarrow{\circ} & X(a;c) \otimes T(s;v)) \\
& \swarrow & & \swarrow
\end{array}
$$

In the values category the arrows composition is associative, therefore we can compose these three cocontinuous diagrams, and we get the continuous diagram for the arrows set mappings of $V$-graphs transport $g : X \Longrightarrow Y$

$$
\begin{array}{cc}
(X(a;b)\otimes T(s;t)) \otimes (X(b;c) \otimes T(t;v)) & \xrightarrow{(g_{a,b}\otimes 1_{T(s;t)})\otimes(g_{b,c}\otimes(1_{T(s;v)}))} \\
\downarrow & \\
X(a;c) \otimes T(s,v) & \xrightarrow[g_{a,c}\otimes 1_{T(s;v)}]{} 
\end{array}
$$

$$
\begin{array}{rl}
\to & (Y(g_0(a);g_0(b)) \otimes T(s;t)) \otimes (Y(g_0(b);g_0(c)) \otimes T(t;v)) \\
& \downarrow \\
\to & Y(g_0(a);g_0(c)) \\
& \nearrow
\end{array}
$$

We remark that in the category $OGrph_V^{\bar{=}}$ compounded by *natural functors* between original $V$-graphs we need to ask that arrows $m_{\langle p,s\rangle,\langle r,t\rangle}$ would define a natural transform.

The *exponential space* $Y^T$ in the category of original $V$-graphs $OGrph_V$ will be defined by the vertex set compounded by all original $V$-graphs transports $\phi : T \Longrightarrow Y$, i. e. these $V$-graphs transports must maintain unit arrow mappings. However in a smaller category of continuous $V$-graphs transports $OGrph_V^-$ we shall take only continuous original $V$-graphs transports $\phi : T \Longrightarrow Y$. Their set will be denoted $Y^{T-}$. Also we define the set of cocontinuous original $V$-graphs transports $Y^{T+}$ in the category $OGrph_V^+$, or the set of original $V$-graphs natural transports which maintain contact arrows mappings $Y^{T=}$ in the category $OGrph_V^{\bar{=}}$.

The edges sets in exponential space of crude functors will be taken the same as in the category $Grph_V$ compounded of all transforms

$$Y^T(\phi;\psi) = \prod_{s \in T_0} Y(\phi_0(s); \psi_0(s)) \ .$$

We shall demand that biproduct in values category $V$ would maintain such product

$$\prod_{s\in T_0}(p_s \otimes r_s) = (\prod_{s\in T_0} p_s) \otimes (\prod_{s\in T_0} r_s) \ .$$

The unit arrow mappings in exponential space are provided by the transform defined with a product of unit arrow mappings in the target $V$-graph

$$1_\phi(*) = \prod_{s\in T_0} 1_{\phi_0(s)}(*) \subset \prod_{s\in T_0} Y(\phi_0(s); \phi_0(s)) \ .$$

The contact arrows mappings also are defined by a product of such mappings in the target $V$-graph $Y$

$$Y^T(\phi;\psi) \otimes Y^T(\psi;\xi) \to (\prod_{s\in T_0} Y(\phi_0(s); \psi_0(s))) \circ (\prod_{s\in T_0} Y(\psi_0(s); \xi_0(s))) =$$

$$\prod_{s\in T_0}(Y(\phi_0(s); \Psi_0(s)) \circ Y(\Psi_0(s); \xi_0(s))) \subset \prod_{s\in T_0} Y(\phi_0(s); \xi_0(s)) \ .$$

Easy to check that chosen unit arrow mappings are neutral for such contact arrows mappings

$$1_\phi(*) \circ Y^T(\phi;\psi) = \prod_{s\in T_0} 1_{\phi_0(s)} \circ Y(\phi_0(s); \psi_0(s)) =$$

$$\prod_{s\in T_0} Y(\phi_0(s); \psi_0(s)) = Y^T(\phi;\psi) \ .$$

In the same way we can check another equality

$$Y^T(\phi;\psi) \circ 1_\psi(*) = Y^T(\phi;\psi) \ .$$

The exponential space will define superfunctors in correspondent categories. For a $V$-graphs transport $g : Y \Longrightarrow Y'$ we take the changing of target space $f^T : Y^T \Longrightarrow (Y')^T$. For each $V$-graph transport $\phi : T \Longrightarrow Y$ we get the new $V$-graphs transport $\phi \circ g : T \Longrightarrow Y'$. Such composition of $V$-graphs transport rests in correspondent categories. For example for the continuous original $V$-graphs transport $\phi : T \Longrightarrow Y$ we get again a continuous original $V$-graphs transport $\phi \circ g : T \Longrightarrow Y'$.

The edges set mappings are defined in old way

$$Y^T(\phi;\psi) \to \prod_{s\in T_0} g_{\phi_0(s),\psi_0(s)}(Y(\phi_0(s); \psi_0(s))) \subset$$

$$\prod_{s\in T_0} Y'(g_0(\phi_0(s)); g_0(\psi_0(s))) = (Y')^T(\phi \circ g; \psi \circ g) \ .$$

Finally we shall check that continuous $V$-graphs transports $f : X \Longrightarrow Y$ provides continuous changing of target space $f^T : X^T \Longrightarrow Y^T$. Let we have continuous quadrat

$$
\begin{array}{ccc}
X(x;y) \otimes X(y;z) & \stackrel{f_{x,y} \otimes f_{y,z}}{\longrightarrow} & Y(f_0(x); f_0(y)) \otimes Y(f_0(y); f_0(z)) \\
\circ \downarrow & & \downarrow \circ \\
X(x;z) & \underset{f}{\longrightarrow} & Y(f_0(x); f_0(z)) \\
& & \nearrow
\end{array}
$$

Then for the changing of target space we get again the continuous quadrat

$$
\begin{array}{ccc}
X^T(\phi;\psi) \otimes X^T(\psi;\xi) & \stackrel{f^T_{\phi;\psi} \otimes f^T_{\psi;\xi}}{\longrightarrow} & Y^T(\phi \circ f, \psi \circ f) \circ Y^T(\psi \circ f; \xi \circ f) \\
\circ \downarrow & & \downarrow \circ \\
Y^T(\phi \circ f; \xi \circ f) & \underset{f^T}{\longrightarrow} & Y^T(\phi \circ f; \xi \circ f) \\
& & \nearrow
\end{array}
$$

The needed inequality is checked for each component with $s \in T_0$. We apply the first quadrat for the vertexes $x = \phi_0(s)$, $y = \phi_0(s)$ and $z = \xi_0(s)$.

We have got a joint pair of superfunctors defined by biproduct $X \otimes T$ and exponential space $Y^T$ in correspondent categories. For such joint pair of superfunctors we can construct the unit transform defined with collection of sections transports

$$\lambda_X : X \Longrightarrow (X \otimes T)^T \ .$$

We need only to check that sections $\lambda_X(a) : T \Longrightarrow X \otimes T$ belongs to the correspondent category. The section for the vertex $s \in T_0$ appoints the couple $\langle a, s \rangle \in X_0 \times T_0$ and has edges set mappings

$$T(s;t) \to 1_a(*) \otimes T(s;t) \subset X(a;a) \otimes T(s;t) \ .$$

For the original $V$-graphs transport we need to check only maintenance of unit arrow mappings

$$1_s(*) \to 1_a(*) \otimes 1_s(*) \subset X(a;a) \otimes T(s;s) \ .$$

We shall check that section's edges set mappings are natural for the contact arrows mappings

$$T(s;t) \otimes T(t;v) \to (1_a(*) \otimes T(s;t)) \otimes (1_a(*) \otimes T(t;v)) \to$$

$$(1_a(*) \otimes 1_a(*)) \otimes (T(s;t) \otimes T(t;v)) \to (1_a(*) \circ 1_a(*)) \otimes (T(s;t) \circ T(t;v)) =$$

$$1_a(*) \otimes (T(s,t) \circ T(t;v)) \subset X(a;a) \otimes T(s;t) \ .$$

This can be express with commuting diagram

$$
\begin{array}{ccc}
T(s;t) \otimes T(t;v) & \to & (1_a(*) \otimes T(s;t)) \otimes (1_a(*) \otimes T(t;v)) \\
\circ \downarrow & & \downarrow \circ \\
T(s;v) & \to & 1_a(*) \otimes T(s;v)
\end{array}
$$

The sections transport also is natural for the contact arrows mappings

$$X(x;y) \otimes X(y;z) \to (\prod_{s \in T_0} X(x;y) \otimes 1_s(*)) \otimes (\prod_{s \in T_0} X(y;z) \otimes 1_s(*)) =$$

$$\prod_{s \in T_0} (X(x;y) \otimes 1_s(*)) \otimes (X(y;z)) \otimes 1_s(*)) \to$$

$$\prod_{s \in T_0} (X(x;y) \otimes X(y;z)) \otimes (1_s(*) \otimes 1_s(*)) \to$$

$$\prod_{s \in T_0} (X(x;y) \circ X(y;z)) \otimes (1_s(*) \circ 1_a(*)) =$$

$$\prod_{s \in T_0} (X(x;y) \circ X(y;z)) \otimes 1_s(*) \subset \prod_{s \in T_0} X(x;z) \otimes T(s;s) \ .$$

This can be express by another commuting diagram

$$\begin{array}{ccc} X(x;y) \otimes X(y;z) & \to & \prod_{s \in T_0} X(x;y) \otimes 1_s(*)) \otimes \prod_{s \in T_0} X(y;z) \otimes 1_s(*) \\ \circ \downarrow & & \downarrow \circ \\ X(x;z) & \to & \prod_{s \in T_0} X(x;z) \otimes 1_s(*) \end{array}$$

We also can check that collection of sections transports defines natural transforms in the category of original $V$-graphs transports $Grph_V$, i. e. for arbitrary original $V$-graphs transport $f : X \Longrightarrow X'$ we have commuting diagram

$$\begin{array}{ccc} X & \xrightarrow{f} & X' \\ \lambda_X \downarrow & & \downarrow \lambda_{X'} \\ (X \otimes T)^T & \xrightarrow[(f \otimes 1_T)^T]{} & (X' \otimes T)^T \end{array}$$

First we shall check the commuting quadrat for appointments of vertexes. For the vertex $x \in X_0$ we get a section $\lambda_X(x) : T \Longrightarrow X \otimes T$ with edges mappings

$$T(s;t) \to 1_x(*) \otimes T(s;t) \subset X \otimes T(\langle x,s \rangle; \langle x,t \rangle) \ .$$

The original $V$-graphs transport provided by changing of target space with the original $V$-graphs transport $f \otimes 1_T : X \otimes T \Longrightarrow X' \otimes T$ appoints another section $T \Longrightarrow X' \times T$ with edges set mappings

$$T(s;t) \to f_{x;x}(1_x(*)) \otimes T(s;t) =$$

$$1_{f_0(x)}(*) \otimes T(s;t) \subset X'(f_0(x); f_0(x)) \otimes T(s;t) \ .$$

So for the vertexes appointment we get the commuting quadrat.

For the edges set $X(x;y)$ the sections mapping $\lambda_X : X \Longrightarrow (X \otimes T)^T$ provides the set of transforms between two sections

$$X(x;y) \to \prod_{s \in T_0} X(x;y) \otimes 1_s(*) \subset \prod_{s \in T_0} X(x;y) \otimes T(s;s) \ .$$

The changing of target space provided with the original $V$-graphs transport $f \times 1_T : X \times T \Longrightarrow X' \otimes T$ produces the edges set mapping

$$X(x;y) \to \prod_{s \in T_0} f_{x,y}(X'(x;y)) \otimes 1_s(*) \subset \prod_{s \in T_0} X'(f_0(x); f_0(y)) \otimes T(s;s) \ .$$

It is the same as the edges set mapping for composition for $V$-graphs transport $f : X \Longrightarrow X'$ with sections transport $\lambda_{X'} : X' \Longrightarrow (X' \otimes T)^T$.

This commuting diagram remains also for the transports $f : X \Longrightarrow X'$ from smaller categories $OGrph_V^-$, $OGrph_V^+$ or $OGrph_V^=$.

The counit transform can be provided by preevaluation or postevaluation transports. They both will be original $V$-graphs transports. For the unit arrow mappings $1_\phi \in Y^T(\phi;\phi)$ from original $V$-graphs transport $\phi : T \Longrightarrow Y$ and unit arrow mapping $1_s(*) \subset T(s;s)$ we get preevaluation with precontact arrows mapping

$$1_{\phi_0(s)}(*) \circ \phi_{s,s}(1_s(*)) = 1_{\phi_0(s)}(*) \circ 1_{\phi_0(s)}(*) = 1_{\phi_0(s)}(*) \subset Y(\phi_0(s); \phi_0(s)) \ .$$

The same checking can be repeated for postcontact arrows mapping.

However we can't define an evaluation transport in smaller categories of continuous or cocontinuous original $V$-graphs transports. So we have the name appointment in such categories

$$OGrph_V^-(X \otimes T; Y) \to OGrph_V^-(X; Y^{T-})$$

but we haven't any realization appointment

$$OGrph_V^-(X; Y^{T-}) \to OGrph_V^-(X \otimes T; Y) \ .$$

Later we shall need some additional properties to get continuity property for $V$-graphs evaluation transports over biproduct $X \otimes T$.

For the original $V$-graphs we may do some useful calculations. We shall call a original $V$-graphs transport $f : X \otimes T \Longrightarrow Y$ *preregular* if the prerealization appointment reverses its name

$$f = f_+(g(f)) \ .$$

For an original $V$-graphs transport $f$ with upper predecomposition we get an inequality with its name's prerealization, which will be preregular $V$-graphs transport

$$f \leq f_+(g(f)) = f^+ \ .$$

This upper bound is a smallest among preregular bounds.

Otherwise we shall call an original $V$-graphs transport $f : X \otimes T \Longrightarrow Y$ *postregular* if the postrealization mapping reverses its name

$$f = f_-(g(f)) \ .$$

For a original $V$-graphs transport $f$ with lower postdecomposition we get an inequality with its name's postrealization, which will be postregular $V$-graphs transport

$$f \geq f_-(g(f)) = f^- .$$

This lower bound is a biggest among postregular lower bounds. Taken $V$-graphs transport $f : X \otimes T \Longrightarrow Y$ has a name corestricted in the space of natural transform exactly when it coincides with its upper preregular and lower postregular bounds

$$f^- = f = f^+ .$$

We shall prove only one of these sentences.

**Proposition 6.** *Let we have a original $V$-graphs transport $f : X \otimes T \Longrightarrow Y$ with upper predecomposition*

$$f_{\langle a,s\rangle,\langle b,s\rangle}(X(a;b) \otimes 1_s(*)) \circ f_{\langle b,s\rangle,\langle b,t\rangle}(1_b(*) \otimes T(s;t)) \geq$$

$$f_{\langle a,s\rangle,\langle b,t\rangle}(X(a;b) \otimes T(s;t)) .$$

*Then appointed prerealization of name $f^+ = f_+(g(f)) \geq f$ is the smallest preregular upper bound.*

**Proof:** We have calculated $f_+(g(f)) : X \times T \Longrightarrow Y$ as the composition

$$X \otimes T \xrightarrow{\lambda_X \otimes 1_T} (X \otimes T)^T \otimes T \xrightarrow{f^T \otimes 1_T} Y^T \otimes T \xrightarrow{\mathbf{ev}_Y^+} Y .$$

First sections transport for a vertex $a \in X_0$ will appoint the section $\lambda_X(a) : T \Longrightarrow X \otimes T$ with edges set mappings

$$\lambda_X(a)_{s,t}(T(s;t)) = 1_a(*) \otimes T(s;t) \subset X(a;a) \otimes T(s;t) .$$

The edges set mappings for a sections transport will be

$$(\lambda_X)_{a,b} : X(a;b) \to \prod_{s' \in T_0} X(a;b) \otimes 1_{s'}(*) \subset \prod_{s' \in T_0} (X(a;b) \otimes T(s';s')) .$$

The changing of target space with $V$-graphs transport $f : X \times T \Longrightarrow Y$ will appoint the $V$-graphs transport $f^T(\lambda_X(a)) : T \Longrightarrow Y$ with edges set mappings

$$f^T(\lambda_X(a))_{s',t'}(T(s';t')) = f_{\langle a,s\rangle,\langle a,t\rangle}(1_a(*) \otimes T(s;t)) \subset Y(f_0(a,s); f_0(a,t)) .$$

The edges set mappings of such changing target space will be

$$(f^T(\lambda_X))_{a,b}(X(a;b)) = \prod_{s' \in T_0} f_{\langle a,s'\rangle,\langle b,s'\rangle}(X(a;b) \otimes 1_{s'}(*)) \subset$$

$$\prod_{s' \in T_0} Y(f_0(a,s'); f_0(b,s')) .$$

The preevaluation transport for a set of transforms between two $V$-graphs transports $\phi, \psi : T \Longrightarrow Y$ and edges set $T(s;t)$ appoints the precontact arrows mapping

$$Y(\phi_0(s); \psi_0(s)) \circ \psi_{s,t}(T(s;t)) \subset Y(\phi_0(s); \psi_0(t)) \ .$$

So in our case we get the contact arrows mappings

$$f^+_{\langle a,s \rangle, \langle b,t \rangle}(X(a;b) \otimes T(s;t)) =$$

$$f_{\langle a,s \rangle, \langle b,s \rangle}(X(a;b) \otimes 1_s(*)) \circ f_{\langle b,s \rangle, \langle b,t \rangle}(1_b(*) \otimes T(s;t)) \ .$$

We can calculate its predecomposition

$$f^+_{\langle a,s \rangle, \langle b,s \rangle}(X(a,b) \otimes 1_s(*)) \circ f^+_{\langle b,s \rangle, \langle b,t \rangle}(1_b(*) \otimes T(s;t)) =$$

$$(f_{\langle a,s \rangle, \langle b,s \rangle}(X(a;b) \otimes 1_s(*)) \circ f_{\langle b,s \rangle, \langle b,s \rangle}(1_b(*) \otimes 1_s(*))) \circ$$

$$(f_{\langle a,s \rangle, \langle a,s \rangle}(1_t(*) \otimes 1_s(*)) \circ f_{\langle b,s \rangle, \langle b,t \rangle}(1_b(*) \otimes T(s;t))) =$$

$$f_{\langle a,s \rangle, \langle b,s \rangle}(X(a;b) \otimes 1_s(*)) \circ f_{\langle b,s \rangle, \langle b,t \rangle}(1_b(*) \otimes T(s;t)) =$$

$$f^+_{\langle a,s \rangle, \langle b,t \rangle}(X(a;b) \otimes T(s;t)) \subset Y(f_0(a,s); f_0(b,t)) \ .$$

$\square$

For a original $V$-graphs and their transports we can announce classical result about bijective name appointment.

For original $V$-graphs transport $f : X \otimes T \Longrightarrow Y$ with decomposition property

$$f_{\langle a,s \rangle, \langle b,s \rangle}(X(a;b) \otimes 1_s(*)) \circ f_{\langle b,s \rangle, \langle b,t \rangle}(1_b(*) \otimes T(s;t)) =$$

$$f_{\langle a,s \rangle, \langle b,t \rangle}(X(a;b) \otimes T(s;t)) \ .$$

we get its name corestricted in the exponential space of natural transforms $g : X \Longrightarrow Y^{(T)}$.

Otherwise for such partial transport of original $V$-graphs corestricted in exponential space of natural transforms we get realization with decomposition property. So name appointment is bijective

$$OGrph_V(X \otimes T; Y) = \to OGrph_V(X; Y^{(T)}) \ .$$

But we haven't naturality for transports in the target space $Y$.

We can generalize the notion of original $V$-graph $X$ to the notion of *lax original $V$-graph* demanding equalities for unit arrow mappings only up unique canonical isomorphism

$$1_a(*) \circ X(a;b) \sim X(a;b) \sim X(a;b) \circ 1_b(*) \ .$$

We get a new category $OGrphl_V$ compounded by all $V$-graphs transport $f : X \Longrightarrow Y$ which maintain the unit arrow mappings up unique canonical isomorphism in ordered space of mappings from the values category $V$

$$f_{a,a}(1_a(*)) \sim 1_{f_0(a)}(*) \subset Y(f_0(a); f_0(a)) \ .$$

It will be interesting in the case when such set of mappings $V(r;s)$ gets its ordering from another kind of arrows of some other category,

We can repeat all consideration about original $V$-graphs. We produce only last sentences. We shall call a $V$-graphs transport $f : X \otimes T \Longrightarrow Y$ *preregular* if the prerealization mapping reverses its name up unique canonical isomorphism

$$f \sim f_+(g(f)) \ .$$

For $V$-graphs transport $f$ with upper predecomposition we get an inequality with its name's prerealization, which will be preregular $V$-graphs transport

$$f \leq f_+(g(f)) = f^+ \ .$$

This upper bound is an initial in a categorical sense among preregular bounds.

Otherwise we shall call a $V$-graphs transport $f : X \otimes T \Longrightarrow Y$ *postregular* if the postrealization mapping reverses its name

$$f \sim f_-(g(f)) \ .$$

For $V$-graphs transport $f$ with the lower postdecomposition we get an inequality with its name's postrealization, which will be the postregular $V$-graphs transport

$$f \geq f_-(g(f)) = f^- \ .$$

This lower bound is a final in a categorical sense among postregular lower bounds.

We can generalize the space of natural transforms $Y^{(T)}$ demanding inverse image of similarity relation defined by existing unique canonical isomorphism. We shall denote such exponential space $Y_l^{(T)}$ and it will be called the exponential space of *laxly natural transforms*. Then taken $V$-graphs transport $f : X \otimes T \Longrightarrow Y$ has name corestricted in the space of natural transform exactly when it coincides up unique canonical isomorphism with its upper preregular and lower postregular bounds

$$f^- \sim f \sim f^+ \ .$$

Now we introduce the *category of $V$-enriched categories*. In an original $V$-graph $X$ we demand that for the contact arrows mappings would be valid associative law.

The first example will be the category $Cat_{Set}$ of $Set$-enriched categories. The $Set$-enriched category will be arbitrary category $X$ with small arrows sets, i. e. these sets will be objects in the category of small sets $Set$. So the properties of such categories will be familiar in the category **CAT** of all categories, only the provings will be adapted to general $V$-enriched categories. Values category $V$ will be taken a small complete monoidal category. The defined biproduct $r \otimes s$ must be a tensor product, i. e. we have natural *associativity isomorphism*

$$s_{r,s,v} : r \otimes (s \otimes v) =\!\!\rightarrow (r \otimes s) \otimes v$$

which meets the *pentagon diagram* from **S. MacLane 1971** part VII Monoids. At this moment we don't demand symmetry for such tensor product, only the *relator's involution*

$$m_{\langle p,s\rangle,\langle r,t\rangle} : (p\otimes s)\otimes(r\otimes t) \to (p\otimes r)\otimes(s\otimes t) .$$

Later we shall ask for these arrows some continuity properties to define a biproduct $X\otimes Y$ of $V$-enriched categories as a superfunctor in categories of continuous or cocontinuous functors $Cat_V^+$, $Cat_V^-$.

The language of enriched categories has loosed the capital place after we have introduced original enriched graphs for study of joint pairs of superfunctors. The associativity property of tensor product in values category $V$ is needed to define for enriched graphs $Y$ and $T$ the exponential space of natural transforms $Y^{(T)}$ as a superfunctor in the category of natural functors $Cat_V^=$. For the exponential space of continuous transforms $Y_+^{(T)}$ we can define such superfunctor only in the category of continuous-cocontinuous functors $Cat_V^\pm$. It will be proved for the tensor product in values category with natural associativity isomorphisms, also it will be true for the associativity property defined by continuous-cocontinuous arrows. To solve this problem for more weak associativity property would be a special task for others values categories.

As earlier we demand existence of an exponential functor $r^s$ with name appointments

$$V(r\otimes s;t) \to V(r;t^s) ,$$

but also we shall use the opposite realization appointments.

Arbitrary $V$-enriched category $X$ will be a set of points $X_0$ with objects in values category $X(x;y) \in V_0$ as *arrows sets* and arrows from values category $V$ as *mappings between arrows sets*.

The structure of $V$-enriched category is defined by *mappings of arrows sets composition*

$$X_{x,y,z} : X(x;y)\otimes X(y;z) \to X(x;z) .$$

We shall simply denote them with the sign of composition for the image set

$$X(x;y)\otimes X(y;z) \to X(x;y)\circ X(y;z) \subset X(x;z) .$$

Such mappings must ensure the associativity of composition

$$s_{X(x;y),X(y;z),X(z;u)} \circ (X_{x,y,z}\otimes 1_{X(z,u)}) \circ X_{x,z,u} =$$

$$(1_{X(x,y)} \otimes X_{y,z,u}) \circ X_{x,y,u} .$$

$$\begin{array}{ccc}
X(x;y)\otimes(X(y;z)\otimes X(z;u)) & \xrightarrow{s} & (X(x;y)\otimes X(y;z))\otimes X(z,;u) \\
{\scriptstyle 1_{X(x;y)}\otimes X_{y,z,u}}\downarrow & & \downarrow{\scriptstyle X_{x,y,z}\otimes 1_{X(z,u)}} \\
X(x;y)\otimes X(y;u) & & X(x;z)\otimes X(z;u) \\
{\scriptstyle X_{x,z,u}}\searrow & & \swarrow {\scriptstyle X_{x,z,u}} \\
& X(x;u) &
\end{array}$$

The choosing of unit arrow
$$X_x : * \to X(x;x)$$
can be denoted as
$$1_x : * \to 1_x(*) \subset X(x;x) \ .$$
The arrow from category's final object $* \to X(x;x)$ usually is understood as choosing of some point in its target space.

This unit arrow choosing must be a neutral for the source of arrows. Such property is expressed demanding the identity mapping for composition
$$1_{X(x;y)} : X(x,y) \xrightarrow{! \otimes 1_{X(x,y)}} * \otimes X(x;y) \xrightarrow{X_x \otimes 1_{X(x;y)}} X(x,x) \otimes X(x;y) \xrightarrow{X_{x,x,y}} X(x,y) \ .$$

The unit arrow choosing also must be neutral for the target of arrows. This property is expressed by demanding the identity mapping for another composition
$$1_{X(x,y)} : X(x;y) \xrightarrow{1_{X(x,y)} \otimes !} X(x;y) \otimes * \xrightarrow{1_{X(x,y)} \otimes X_y} X(x,y) \otimes X(y;y) \xrightarrow{X_{x,y,y}} X(x,y) \ .$$

The *natural functors between $V$-enriched categories* $f : X \Longrightarrow Y$ will be defined by appointment of category's points $f_0 : X_0 \to Y_0$ and mappings for the arrows sets
$$f_{x,y} : X(x,y) \to Y(f(x);f(y)) \ .$$
For such mappings we shall use also more detailed notation
$$f_{x,y} = X(x;y) \to f_{x,y}(X(x;y)) \subset Y(f(x);f(y)) \ .$$

These arrows sets mappings must allow the coherence properties. For choosing of unit arrow we have commuting diagram
$$\begin{array}{ccc} * & \xrightarrow{X_x} & X(x;x) \\ Y_{f_0(x)} \downarrow & \swarrow & f_{x,x} \\ Y(f_0(x);f_0(x)) & & \end{array} \ .$$

Also we have commuting diagrams for arrows composition mappings
$$\begin{array}{ccc} X(x,y) \otimes X(y;z) & \xrightarrow{X_{x,y,z}} & X(x,z) \\ f_{x,y} \otimes f_{y,z} \downarrow & & \downarrow f_{x,z} \\ Y(f_0(x);f_0(y)) \otimes Y(f_0(y);f_0(z)) & \xrightarrow[Y_{f_0(x),f_0(y),f_0(z)}]{} & Y(f_0(x);f_0(z)) \end{array} \ ,$$

The natural functors between $V$-enriched categories compound a category. The composition of two functors $f : X \Longrightarrow Y$ and $g : Y \Longrightarrow Z$ is defined as usual composition of points appointments
$$(f \circ g)_0(x) = g_0(f_0(x)) \uparrow x \in X$$

and the composition of mappings for arrows sets will be

$$(f \circ g)_{x,y} = f_{x,y} \circ g_{f_0(x),f_0(y)} : X(x;y) \to Z(g_0(f_0(x));g_0(f_0(y))) \ .$$

Obviously such composition is associative.

The identity functors $\text{Id}_X : X \Longrightarrow X$ are defined with identity appointment of points $\text{Id}_X(x) = x \uparrow x \in X$ and the identity mappings for arrows sets $\text{Id}_{x,y} : X(x;y) \to X(x;y)$. They will provide the *left and right units* for the composition of functors between $V$-enriched categories.

We shall call such category as a *category of natural functors between $V$-enriched categories* and it will be denoted $Cat_V^{\overline{=}}$. So we have embedding of such category in the category of natural functors between original $V$-graphs

$$Cat_V^{\overline{=}} \subset OGrph_V^{\overline{=}} \ .$$

For values category $V$ with ordered arrows sets $V(p;r)$ we can define a more larger stock of continuous functors. For them we get an imbedding in the category of continuous functors between original $V$-graphs

$$Cat_V^{-} \subset OGrph_V^{-} \ .$$

Also we can define the category of cocontinuous functors. For them we get an imbedding in the category of original $V$-graphs with cocontinuous functors

$$Cat_V^{+} \subset OGrph_V^{+} \ .$$

We again shall work with functors between constructed new categories. Such functors will be called *superfunctors*, as they are defined over larger category.

The superfunctor within the category $Cat_V$ will be defined by mapping $F$ over the set of all functors between $V$-enriched categories. It maps arbitrary $V$-enriched category $X$ to another $V$-enriched category $X^F$ and arbitrary functor $f : X \Longrightarrow Y$ of $V$-enriched categories to another functor $f^F : X^F \Longrightarrow Y^F$ between $V$-enriched categories.

It must maintain the composition of functors

$$(f \circ g)^F = f^F \circ g^F \ ,$$

and the appointment of unit arrows

$$(\text{Id}_X)^F = \text{Id}_{X^F} \ .$$

Now we shall define some concrete superfunctors. For two $V$-enriched categories $X$ and $Y$ we shall define a biproduct $X \otimes Y$ as a category with Carte tensor product of objects sets $X \times Y$ and biproduct in value category $V$ of arrows sets

$$X \otimes Y(\langle x,y\rangle;\langle x',y'\rangle) := X(x;x') \otimes Y(y;y') \ .$$

Such biproduct has defined superfunctor in the categories of functors between original graphs. For $V$-enriched categories we need additionally to check that in the biproduct $X \otimes Y$ contact arrows mappings

$$X \otimes Y(\langle a, a' \rangle, \langle b, b' \rangle) \circ X \otimes Y(\langle b, b' \rangle' \langle c, c' \rangle) =$$

$$(X(a;b) \circ X(b;c)) \otimes (Y(a';b') \circ Y(b';c')) \subset X(a;c) \otimes Y(a';c')$$

are associative

$$(X(a;b) \circ (X(b;c) \circ X(c;d))) \otimes (Y(a';b') \circ (Y(b';c') \circ Y(c';d'))) =$$

$$((X(a;b) \circ X(b;c)) \circ X(c;d)) \otimes ((Y(a';b') \circ Y(b';c')) \circ Y(c';d')) \ .$$

This is valid for each component, and we must remind that relator's involution are isomorphic.

In the category of small sets $Set$ we have used the Carte tensor product. We can check that new tensor product for $Set$-enriched categories rests Carte tensor product, i. e. its projections $p_1, p_2 : X \otimes Y \Longrightarrow X$ defines *final cone*. Such projections obviously are natural functors of $Set$-enriched categories, i. e. maintain the composition of arrows sets and unit arrows. For arbitrary other cone of two natural functors $f : Z \Longrightarrow X$ and $g : Z \Longrightarrow Y$ we unanimously get a natural functor of $Set$-enriched categories $\langle f \otimes g \rangle : Z \Longrightarrow X \otimes Y_0$. We get the unique appointment for every point $z \in Z$

$$\langle f_0(z), g_0(z) \rangle \in X_0 \times Y_0 \ .$$

For the mappings of arrows sets

$$f_{z_1, z_2} : Z(z_1; z_2) \to X(f_0(z_1); f_0(z_2)) \ , \quad g_{z_1; z_2} : Z(z_1; z_2) \to Y(g_0(z_1); g_0(z_2))$$

we unanimously get a mapping

$$\langle f_{z_1, z_2} \times g_{z_1.z_2} \rangle : Z(z_1; z_2) \to X(f_0(z_1); f_0(z_2)) \times Y(g_0(z_1); g_0(z_2)) \ .$$

Such mappings of arrows sets maintain the arrows sets composition

$$f_{z_1, z_2} \times g_{z_1.z_2} : Z(z_1; z_2) \circ Z(z_2; z_3) \to$$

$$(f_{z_1, z_2}(Z(z_1; z_2)) \times g_{z_1, z_2}(Z(z_1; z_2))) \times (f_{z_2, z_3}(Z(z_2; z_3)) \times g_{z_2, z_3}(Z(z_2; z_3))) =$$

$$(f_{z_1, z_2}(Z(z_1; z_2)) \times f_{z_2, z_3}(Z(z_2; z_3)) \times (g_{z_1, z_2}(Z(z_1; z_2)) \times g_{z_2, z_3}(Z(z_2; z_3))) \to$$

$$(f_{z_1, z_2}(Z(z_1; z_2)) \circ f_{z_2, z_3}(Z(z_2; z_3))) \times (g_{z_1, z_2}(Z(z_1; z_2)) \circ g_{z_2, z_3}(Z(z_2; z_3))) =$$

$$f_{z_1, z_3}(Z(z_1; z_2) \circ Z(z_2; z_3)) \times (g_{z_1, z_3}(Z(z_1; z_2) \circ g_{z_2, z_3}(Z(z_2; z_3)) \subset$$

$$f_{z_1, z_3}(Z(z_1; z_3)) \times g_{z_1, z_3}(Z(z_1; z_3)) \ .$$

and also maintain the equalities for chosen unit arrows

$$1_z(*) \circ f_{z,z} = 1_{f_0(z)}(*) \ ,$$

$$1_z(*) \circ g_{z,z} = 1_{g_0(z)}(*) \ .$$

We get the equality

$$1_z(*) \circ \langle f_{z,z} \times g_{z,z} \rangle) = 1_{\langle f_0(z), g_0(z) \rangle}(*) \ .$$

This ends the proving that defined tensor product of $Set$-enriched categories is a Carte tensor product in the category $Cat_{Set}$.

For the general values category $V$ we also can take the Carte tensor product $p \times r$, but another biproducts. can be also useful. Now we want to define the coadjoint exponential superfunctor $Y^T$ at first in the category of $Set$-enriched categories, i. e. for the categories with small arrows sets. This will be possible only for a small assistant category $T$.

We shall say that $Set$-enriched category $T$ is a *small category* if its points set $T_0$ is a small set. For such small set $T_0$ we shall have a products of small sets $S_s \in Set$

$$\prod_{s \in T_0} S_s \in Set \ ,$$

as for small set $T_0$ such product remain a small set, cl. **S. MacLane 1971** part 1, §6, exercise 1. We shall say that the category of small sets $Set$ is a *complete for small products*.

The exponential category of natural transforms $Y^{(T)}$ will be defined as some $Set$-enriched category only by means of *structural mappings* $Y_{x,y,z}$ and $Y_x$. For a $Set$-enriched category $Y$ and a small assistant $Set$-enriched category $T$ the exponential category $Y^{(T)}$ is defined with the objects set compounded of all functors between $Set$-enriched categories $\phi : T \Longrightarrow Y$. The arrows set between two functors $\phi : T \Longrightarrow Y$ and $\psi : T \Longrightarrow Y$ will be a set of natural transforms defined as an equalizer in the product of arrows sets

$$\to Y^{(T)}(\phi;\psi) \to \prod_{s \in T_0} Y(\phi_0(s);\psi_0(s)) \Rightarrow \prod_{s,t \in T_0} Y(\phi_0(s);\psi_0(t))^{T(s;t)} \ .$$

We shall use the property of truly joint functors $r \times s$ and $t^s$ in values category $Set$ to define such arrows. In values category $Set$ for the mapping

$$f : Y(\phi_0(s);\psi_0(s)) \otimes T(s;t) \to Y(\phi_0(s);\psi_0(t))$$

we get its name

$$g = \lceil f \rceil : Y(\phi_0(s);\psi_0(s)) \to Y(\phi_0(s);\psi_0(t))^{T(s;t)} \ .$$

This allows us to construct two different mappings over the product from values category $Set$

$$\prod_{s \in T_0} Y(\phi_0(s);\psi_0(s)) \ .$$

First we shall take an arrow set from source category $T(s;t)$. One mapping is defined by precontact arrows mapping

$$Y(\phi_0(s);\psi_0(s)) \otimes T(s;t) \to Y(\phi_0(s);\psi_0(s)) \circ \psi_{s,t}(T(s;t)) \subset Y(\phi_0(s);\psi_0(t)) \ ,$$

and another mapping is defined by another postcontact arrows mapping

$$Y(\phi_0(t); \psi_0(t)) \otimes T(s;t) \to \phi_{s,t}(T(s;t)) \circ Y(\phi_0(t); \psi_0(t)) \subset Y(\phi_0(s); \psi_0(t)) \ .$$

So we have got two mappings over the product of two arrows sets

$$Y(\phi_0(s); \psi_0(s)) \times Y(\phi_0(t); \psi_0(t)) \Rightarrow (Y(\phi_0(s); \psi_0(t)))^{T(s;t)} \ ,$$

yet more we get two needed arrows over the whole product

$$\prod_{s \in T_0} Y(\phi_0(s); \psi_0(s)) \Rightarrow Y(\phi_0(s); \psi_0(t))^{T(s;t)} \ ,$$

and they all together produces the correspondent arrows to the product of target spaces

$$\prod_{s \in T_0} Y(\phi_0(s); \psi_0(s)) \Rightarrow \prod_{s,t \in T_0} Y(\phi_0(s); \psi_0(t))^{T(s;t)} \ .$$

The notion of *small category* is meaningful also for arbitrary $V$-enriched category. We demand that the set of points $T_0$ for the category $T$ would be a small set. Let values category $V$ is *small complete*, i. e. we always have the small products of its objects. Then for a small assistant category $T$ we can apply the same construction of various exponential categories. For values category $V$ with ordered arrows sets $V(s;r)$ we are yet defined the set of continuous transforms $Y_{-}^{(T))}(\phi; \psi)$ or cocontinuous transforms $Y_{+}^{(T)}(\phi; \psi)$ in the category of original $V$-graphs $Grph_V$. Now we intend to define the correspondent superfunctors in the category of $V$-enriched categories $Cat_V$. We shall probe such construction first for the set of continuous transforms $Y_{-}^{(T)}(\phi; \psi)$.

We have an exponential superfunctor $Y^T$ defined by the sets of all transforms $Y^T(\phi; \psi)$ in the category of original $V$-graphs transports $Grph_V$. The associative composition mappings in the target $V$-graph $Y$ will define the associative composition mappings in the sets of all transforms $Y^T$. The equality of mappings

$$Y(\phi_0(s); \psi_0(s)) \circ (Y(\psi_0(s); \xi_0(s)) \circ Y(\xi_0(s); \chi_0(s))) =$$

$$(Y(\phi_0(s); \psi_0(s)) \circ Y(\psi_0(s); \xi_0(s))) \circ Y(\xi_0(s); \chi_0(s))$$

produces the equality for their products

$$Y^T(\phi; \psi) \circ (Y^T(\psi; \xi) \circ Y^T(\xi; \chi)) = (Y^T(\phi; \psi) \circ Y^T(\psi; \xi)) \circ Y^T(\xi; \chi) \ .$$

We need to check that for the sets of continuous transforms such composition produces a mapping to the correspondent set of continuous transforms, i. e. lies in the inverse image of order relation from correspondent product

$$\prod_{s,t \in T_0} Y(\phi_0(s); \xi_0(t))^{T(s;t)} \ .$$

Having the inequalities

$$Y(\phi_0(s); \psi_0(t)) \circ \psi_{s,t}(T(s;t)) \leq \phi_{s,t}(T(s;t)) \circ Y(\phi_0(t); \psi_0(t)) \ ,$$

$$Y(\psi_0(s); \xi_0(t)) \circ \xi_{s,t}(T(s;t)) \leq \psi_{s,t}(T(s;t)) \circ Y(\psi_0(t); \xi_0(t)) \ ,$$

we shall check the inequality for composition

$$(Y(\phi(s); \psi(s)) \circ Y(\psi(s); \xi(s))) \circ \xi(T(s;t)) \leq$$
$$\phi(T(s;t)) \circ (Y(\phi(t); \psi(t)) \circ Y(\psi(t); \xi(t))) \ .$$

It can be checked by steps applying an associative composition mappings in the target space $Y$

$$(Y(\phi_0(s); \psi_0(s)) \circ Y(\psi_0(s); \xi_0(s))) \circ \xi(T(s;t)) \leftarrow=$$
$$Y(\phi_0(s); \psi_0(s)) \circ (Y(\psi_0(s); \xi_0(s)) \circ \xi(T(s;t))) \leq$$
$$Y(\phi_0(s); \psi_0(s)) \circ (\psi(T(s;t)) \circ Y(\psi_0(t); \xi_0(t))) =\to$$
$$(Y(\phi_0(s); \psi_0(s)) \circ \psi(T(s;t))) \circ Y(\psi_0(t); \xi_0(t)) \leq$$
$$(\phi(T(s;t)) \circ Y(\phi_0(t); \psi_0(t))) \circ Y(\psi_0(t); \xi_0(t)) \leftarrow=$$
$$\phi(T(s;t)) \circ (Y(\phi_0(t); \psi_0(t)) \circ Y(\psi_0(t); \xi_0(t))) \ .$$

So we have checked that the set of continuous transforms $Y_-^{(T)}$ is again a $V$-enriched category. This proving is valid also for an equalizer defining the set of natural transforms $Y^{(T)}$, as the equality is only a partial case of an inequality.

We shall check when for a small assistant $V$-enriched category $T$ the appointment of new narrower exponential category $Y_-^{(T)}$ could define a superfunctor within some category of $V$-enriched categories. For the $V$-graphs transport $f : Y \Longrightarrow Y'$ we have appointed the changing of the target space $f^T : Y^T \Longrightarrow (Y')^T$ between functor spaces $Y^T$ and $(Y')^T$.

However we can't solve this problem in the categories $Cat_V^-$ or $Cat_V^+$. We shall work in a smaller category of $V$-enriched categories of functors which are continuous and cocontinuous at once

$$Cat_V^\mp = Cat_V^+ \cap Cat_V^- \ .$$

We also demand that arrows $m_{\langle p,s \rangle, \langle r,t \rangle}$ would define the transform being continuous and cocontinuous at once. Such property is weaker than the property of natural functors defined by equality relation, i. e. we have inclusion of categories

$$Cat_V^= \subset Cat_V^\mp \ .$$

We shall prove that the changing of target space with *continuous-cocontinuous functors* $f \in Cat_V^\pm$ maintains the exponential space of continuous transforms $Y_-^{(T)}$.

**Proposition 7.** *The changing of target space with continuous-cocontinuous functors $f \in Cat_V^{\pm}$ maintains the exponential space of continuous transforms $Y_-^{(T)}$.*

**Proof:**

Let we have inequalities for arbitrary V-graphs transports $\phi, \psi \in X^T$

$$X(\phi_0(s); \psi_0(s)) \circ \psi_{s,t}(T(s;t)) \leq \phi_{s,t}(T(s;t)) \circ X(\phi_0(t); \psi_0(t)) .$$

Then we can get such inequalities again for changed V-graphs transports $\phi \circ f, \psi \circ f \in Y^T$

$$f_{\phi_0(s), \psi_0(s)}(X(\phi_0(s); \psi_0(s))) \circ f_{\psi_0(s), \psi_0(t)}(\psi_{s,t}(T(s;t))) \leq$$

$$f_{\phi_0(s), \psi_0(t)}(X(\phi_0(s); \psi_0(s)) \circ \psi_{s,t}(T(s;t))) \leq$$

$$f_{\phi_0(s), \psi_0(t)}(\phi_{s,t}(T(s;t)) \circ X(\phi_0(t); \psi_0(t))) \leq$$

$$f_{\phi_0(s), \phi_0(t)}(\phi_{s,t}(T(s;t)))) \circ f_{\phi_0(t), \psi_0(t)}(X(\phi_0(t); \psi_0(t))) .$$

So the changed set of continuous transforms $Y_-^{(T)}(\phi; \psi)$ is included in the inverse image defining the set of continuous transforms between changed V-graphs transports $Y_-^{(T)}(\phi \circ f; \psi \circ f)$ .

To show the proving is easier for the case of equality, defining exponential space of natural transforms as an equalizer

$$\rightarrow X^{(T)}(\phi; \psi) \rightarrow \prod_{s \in T_0} X(\phi_0(s); \psi_0(s)) \rightrightarrows \prod_{s,t \in T_0} X(\phi_0(s); \psi_0(t))^{T(s;t)} .$$

The changing of target space defines mapping

$$\prod_{s \in T_0} X(\phi_0(s); \psi_0(s)) \rightarrow \prod_{s \in T_0} f_{\phi_0(s), \psi_0(s)}(X(\phi_0(s); \psi_0(s))) \subset$$

$$\prod_{s \in T_0} Y(f_0(\phi_0(s)); f_0(\psi_0(t))) .$$

The inclusion of $X^{(T)}(\phi; \psi)$ into the equalizer of name appointments provides equals realizations

$$X^{(T)}(\phi; \psi) \otimes T(s;t) \rightarrow X(\phi_0(s); \psi_0(s)) \circ \psi_{s,t}(T(s;t)) \subset X(\phi_0(s); \phi_0(t)) ,$$

$$X^{(T)}(\phi; \psi) \otimes T(s;t) \rightarrow \phi_{st}(T(s;t)) \circ X(\phi_0(t); \psi_0(t)) \subset X(\phi_0(s); \phi_0(t)) .$$

We have yet checked that these equalities provide the equalities also for the sets of transforms between changed V-graphs transports $\phi \circ f, \psi \circ f : T \Longrightarrow Y$

$$X^{(T)}(\phi; \psi) \otimes T(s;t) \rightarrow f_{\phi_0(s), \psi_0(s)}(X(\phi_0(s); \psi_0(s))) \circ f_{\psi_0(s), \psi_0(t)}(\psi_{s,t}(T(s;t))) \subset$$

$$Y(f_0(\phi_0(s)); f_0(\psi_0(t))) ,$$

$$X^{(T)}(\phi;\psi) \otimes T(s;t) \to f_{\phi_0(s),\phi_0(t)}(\phi_{st}(T(s;t))) \circ f_{\phi_0(t),\psi_0(t)}(X(\phi_0(t);\psi_0(t))) \subset$$
$$Y(f_0(\phi_0(s)); f_0(\psi_0(t))) ,$$

So we get an equal name arrows, and an arrow to their equalizer

$$\begin{array}{ccc} X^{(T)}(\phi;\psi) & \to & \prod_{s\in T_0} X(\phi_0(s);\psi_0(s)) & \Rightarrow \\ \downarrow & & \downarrow & \\ Y^{(T)}(\phi \circ f; \psi \circ f) & \to & \prod_{s\in T_0} Y(f_0(\phi_0(s)); f_0(\psi_0(s))) & \Rightarrow \end{array}$$

For the exponential space of continuous transforms $Y_-^{(T)}$ we need work with product mappings between product spaces

$$X_-^{(T)}(\phi;\psi) \times X_-^{(T)}(\phi;\psi) \to \prod_{s\in T_0} X(\phi_0(s);\psi_0(s)) \times \prod_{s\in T_0} X(\phi_0(s);\psi_0(s))$$

$$\to \prod_{s,t\in T_0} X(\phi_0(s);\psi_0(t))^{T(s;t)} \times \prod_{s,t\in T_0} X(\phi_0(s);\psi_0(t))^{T(s;t)} .$$

to define inverse image of order relation in the last product

$$(\leq) \subset \prod_{s,t\in T_0} X(\phi_0(s);\psi_0(t))^{T(s;t)} \times \prod_{s,t\in T_0} X(\phi_0(s);\psi_0(t))^{T(s;t)} .$$

□

With a small $V$-enriched assistant category $T$ we have constructed two superfunctors in the category of $V$-enriched categories: the biproduct $X \otimes T$ and the *exponential space of continuous transforms between continuous-cocontinuous functors* $Y_-^{(T)\mp}$. We want to check have we got truly joint pair of superfunctors

$$\langle F, G \rangle : Cat_V^{\mp} \rightleftarrows Cat_V^{\mp} .$$

The unit transform for such joint pair of superfunctors is defined by collection of *sections functors* of $V$-enriched categories

$$\lambda_X : X \Longrightarrow (X \otimes T)_{\underline{=}}^{(T)=} \subset (X \otimes T)_-^{(T)\pm} .$$

Also the unit arrow mapping is corestricted in the set of naturals transforms

$$1_{\lambda_X(x)}(*) = \prod_{s\in T_0} (1_x(*) \otimes 1_s(*)) \subset (X \otimes T)_{\underline{=}}^{(T)=}(\lambda_X(x); \lambda_X(x)) .$$

So indeed we have defined natural functor $\lambda_X : X \Longrightarrow (X \otimes T)_{\underline{=}}^{(T)=}$.

These sections functors $\lambda_X : X \to (X \times T)^T$ has defined a natural transform between superfunctors in the category of original $V$-graphs $OGrph_V$, so the commuting quadrats remains in imbedded smaller category of continuous-cocontinuous functors between $V$-enriched categories $Cat_V^{\mp}$ or the category of natural functors between $V$-enriched categories $Cat_V^=$, i. e. for a continuous-cocontinuous functor between $V$-enriched categories $f : X \to X'$ we get commutating quadrats

$$\begin{array}{ccc} X & \xrightarrow{f} & X' \\ \lambda_X \downarrow & & \downarrow \lambda_{X'} \\ (X \otimes T)_{\mp}^{(T)\mp} & \xrightarrow[(f \otimes 1_T)^{(T)\mp}]{} & (X' \otimes T)_{\mp}^{(T)\mp} \end{array}$$

In the category of original $V$-graphs $Grph_V$ for the joint pair of superfunctors defined with the biproduct of assistant $V$-graph $X \otimes T$ and exponential space of all transforms between original $V$-graphs transports $Y^T$ the counit transform was provided in two ways. One is a collection of preevaluation functors

$$\mathbf{ev}_Y^+ : Y^T \otimes T \Longrightarrow Y .$$

They for a original $V$-graphs transport $\phi : T \Longrightarrow Y$ and vertex $s \in T$ appoints the value $\phi_0(s) \in Y$ and for the biproduct of edges sets $Y^T(\phi; \psi) \otimes T(a; b)$ they appoint the set of *precontact arrow*

$$Y(\phi_0(s); \psi_0(s)) \circ \psi_{s,t}(T(s;t)) \subset Y(\phi_0(s); \psi_0(t)) .$$

Another is a collection of postevaluation functors

$$\mathbf{ev}_Y^- : Y^T \otimes T \Longrightarrow Y .$$

They for a original $V$-graphs transport $\phi : T \Longrightarrow Y$ and vertex $s \in T$ appoints the same value $\phi_0(s) \in Y$ and for the biproduct of edges sets $Y^T(\phi; \psi) \otimes T(a;b)$ they appoint the set of *postcontact arrow*

$$\phi_{s,t}(T(s;t)) \circ Y(\phi_0(t); \psi(_0t)) \subset Y(\phi_0(s); \psi_0(t)) .$$

Then with correspondent prerealization or postrealization appointments of name from any transport $f : X \otimes T \Longrightarrow Y$ with the properties of upper predecomposition and lower postdecomposition

$$f_{\langle a,s \rangle, \langle b,s \rangle}(X(a;b) \otimes 1_s((*)) \circ f_{\langle b,s \rangle, \langle b,t \rangle}(1_b(*) \otimes T(s;t)) \geq$$

$$f_{\langle a,s \rangle, \langle b,t \rangle}(X(a;b) \otimes T(s;t)) \geq$$

$$f_{\langle a,s \rangle, \langle a,t \rangle}(1_a(*) \otimes T(s;t)) \circ f_{\langle a,t \rangle, \langle b,t \rangle}(X(a;b) \otimes 1_t(*))$$

we have got the inequalities

$$f^+ \geq f \geq f^- .$$

If contact arrows mappings in $V$-enriched category $Y$ were monotone, then these inequalities can be interpreted as unit and counit transforms for a newly defined joint pair of monotone superfunctors.

For the categories of $V$-enriched categories we shall construct more usual truly joint pair of superfunctors. In the category compounded by continuous functors between $V$-enriched categories $Cat_V^-$ we can check the continuity of preevaluation transport over the exponential spaces of continuous transforms between continuous functors

$$\mathbf{ev}_Y^+ : Y_-^{(T)-} \otimes T \Longrightarrow Y .$$

. Let we have three continuous transports of $V$-enriched graphs $\phi, \psi, \xi : T \Longrightarrow Y$ with the correspondent sets of continuous transforms $Y_-^{(T)}(\phi; \psi)$, $Y_-^{(T)}(\psi; \xi)$, $Y_-^{(T)}(\phi; \xi)$. For them we can draw the diagram

$$\begin{array}{ccccc} \phi(s) & & & & \\ \downarrow & \searrow & & & \\ \psi(s) & \to & \psi(t) & & \\ \downarrow & & \downarrow & \searrow & \\ \xi(s) & \to & \xi(t) & \to & \xi(v) \end{array}$$

Then we can write some equalities with associativity isomorphisms and inequalities for the sets of continuous transforms between continuous functors

$$(Y_-^{(T)}(\phi;\psi)_s \circ \psi_{s,t}(T(s;t))) \circ (Y_-^{(T)}(\psi;\xi)_t \circ \xi_{t,v}(T(t;v))) \longleftarrow=$$

$$Y_-^{(T)}(\phi;\psi)_s \circ (\psi_{s,t}(T(s;t)) \circ (Y_-^{(T)}(\psi;\xi)_t \circ \xi_{t,v}(T(t;v)))) =\to=$$

$$Y_-^{(T)}(\phi;\psi)_s \circ ((\psi_{s,t}(T(s;t)) \circ Y_-^{(T)}(\psi;\xi)_t) \circ \xi_{t,v}(T(t;v))) \leq$$

$$Y_-^{(T)}(\phi;\psi)_s \circ ((Y_-^{(T)}(\psi;\xi))_s \circ \xi_{s,t}(T(s;t))) \circ \xi_{(t,v}(T(t;v))) \longleftarrow=$$

$$Y_-^{(T)}(\phi;\psi)_s \circ (Y_-^{(T)}(\psi;\xi)_s) \circ (\xi_{s,t}(T(s;t)) \circ \xi_{t,v}(T(t;v)))) =\to$$

$$(Y_-^{(T)}(\phi;\psi)_s \circ (Y_-^{(T)}(\psi;\xi)_s) \circ (\xi_{s,t}(T(s;t) \circ \xi_{t,v}(T(t;v))) .$$

Otherwise in the category compounded by cocontinuous factors between $V$-enriched categories $Cat_V^+$ we can check the cocontinuity of postevaluation transport $\mathbf{ev}_Y^- : Y_+^{(T)+} \otimes T \Longrightarrow T$

$$(\phi_{s,t}(T(s;t)) \circ Y_+^{(T)}(\phi;\psi)_t) \circ (\psi_{t,v}(T(t;v)) \circ Y_+^{(T)}(\psi;\xi)_v) =$$

$$\phi_{s,t}(T(s;t)) \circ ((Y_+^{(T)}(\phi;\psi)_t) \circ (\psi_{t,v}(T(t;v))) \circ Y_+^{(T)}(\psi;\xi)_v) \geq$$

$$\phi_{s,t}(T(s;t)) \circ ((\phi_{t,v}(T(t;v)) \circ Y_+^{(T)}(\phi;\psi)_v) \circ Y_+^{(T)}(\psi;\xi)_v) =$$

$$(\phi_{s,t}(T(s;t)) \circ \phi_{t,v}(T(t;v))) \circ (Y_+^{(T)}(\phi;\psi)_v \circ Y_+^{(T)}(\psi;\xi)_v) .$$

At this time I haven't written used isomorphic associativity arrows.

So in the category of continuous-cocontinuous functors between $V$-enriched categories $Cat_\mp^{(T)\mp}$ we get both inequalities between preevaluation and postevaluation transports over the exponential category of continuous-cocontinuous transforms between continuous-cocontinuous functors of $V$-enriched categories. We get two continuous-cocontinuous evaluation transports

$$\mathbf{ev}_Y^\pm : Y_\mp^{(T)\mp} \otimes T \Longrightarrow Y \ .$$

For them we get *double inequalities*

$$f^+ \geq\leq f \geq\leq f^- \ .$$

Otherwise in the category of natural functors between $V$-enriched categories $Cat_V^=$ we get natural evaluation functor over exponential space of natural transforms between natural functors

$$\mathbf{ev}_Y : Y_=^{(T)=} \otimes T \Longrightarrow Y \ .$$

We want to check that such transforms are natural, i. e. we have commuting quadrats for arbitrary functor from correspondent categories between $V$-enriched categories $g : Y \Longrightarrow Y'$

$$\begin{array}{ccc} Y^T \otimes T & \stackrel{g^T \otimes T}{\longrightarrow} & Y'^T \otimes T \\ \mathbf{ev}_Y \downarrow & & \downarrow \mathbf{ev}_{(Y')} \\ Y & \underset{g}{\longrightarrow} & Y' \end{array}$$

We have shown such commuting diagrams even in the category of $V$-enriched graphs for the transports between $V$-enriched graphs $g : Y \Longrightarrow Y'$. Now for the continuous-cocontinuous $V$-graphs transports $g : Y \Longrightarrow Y'$ we have maintained the exponential space of continuous transforms

$$g_{\phi,\psi}^T(Y_-^{(T)\mp}(\phi;\psi)) \subset (Y')_-^{(T)\mp}(\phi \circ g; \psi \circ g) \ .$$

The same is true for the exponential space of cocontinuous transforms between continuous -cocontinuous $V$-graphs transports

$$g_{\phi,\psi}^T(Y_+^{(T)\mp}(\phi;\psi)) \subset (Y')_+^{(T)\mp}(\phi \circ g; \psi \circ g) \ ,$$

also for the intersection of such exponential spaces

$$g_{\phi,\psi}^T(Y_\mp^{(T)\mp})(\phi;\psi)) \subset (Y')_\mp^{(T)\mp}(\phi \circ g; \psi \circ g) \ ,$$

For the cocontinuous transform defined with arrows $M_{\langle p,s\rangle,\langle r,t\rangle}$ we get a continuous biproduct $g^T \otimes \mathtt{Id}_T : Y^T \otimes T \Longrightarrow (Y')^T \otimes T$ for the exponential spaces $Y^T$ of all transforms between arbitrary $V$-graphs transports. This is true also for smaller exponential spaces $Y_-^{(T)\mp}$ of continuous transforms between continuous-cocontinuous functors of $V$-enriched categories.

And for the cocontinuous transform defined with arrows $M_{\langle p,s\rangle,\langle r,t\rangle}$ we get cocontinuous biproduct $g^T \otimes \mathtt{Id}_T : Y^T \otimes T \Longrightarrow (Y')^T \otimes T$ for exponential spaces $Y^T$ of all transforms between arbitrary $V$-graphs transports. This is also true for smaller exponential spaces $Y_+^{(T)\mp}$ cocontinuous transforms between continuous-cocontinuous functors between $V$-enriched categories.

We shall get usual results for natural functors. Sometimes is useful to apply more complicate notation for the exponential space of natural transforms

$$Y^{(T)} = Y_{\underline{\underline{\phantom{x}}}}^{(T)} \ .$$

**Proposition 8.** *Let we have natural transform provided by relator's involution*

$$m_{\langle p,s \rangle, \langle r,t \rangle} : (p \otimes s) \otimes (r \otimes t) =\to (p \otimes r) \otimes (s \otimes t) \ .$$

*Then for the category of natural functors between $V$-enriched categories $Cat_V^{\underline{\underline{\phantom{x}}}}$ we get all properties of usual trey joint pair of superfunctors*

$$Cat_V^{\underline{\underline{\phantom{x}}}} \rightleftarrows Cat_V^{\underline{\underline{\phantom{x}}}} \ .$$

*The biproduct $X \otimes T$ defines a superfunctor.*

*For a natural functor $g : X \Longrightarrow Y$ it appoints the natural functor provided by biproduct with identity functor of assistant $V$- category $g \otimes \mathtt{Id}_T : X \otimes T \Longrightarrow Y \otimes T$.*

*The exponential space of natural transforms between natural functors $Y_{\underline{\underline{\phantom{x}}}}^{(T)=}$ defines coadjoint superfunctor.*

*For a natural functor between $V$-enriched categories $f : X \Longrightarrow X$ the changing of target source $f^T : X^T \Longrightarrow Y^T$ maintains the imbedded exponential space of natural transforms, so we get the transport of such exponential spaces $g_{\underline{\underline{\phantom{x}}}}^{(T)=} : X_{\underline{\underline{\phantom{x}}}}^{(T)=} \Longrightarrow Y_{\underline{\underline{\phantom{x}}}}^{(T)=}$.*

*The unit transform is defined by natural sections functors*

$$\lambda_X : X \Longrightarrow X_{\underline{\underline{\phantom{x}}}}^{(T)=} \ .$$

They are natural even for arbitrary original $V$-graphs transport $f : X \Longrightarrow Y$.

*The counit transform is defined by natural evaluation functor*

$$\mathbf{ev}_Y : Y_{\underline{\underline{\phantom{x}}}}^{(T)=} \otimes T \Longrightarrow Y \ .$$

It is also natural for natural functors of $V$-enriched categories $g : X \Longrightarrow Y$.

*The name appointment $Cat_V^{\underline{\underline{\phantom{x}}}}(X \otimes T; Y) \to Cat_V^+(X; Y_{\underline{\underline{\phantom{x}}}}^{(T)=})$. for a natural functor $f : X \otimes T \Longrightarrow Y$ appoints the composition*

$$g(f) = f^\# : X \xrightarrow{\lambda_X} (X \otimes T)^T \xrightarrow{f^T} Y^T \ .$$

*The realization appointment $Cat_V^{\underline{\underline{\phantom{x}}}}(X; Y_{\underline{\underline{\phantom{x}}}}^{(T)=}) \to Cat_V^{\underline{\underline{\phantom{x}}}}(X \otimes T; Y)$ for a natural functor $g : X \Longrightarrow Y_{\underline{\underline{\phantom{x}}}}^{(T)=}$ appoints the composition*

$$f(g) = g^\flat : X \otimes T \xrightarrow{g \otimes 1_T} Y_{\underline{\underline{\phantom{x}}}}^{(T)=} \otimes T \xrightarrow{\mathbf{ev}_Y} Y \ .$$

*The name and realization appointments are one to another inverse bijections, also natural for natural functors $f : X \Longrightarrow X'$ and $g : Y \Longrightarrow Y'$ between $V$-enriched categories.*

*So the first triangular equality is get as name of the identity mapping* $\text{Id}_{X \otimes T} : X \otimes T \Longrightarrow X \otimes T$, *and the second triangular equality is get as realization of the identity mapping* $\text{Id}_{Y_{\underline{\underline{\cong}}}^{(T)=}} : Y_{\underline{\underline{\cong}}}^{(T)=} \Longrightarrow Y_{\underline{\underline{\cong}}}^{(T)=}$.

**Proof:** This is only partial case of truly joint pair of superfunctors in the category $Cat_V^{\mp}$ of continuous-cocontinuous functors between $V$-enriched categories.

□

We can introduce lax categories demanding that in some lax original $V$-graph $X$ the contact arrows mappings would define the associative composition up to unique canonical isomorphism. Such assumption may be almost possible. May be we shall need some coherence theorems, as it was for tensor product in monoidal categories done by MacLane 1970.

Now we could work with lax functors which maintained contact arrows mappings up unique canonical isomorphisms. They will compound a new category of *lax functors* $Catl_V^{\overline{=}}$. Or we could work with *weakly lax functors* which maintain contact arrows nap pings up some isomorphism. Their category we shall denote $Catw_V$.

The biproduct with assistant lax category $X \otimes T$ again defines a superfunctor in the category of lax functors. For the lax functors $\phi, \psi : T \Longrightarrow Y$ we can define the set of *lax transforms* $Y_{l=}^{(T)l=}(\phi; \psi)$. Such sets defines superfunctor $Y_{l=}^{(T)l=}$ in the category of lax functors $Catl_V^{\overline{=}}$.

We again get a joint pair of superfunctors in the category of lax functors

$$Catl_V^{\overline{=}} \rightleftarrows Catl_V^{\overline{=}}$$

with unit transform provided by sections functors

$$\lambda_X : X \Longrightarrow (X \otimes)_{\underline{\underline{\cong}}}^{(T)=}$$

and two counit transforms provided by preevaluation or postevaluation functors

$$\mathbf{ev}_Y^p m : Y_{l=}^{(T)l=} \otimes T \Longrightarrow Y \ .$$

They both will be lax functors coinciding up Canonical isomorphism. For proving we need to repeat earlier proposition about continuous-cocontinuous functors between $V$-enriched categories. The name appointment can't be bijective. We shall recognize the equivalence between forms $f : X \otimes T \Longrightarrow Y$ having the same name form.

**Proposition 9.** *Let we have laxly natural transform provided by relator's involution*

$$m_{\langle p,s \rangle, \langle r,t \rangle} : (p \otimes s) \otimes (r \otimes t) =\to (p \otimes r) \otimes (s \otimes t) \ .$$

*Then in the category* $Catl_V^{\overline{=}}$ *of laxly natural functors between $V$-enriched categories we get the joint pair of superfunctors. One adjoint is the biproduct* $X \otimes T$. *For a laxly natural functor* $g : X \Longrightarrow Y$ *it appoints the laxly natural* functor provided by biproduct with identity functor of assistant $V$- category $g \otimes \text{Id}_T : X \otimes T \Longrightarrow Y \otimes T$. For a monotone biproduct from values category $p \otimes r$, we get this superfunctor $X \otimes T$ being monotone.

The exponential space of laxly natural transforms between laxly natural functors $Y_{l=}^{(T)l=}$ defines coadjoint superfunctor. For a laxly natural functor between $V$-enriched categories $f : X \Longrightarrow X$ the changing of target source $f^T : X^T \Longrightarrow Y^T$ maintains the imbedded exponential space of laxly natural transforms, so we get the transport of such exponential spaces $g_{l=}^{(T)l=} : X_{\underline{=}}^{(T)=} \Longrightarrow Y_{\underline{=}}^{(T)=}$. For the monotone contact arrows mappings in target space $Y$, we get again this superfunctor monotone one.

The unit transform is defined by natural sections functors

$$\lambda_X : X \Longrightarrow X_{\underline{=}}^{(T)=} \ .$$

They are natural even for arbitrary original $V$-graphs transport $f : X \Longrightarrow Y$.

We have two counit transform defined by evaluation functors

$$\mathbf{ev}_Y^{\pm} : Y_{l=}^{(T)l=} \otimes T \Longrightarrow Y \ .$$

It is also natural for laxly natural functors between $V$-enriched categories $g : X \Longrightarrow Y$.

The name appointment $Cat_V^=(X \otimes T; Y) \to Cat_V^+(X; Y_{\underline{=}}^{(T)=})$. is corestricted in the space $Y_{l=}^{(T)l=}$ of laxly natural transforms between laxly natural functors exactly when we have equalities up unique canonical isomorphism with both realizations defined by the preevaluation and the postevaluation functors

$$f^+ \sim f \sim f^- \ .$$

**Proof:** This is partial case of the case with continuous-cocontinuous functors and their continuous-cocontinuous transforms.

□

We shall try to recognize $V$-graphs for some concrete values categories. For the category of small sets $Set$ we get usual graphs. The $Set$-enriched graph will be arbitrary graph $X$ with small edges sets $X(a; b)$. So we get an original graphs or categories. The $Set$-graphs transport $f : X \Longrightarrow Y$ will be transport of usual graphs with edges set mappings

$$f_{a,b} : X(a; b) \to Y(f_0(a); f_0(b)) \ .$$

So all properties of $Set$-enriched graphs or categories can be verified for usual graphs or categories. For the functors of categories $\phi, \psi : T \Longrightarrow Y$ we get the same set of natural transforms $Y^{(T)}(\phi; \psi)$. These sets defines an exponential superfunctor $Y^{(T)}$ which is coadjoint to the superfunctor defined with biproduct $X \otimes T$.

More interesting opportunities are provided by a larger values category $Cat$. The points of such category will be categories. They will be denoted by hind Latin letters $p, q, r, s, t \in Cat_0$. Further the points in such categories will be denoted $u, v, x, y, z \in p_0$. In the category $Cat$ the arrows $f : p \to q$ will be functors between such categories. They will be denoted by fore Latin letters $f, g, h, k \in Cat$. The usual composition of functors defines the appointment

$$f \circ g \in Cat(p; q \uparrow \langle f, g \rangle \in Cat(p; q) \times Cat(q; r)$$

over *pullbacked product*

$$Cat \times_{Cat_0} Cat = \{\langle f, g \rangle : p_2(f) = p_1(g)\} \subset Cat \times Cat \ .$$

The identity functors are *neutral arrows* for the source and the target of arrows

$$\mathtt{Id}_p \circ f = f = f \circ \mathtt{Id}_q \uparrow f \in Cat(p; q) \ .$$

The arrows from $Cat(p; q)$ will be called parallel. For parallel arrows we can define order relation with existing natural transform between functors

$$\alpha : f \to g \uparrow f, g \in Cat(p; q) \ .$$

The natural transforms will be denoted by initial Greek letters $\alpha, \beta, \gamma, \theta$. We shall say $f \geq g$ if exists a natural transform $\alpha : f \to g$. For functors between categories it is possible to construct concrete natural transforms, but we shall work with notion of continuity which uses only order relation between parallel functors defined by existing natural transform.

We choose Carte tensor product of categories as a biproduct

$$p \otimes q := p \times q \ .$$

The final category with unique point $* = \{O\}$ will be a *neutral point* for such biproduct, i. e. we have uniquely defined canonical isomorphisms of biproducts

$$* \otimes r = r = r \otimes * \ .$$

We have an expression of functors to the product of categories

$$Cat(p; q \otimes r) \to Cat(p; q) \times Cat(p; r) \ .$$

The Carte product of natural transforms $\alpha : f \to g \uparrow f, g \in Cat(p; q)$ and $\beta : h \to k \uparrow h, k \in Cat(p; r)$ provides the natural transform between the functors to the Carte product of target spaces

$$\alpha \times \beta : f \times h \to g \times k \uparrow \langle f \times h \rangle, \langle g \times k \rangle \in Cat(p; q \times r) \ .$$

This shows that the biproduct of categories $p \times q$ defines monotone bisuperfunctor. For a natural transform $\alpha : f \to g$ between two functors $f, g \in Cat(p; q)$ we get the natural transform between biproducts of these functors with any functor $h \in Cat(r; s)$

$$\alpha \times 1_h : f \times h \to g \times h \ .$$

Also for a natural transform $\beta : h \to k$ between two functors $h, k \in Cat(r; s)$ we get the natural transform between biproducts of these functors with any functor $f \in Cat(p; q)$

$$1_f \times \beta : f \times h \to f \times k .$$

The same is true for arbitrary Carte product of categories

$$\prod_u r_u .$$

It defines monotone superfunctors with edges set mappings

$$\prod_u Cat(p; r_s) \Longrightarrow Cat(p; \prod_u r_u) .$$

The chosen biproduct $r \otimes s$ with an assistant category $s \in Cat_0$ defines a superfunctor within the $Cat$. For a functor $f \in Cat(p; r)$ we shall appoint the biproduct with identity functor $\text{Id}_s \in Cat(s; s)$

$$f \otimes \text{Id}_s \in Cat(p \otimes s; r \otimes s) .$$

Such superfunctor will be adjoint one. The coadjoint superfunctor is provided by an exponential space of natural transforms $p^s(\phi; \psi)$. The exponential space of natural transforms $p^s$ is a category on the set of parallel functors $Cat(s; p) = (p^s)_0$ having as arrows natural transforms $\alpha : \phi \to \psi$ between some functors $\phi; \psi \in Cat(s; p)$. In such category we take the *vertical composition* of natural transforms

$$\alpha \circ \beta : f \to h \uparrow \alpha : f \to g, \beta : g \to h .$$

The identity transforms $\text{Id}_f : f \to f$ will be *neutral units* for such composition

$$\text{Id}_f \circ \alpha = \alpha = \alpha \circ \text{Id}_g \uparrow \alpha : f \to g .$$

The exponential space $p^s$ indeed defines a superfunctor. It for a functor $f \in Cat(p; q)$ appoints the changing of target space $f^s \in Cat(p^s; q^s)$. As usually it maintains the identity functor and the composition of functors.

We have another *horizontal composition* of natural transforms

$$\alpha * \alpha' : f \circ h \to g \circ k \uparrow \alpha : f \to g, \alpha' : h \to k .$$

It will be a bifunctor over the exponential spaces of natural transforms

$$* \in Cat(q^p \times r^q; r^p) .$$

Indeed. We have *middle four exchange* equality for both vertical and horizontal compositions of natural transforms

$$(\alpha \circ \beta) * (\alpha' \circ \beta') = (\alpha * \alpha') \circ (\beta * \beta') .$$

Therefore we have commuting diagram

$$
\begin{array}{ccc}
(q^p \times r^q) \times (q^p \times r^p) & \xrightarrow{* \times *} & r^p \times r^p \\
\circ \downarrow & & \downarrow \circ \\
q^p \times r^q & \xrightarrow{*} & r^p
\end{array}
$$

That means that horizontal composition is graphs transport natural for vertical composition. Also it maintains identity transforms

$$\mathtt{Id}_f * \mathtt{Id}_h = \mathtt{Id}_{f \circ h} \ .$$

So horizontal composition is a natural functor over Carte tensor product of categories $q^p \times r^q$ with vertical composition of arrows.

We are ready to prove that superfunctor of exponential space $p^s$ is monotone for the weaker order relation in exponential space $r^s$ defined between *secondarily parallel functionals*.If we have an inequality between primarily parallel functors

$$f \geq g : p \to r$$

we get the inequality between secondarily parallel functionals $f^s, g^s : p^s \to r^s$. We have coinciding secondarily projections for changing of target space

$$(f_0)^{s_0} = (g_0)^{s_0} : (p_0)^{s_0} \to (r_0)^{s_0} \ ,$$

i. e. for arbitrary functor $\phi \in p^s$ we have the equality of composition

$$\phi_0 \circ f_0 = \phi_0 \circ g_0 \ .$$

Also we get the inequalities of arrows set mappings

$$(\phi \circ f)_{u,v} \geq (\phi \circ g)_{u,v} \uparrow u, v \in s \ ,$$

$$f^s_{\phi;\psi} \geq g^s_{\phi,\psi} \uparrow \phi, \psi \in p^s \ .$$

For a functor $\phi \in Cat(s;p)$ the changing of target space with functor $f \in Cat(p;r)$ appoints the new functor $\phi \circ f \in Cat(s;r)$, or the changing of target space with another functor $g \in Cat(p;r)$ appoints another new functor $\phi \circ g \in Cat(s;r)$. The horizontal composition $1_\phi * \alpha : \phi \circ f \to \phi \circ g$ provides needed natural transform for the first inequality. This can be understood as transferred natural transform by source functor $\phi \star \alpha$. For each point $u \in s_0$ it produces the arrow in the target category $r$

$$\alpha_{\phi_0(u)} \in r(f_0(\phi_0(u)); f_0(\psi_0(u)))$$

with commuting diagram for arrows in assistant category $\theta \in s(u;v)$ between the points $u, v \in s_0$

$$
\begin{array}{ccc}
f_0(\phi_0(u)) & \xrightarrow{\alpha_{\phi_0(u)}} & g_0(\phi_0(u)) \\
f(\phi(\theta)) \downarrow & & \downarrow f(g(\theta)) \\
f_0(\phi_0(v)) & \xrightarrow[\alpha_{\phi_0(v)}]{} & g_0(\phi_0(v))
\end{array}
$$

So the collection of these arrows $\alpha_{\phi_0(u)} \uparrow u \in s_0$ provides a natural transform $\phi \times \alpha_u : \phi \circ f \to \phi \circ g$ between two functors $\phi \circ f, \phi \circ g \in (r^s)_0$, i. e. an arrow in exponential space of natural transforms

$$\phi \star \alpha \in r^s(\phi \circ f; \phi \circ g) \ .$$

The collection of such arrows defines the transform between two functors defined by changing of target space $\alpha^s : f^s \to g^s$. We need to check that this transform is natural for the arrows in the category $p^s$, i. e. for arbitrary natural transforms $\theta : \phi \to \psi$ between two functors $\phi, \psi \in Cat(s;p)$ we have commuting diagram

$$\begin{array}{ccc} \phi \circ f & \xrightarrow{\phi \times \alpha} & \phi \circ g \\ \theta \star f \downarrow & & \downarrow \theta \star g \\ \psi \circ f & \xrightarrow[\psi \times \alpha]{} & \psi \circ g \end{array}$$

This can be proved by equalities

$$(\phi \star \alpha) \circ (\theta \star g) = (1_\phi \ast \alpha) \circ (\theta \ast 1_g) =$$

$$(1_\phi \circ \theta) \ast (\alpha \circ 1_g) = \theta \ast \alpha \ .$$

For each point $x \in p_0$ in the category $p$ we have chosen a unit arrow $1_x \in p(x;x)$, so we can define the sections functor $\lambda_p \in Cat(p;(p \otimes s)^s)$ and the name appointment

$$Cat(p \otimes s; r) \to Cat(p; r^s)$$

provided by composition of functors $f^\# = \lambda_p \circ f^s$.

We can check that the name appointment remains monotone for the weaker order relation in exponential space $r^s$ defined with natural transforms between secondarily parallel functors.

For a point $x \in p_0$ we get the section $\lambda_p(x) \in (p \otimes s)^s$ which has edges set mappings

$$s(u;v) \to \lambda_p(x)(s(u;v)) = 1_x \otimes s(u;v) \subset p(x;x) \otimes s(u;v) \ .$$

The changing of target space with functor $f \in Cat(p \otimes s; r)$ provides the functor $\lambda_p(x) \circ f \in Cat(p; r^s)$ with edges set mappings

$$s(u;v) \to f_{\langle x,u \rangle, \langle x,v \rangle}(1_x \otimes s(u;v)) \subset r(f_0(x,u); f_0(x,v)) \ .$$

Having a natural transform $\alpha : f \to g$ we get the natural transform

$$\lambda_p(x) \times \alpha : \lambda_p(x) \circ f \to \lambda_p(x) \circ g \ .$$

between appointed functors $\lambda_p(x) \circ f, \lambda_p(x) \circ g \in Cat(p; r^s)$ with edges set mappings

$$s(u;v) \to f_{\langle x,u \rangle, \langle x,v \rangle}(1_x \otimes s(u;v)) \geq$$

$$g_{\langle x,u\rangle,\langle x,v\rangle}(1_x \otimes s(u;v)) \subset r(g_0(x,u); g_0(x,v)) \ .$$

Otherwise we get a natural transform also for the edges set mappings. For the sections functor $\lambda_p \in Cat(p; (p \otimes s)^s)$ the edges set mappings will be

$$p(x;y) \to \prod_{u \in s_0} p(x;y) \otimes 1_u \subset \prod_{u \in s_0} p(x;y) \otimes s(u;u) \ .$$

The changing of target space with functor $f \in Cat(p;r)$ provides a new edges set mapping

$$p(x;y) \to \prod_{u \in s_0} f_{\langle x,u\rangle,\langle y,u\rangle}(p(x;y) \otimes 1_u) \subset r(f_0(x,u); f_0(y,u)) \ ,$$

and a natural transform $\alpha : f \to g$ provides the natural transform for the edges set mappings $\lambda_p \times \alpha^s : \lambda_p \circ f^s \to \lambda_p \circ g^s$

$$p(x;y) \to \prod_{u \in s_0} f_{\langle x,u\rangle,\langle y,u\rangle}(p(x;y) \otimes 1_u) \geq$$

$$\prod_{u \in s_0} g_{\langle x,u\rangle,\langle y,u\rangle}(p(x;y) \otimes 1_u) \subset \prod_{u \in s_0} r(g_0(x,u); g_0(y,u)) \ .$$

The composition of arrows in a category $r$ provides the contact arrows, and we get the evaluation functor over exponential space of natural transforms $\mathbf{ev}_r \in Cat(r^s \otimes s; r)$. The composition with this functor

$$g^\flat = (g \otimes 1_s) \circ \mathbf{ev}_r \in Cat(p \otimes s; r) \uparrow g \in Cat(p; r^s)$$

provides the realization of functional

$$Cat(p; r^s) \to Cat(p \otimes s; r)$$

which is monotone from the functionals to the exponential space of natural transforms with order relation defined by existing natural transform between the secondarily parallel mappings to the forms into the target space $r$ with an order relation defined by existing natural transform between (primarily ) parallel functors.

Let we have two secondarily parallel functors $f, g \in Cat(p; r^s)$ with a natural transform $\alpha : f \to g$. defined by collection of arrows

$$\alpha_x : f_0(x) \to g_0(x) \uparrow x \in p_0$$

between functors $f_0(x), g_0(x) \in (r^s)_0$. We have commuting diagrams for edges set mappings for every arrow in source category $\theta \in p(x;y)$

$$\begin{array}{ccc} f_0(x) & \xrightarrow{\alpha_x} & g_0(x) \\ f_{x,y}(\theta) \downarrow & & \downarrow g_{x,y}(\theta) \\ f_0(y) & \xrightarrow[\alpha_x]{} & g_0(y) \end{array}$$

They are arrows in the exponential space of natural transforms $\alpha_x \in r^s$, so they will be the natural transforms between functors, and they are defined as Carte product of arrows in the target space $r$

$$\alpha_{x,u} \in r((f_0(x))_0(u); (g_0(x))_0(u)) \uparrow u \in s_0 \ .$$

These arrows allows the commuting diagram for any arrow in the source category $\theta \in p(x; y)$

$$\begin{array}{ccc} (f_0(x))_0(u) & \xrightarrow{\alpha_{x,u}} & (g_0(x))_0(u) \\ f_{x,y}(\theta)_u \downarrow & & \downarrow g_{x,y}(\theta)_u \\ (f_0(y))_0(u) & \xrightarrow{\alpha_{y,u}} & (g_0(y))_0(u) \end{array}$$

For a couple of arrows $\theta \in p(x, y)$ and $\kappa \in s(u; v)$ the evaluation functor $\mathbf{ev}_r \in Cat(r^s \otimes s; r)$ can be defined by precontact or postcontact arrows

$$\begin{array}{ccc} f_0(x, u) & \xrightarrow{f_0(x)_{u,v}(\kappa)} & f_0(x, v) \\ f_{\langle x,u\rangle,\langle y,u\rangle}(\theta) \downarrow & \searrow & \downarrow f_{\langle x,v\rangle,\langle y,v\rangle}(\theta) \\ f_0(y, u) & \xrightarrow[f_0(y)_{u,v}(\kappa)]{} & f_0(y; v) \end{array}$$

For the precontact arrow we have a diagram

$$\begin{array}{cc} f_0(x, u) & \\ f_{\langle x,u\rangle,\langle y,u\rangle}(\theta) \downarrow & \searrow \\ f_0(y, u) & \xrightarrow[f_0(y)_{u,v}(\kappa)]{} \quad f_0(y; v) \end{array}$$

All mappings are arrows in the target category $r$, therefore we have for them associative composition, so the commuting transforms for two arrows will provide commuting transform for their composition

$$\begin{array}{ccccc} f_0(x, u) & \xrightarrow{f_{\langle x,u\rangle,\langle y,u\rangle}(\theta)} & f_0(y, u) & \xrightarrow{f_0(y)_{u,v}(\kappa)} & f_0(y, v) \\ \alpha_{\langle x,u\rangle} \downarrow & & \downarrow \alpha_{\langle y,u\rangle} & & \downarrow \alpha_{\langle y,v\rangle} \\ g_0(x, u) & \xrightarrow[g_{\langle x,u\rangle,\langle y,u\rangle}(\theta)]{} & g_0(y, u) & \xrightarrow[g_0(y)_{u,v}(\kappa)]{} & g_0(y, v) \end{array}$$

So we get a natural transform between the realizations of two functors $f, g \in Cat(p \otimes s; r)$, and we have shown that the realization appointment is monotone

$$Cat(p; r^s) \to Cat(p \times s; r) \ .$$

Next we shall prefer to use "small" terminology. We shall say that we have a category $Cat$ with points $p \in Cat_0$ and arrows $f : p \to q$. For them we have truly joint pair of monotone functors. The adjoint functor defined with biproduct $p \otimes s$ is monotone for order relation defined with existing of natural transform between the primarily parallel arrows. The coadjoint functor defined with exponential space of natural transforms $r^s$ is monotone for the order relation defined by existing natural transform between secondarily parallel arrows in exponential space $r^s$. For them we have the name appointment isomorphic for the order relation defined by existing natural transform between the primarily parallel arrows $Cat(p \times s; r)$ and the secondarily parallel arrows $Cat(p; r^s)$ in exponential space also defined by natural transform between two functionals

$$Cat(p \times s; r) =\rightarrow Cat(p; r^s) \ ,$$

Now we shall study the *Cat-enriched spaces*. We shall call the *Cat*-enriched graph $X$ also as a *Cat*-graph. The edges sets between vertexes $a, b \in X_0$ will be categories $X(a;b) \in Cat_0$ and functors between such categories will be taken as mappings.

The values category $Cat$ is the first example of $Cat$-enriched space. It provides useful experience with new notions of $Cat$-enriched spaces. The enriched $Cat$-space has another structure as the category $Cat$ itself. Therefore we shall use another notation for it $Cat^e$. The letter $e$ will remind that we are dealing with enriched space. For arbitrary categories $p, r \in Cat_0$ we take the exponential space as an edges set

$$Cat^e(p; r) = r^p \in Cat_0 \ .$$

The set of arrows between two points $Cat(p; r)$ defines the structure of *Set*-enriched space. We have projection

$$Cat(p; r) = Cat^e(p; r)_0 \ .$$

We shall say that these projections are provided by *translation* of values category. For $Cat$-enriched graphs we can also take arbitrary graphs transport $Cat \Longrightarrow V$. For a category $p \in Cat_0$ the projection to its vertexes set $p_0 \in Set_0$ will define very simple translation $Cat \Longrightarrow Set$ of $Cat$-enriched spaces with many good properties, but we are conscious to know what exact properties we are using in each investigation's step.

The $Cat$-graphs transport $f : X \Longrightarrow Y$ will be defined by an appointment of vertexes $f_0 : X_0 \to Y_0$ and edges set mappings

$$f_{a,b} : X(a;b) \to Y(f_0(a); f_0(b)) \ .$$

We get a category $Grph_{Cat}$ compounded of all $Cat$-graphs transports. Arbitrary graphs transport $D : Cat \Longrightarrow V$ defines the superfunctor of categories

$$(D) : Grph_{Cat} \Longrightarrow Grph_V \ .$$

For a $Cat$-enriched graph $X$ it appoint the new $V$-enriched graph $X^{(D)}$ with the same vertexes and their appointments, and new mappings for graphs transports $f^{(D)} : X^{(D)} \Longrightarrow Y^{(D)}$

$$f^{(D)}_{a,b} : X^{(D)}(a,b) \to Y^{(D)}(f_0(a); f_0(b)) \ .$$

We shall call two $Cat$-graphs transports $f, g : X \Longrightarrow Y$ *primarily parallel* if they have the same appointments of vertexes $f_0 = g_0 : X_0 \to Y_0$. For a primarily parallel $Cat$-graphs transports $f, g : X \implies Y$ we get an order relation from values category $Cat$. For mappings between edges sets

$$f_{a,b} : X(a;b) \to f_{a,b}(X(a;b)) \subset Y(f_0(a); f_0(b)) \ ,$$

$$g_{a,b} : X(a;b) \to g_{a,b}(X(a;b)) \subset Y(g_0(a); g_0(b))$$

the inequality $f_{a,b} \geq g_{a,b}$ means that we can find some natural transform $\alpha_{a,b} : f_{a,b} \to g_{a,b}$ between functors of correspondent categories. We shall come down to the category $X(a;b) \in Cat_0$. For every point $u \in X(a;b)_0$ these functor appoints the couple of points $(f_{a,b})_0(u) \in Y(f_0(a); f_0(b))$ and $(g_{a,b})_0(u) \in Y(g_0(a); g_0(b))$ and such natural transform is defined by collection of arrows from category $Y(f_0(a); f_0(b)) = Y(g_0(a); g_0(b))$

$$\alpha_u : (f_{a,b})_0(u) \to (g_{a,b})_0(u) \uparrow u \in X(a;b)_0$$

We shall define a joint pair of superfunctors within the category $Grph_{Cat}$. First we define biproduct of $Cat$-graphs $X \otimes Y$. We take the Carte tensor product of vertexes sets

$$(X \otimes Y)_0 = X_0 \times Y_0$$

and biproduct of edges sets

$$X \otimes Y(\langle x, y \rangle; \langle x', y' \rangle) = X(x; x') \otimes Y(y; y') \ .$$

The biproduct of $Cat$-graphs transports $f : X \implies X'$ and $g : Y \implies Y'$ will be defined by biproduct of edges set mappings

$$f \otimes g_{\langle x,y \rangle, \langle x', y' \rangle} : X(x; x') \otimes Y(y; y') \to f_{x,x'}(X(x;x')) \otimes g_{y,y'}(Y(y;y')) \subset$$

$$X'(f_0(x); f_0(x')) \otimes Y'(g_0(y); g_0(y')) \ .$$

Such bisuperfunctor is monotone, as we have monotone biproduct in the values category $Cat$. We need to check this property only for the edges set mappings. We shall not be lazy to come down inside the categories $X(x; x')$ and $Y(y; y')$.

Between two parallel functors $f_{x,x'}, f'_{x,x'} : X(x; x') \to X'(f_0(x); f_0(x'))$ a natural transform $\alpha_{x,x'} : f_{x,x'} \to f'_{x,x'}$ is defined by some collection of arrows from the category $X'(f_0(x); f_0(x'))$

$$\alpha_u : (f_{x,x'})_0(u) \to (f'_{x,x'})_0(u) \uparrow u \in X(x; x')_0 \ .$$

Also a natural transform $\beta_{y,y'} : g_{y,y'} \to g'_{y,y'}$ between two parallel functors $g_{y,y'}, g'_{y,y'} : Y(y; y') \to Y'(g_0(y); g_0(y'))$ is defined by another collection of arrows from category $Y'(g_0(y); g_0(y'))$

$$\beta_v : (g_{y,y'})_0(v) \to (g_{y,y'})_0(v) \uparrow v \in Y(y; y')_0 \ .$$

The Carte product of such natural transforms will be defined by collection of arrows couples from the category $X'(f_0(x); f_0(x')) \times Y'(g_0(y); g_0(y'))$

$$\langle \alpha_u, \beta_v \rangle : \langle (f_{x,x'})_0(u), g_{y,y'})_0(v) \rangle \to \langle (f'_{x,x'})_0(u), (g'_{y,y'})_0(v) \rangle \ .$$

And it will be a natural transform between biproducts of taken functors

$$\alpha_{x,x'} \otimes \beta_{y,y'} : f_{x,x'} \otimes g_{y,y'} \to f'_{x,x'} \otimes g'_{y,y'} .$$

We shall work with a for chosen assistant $Cat$-graph $T$ defined superfunctor. For a $Cat$-graph $X$ it appoints the biproduct with the assistant $Cat$-graph $X \otimes T$, and for a $Cat$-graphs transport $f : X \Longrightarrow X'$ it appoints the biproduct with identity $Cat$-graphs transport of the assistant $Cat$-graph $T$

$$f \otimes \mathtt{Id}_T : X \times T \Longrightarrow X' \otimes T .$$

We shall take it as an adjoint superfunctor within the category of $Cat$-graphs $Grph_{Cat}$.

The coadjoint superfunctor will be defined with exponential space $Y^T$ of all transforms $\alpha : \phi \to \psi$ between $Cat$-graphs transports $\phi, \psi : T \Longrightarrow Y$. The vertexes will be all $Cat$-graphs transports $\phi : T \Longrightarrow Y$ and edges sets are provided by the Carte product of edges sets in the target space $Y$

$$Y^T(\phi; \psi) = \prod_{s \in T_0} Y(\phi_0(s); \psi_0(s)) .$$

Obviously an assistant $Cat$-graph $T$ would have a vertexes set $T_0$ small enough that such Carte product could exist.

For a $Cat$-graphs transport $f : Y \Longrightarrow Y'$ such superfunctor appoints the changing of target space with taken $Cat$-graphs transport

$$f^T : Y^T \Longrightarrow (Y')^T .$$

For a $Cat$-graphs transport $\phi : T \Longrightarrow Y$ it appoints the composition of $Cat$-graphs transports $\phi \circ f : T \Longrightarrow Y'$ and for an edges set $Y^T(\phi; \psi)$ the mapping is defined with the Carte product of edges set mappings in the target space

$$f^T_{\phi,\psi} : Y^T(\phi; \psi) = \prod_{s \in T_0} Y(\phi_0(s); \psi_0(s)) \to$$

$$\prod_{s \in T_0} f_{\phi_0(s);\psi_0(s)}(Y(\phi_0(s); \psi_0(s))) \subset \prod_{s \in T_0} Y'(f_0(\phi_0(s)); f_0(\psi_0(s))) .$$

For primarily parallel $Cat$-graphs transports $f, g : Y \Longrightarrow Y'$ we get only secondarily parallel $Cat$-graphs transports $f^T, g^T : Y^T \Longrightarrow (Y')^T$, i. e. we don't have coincident appointments of vertexes $f_0^T(\phi) \neq g_0^T(\phi)$, but we get coincidence only for the projection to vertexes appointment

$$f_0(\phi_0(s)) = g_0(\phi_0(s)) \uparrow \phi \in (Y^T)_0, s \in T_0 .$$

For a natural transform between the edges set mappings $f_{x,y} \geq g_{x.y}$ we get the natural transform between vertexes appointments

$$T(s;t) \to f_{\phi_0(s);\psi_0(t)}(\phi_{s,t}(T(s;t)) \geq g_{\phi_0(s);\psi_0(t)}(\phi_{s,t}(T(s;t)) ,$$

and another natural transform between the edges set mappings

$$Y(x;y) \to \prod_{s \in T_0} f_{\phi_0(s);\psi_0(s)}(Y(\phi_0(s);\psi_0(s))) \geq$$

$$\prod_{s \in T_0} g_{\phi_0(s);\psi_0(s)}(Y(\phi_0(s);\psi_0(s))) \subset \prod_{s \in T_0} Y'(g_0(\phi_0(s));g_0(\psi_0(s))) \ .$$

We needed to check this property only for the target space changing in the values category $Cat$. We shall come down inside the categories $X(x;y)$ and $T(s;t)$ again. Let we have a natural transform $\alpha_{x,y} : f_{x,y} \to g_{x,y}$ defined by collection of arrows from the category $Y(f_0(x); f_0(y))$

$$\alpha_w : (f_{x,x'})_0(w) \to (g_{x,x'})_0(w) \uparrow w \in X(x;y)_0 \ .$$

For some $Cat$-graphs transport $\phi : T \Longrightarrow X$ the changing of target space with another $Cat$-graphs transport $f : X \Longrightarrow Y$ will appoint the composition of $Cat$-graphs transports $\phi \circ f : X \Longrightarrow Y$ which will have the edges set mappings

$$T(s;t) \to f_{\phi_0(s),\phi_0(t)}(\phi_{s,t}(T(s;t))) \subset Y(f_0(\phi_0(s)); f_0(\phi_0(t))) \ .$$

The natural transforms $\alpha_{x,y} : f_{x,y} \to g_{x,y}$ now provides a new natural transform $\alpha^{\phi_{s,t}(T(s;t))} : \phi_{s,t} \circ f_{\phi_0(s),\phi_0(t)} \to \phi_{s,t} \circ g_{\phi_0(s);\phi_0(t)}$ defined with collection of arrows from the category $Y(f_0(\phi_0(s)); f_0(\phi_0(t))) = Y(g_0(\phi_0(s)); g_0(\phi_0(t)))$

$$\alpha_{(\phi_{s,t})_0(u)} : (f_{\phi_0(s),\phi_0(t)})_0((\phi_{s,t})_0(u)) \to (g_{\phi_0(s),\phi_0(t)})_0((\phi_{s,t})_0(u))$$

for $u \ in T(s;t)_0$. The notation $\alpha^{\phi_{s,t}(T(s;t))}$ denotes the collection of arrows for the points in the category $\phi_{s,t}(T(s;t))$.

For a edges set $X^T(\phi;\psi)$ between $Cat$-graphs transports $\phi,\psi : T \Longrightarrow X$ the changing of target space with $Cat$-graphs transport $f : X \Longrightarrow Y$ provides the edges set mapping as the Carte product of mappings

$$\prod_{s \in T_0} X(\phi_0(s);\psi_0(s)) \to \prod_{s \in T_0} f_{\phi_0(s),\psi_0(s)}(X(\phi_0(s);\psi_0(s))) \subset$$

$$\prod_{s \in T_0} Y(f_0(\phi_0(s)); f_0(\psi_0(s))) \ .$$

The natural transforms $\alpha_{x,y} : f_{x,y} \to g_{x,y}$ now provides a new natural transform

$$\alpha^{X^T(\phi;\psi)} : \prod_{s \in T_0} f_{\phi_0(s),\psi_0(s)} \to \prod_{s \in T_0} g_{\phi_0(s),\psi_0(s)}$$

with Carte product of arrows

$$(f_{\phi_0(s),\psi_0(s)})_0(w)) \to (g_{\phi_0(s),\psi_0(s)})_0(w) \uparrow w \in X(\phi_0(s);\psi_0(s))_0, s \in T_0 \ .$$

We should want to take arbitrary arrow $e_x \in X(x;x)$ as a unit arrow. Let $e_x : u \to v$ is an arrow between two points $u, v \in X(x;x)_0$. The choosing of unit arrow would allow us to define sections transport $\lambda_X : X \Longrightarrow (X \otimes T)^T$. The vertexes appointment for a vertex $x \in X_0$ appoints the section $\lambda_X(x) \in X^T$. It will be a $Cat$-graphs transport for a vertex $s \in T$ appointing a vertex in the biproduct $\lambda_X(x,s) = \langle x, s \rangle \in X_0 \otimes T_0$ and having edges set mappings

$$\lambda_X(x)_{s,t} : T(s;t) \to e_x(*) \otimes T(s;t) \subset X(x;x) \otimes T(s;t) .$$

The edges set mappings for sections transport will be

$$(\lambda_X)_{x,y} : X(x;y) \to (X \otimes T)^T(\lambda_X(x); \lambda_X(y)) = \prod_{s \in T_0} X(x;y) \otimes e_s(*) \subset$$

$$\prod_{s \in T_0} X(x;y) \otimes T(s;s) .$$

These mappings for $Cat$-enriched graphs $X$ and $T$ must be functors between correspondent categories. However the edges set mappings between section $\Lambda_X(x)$

$$\Lambda_X(x)_{s,t} : T(s;t) \to e_x(*) \otimes T(s;t)$$

doesn't maintain composition of arrows from category $T(s;t)$ if the unit arrow doesn't allow the equality of idempotent

$$e_x \circ e_x = e_x .$$

Otherwise for an arrow $e_x : v \to w$ with the different points $v \in X(x;x)_0$ and $v \in X(x;x)_0$ defined mapping $T(s;t) \to e_x(*) \otimes T(s;t)$ can't be functor, as for some point $u \in T(s;t)_0$ we get an arrow $\langle e_x, u \rangle :\in X(x;x) \otimes T(s;t)$ between different couples $\langle v, u \rangle \neq \langle w, u \rangle$. We need extend the stock of mappings allowing new functors $f : p \to q^{\mathbf{2}}$. We have denoted the category $\mathbf{2} = \{\cdot \to \cdot\}$ with two points and three arrows. The composition of such arrows will extends the stock of arrows. We must include also the arrows $f : p \to q^{\mathbf{2}^n}$ with arbitrary finite number $n \geq 0$. So we need to extend the values category $Cat$ with graphs transports $p \to q^{\mathbf{2}^n}$ between categories $p, q \in Cat_0$. Such categories are original graphs, and graphs transports $p \to q^{\mathbf{2}^n}$ can be taken maintaining unit arrows. Such new category isn't closed one, but has name and realization mappings with simply properties.

At this moment we shall work with our taken category $Cat$. The collection of sections transports defines the name appointment from the set of $Cat$-graphs transports over biproduct $X \otimes T$ with the assistant $Cat$-graph $T$ to the set of $Cat$-graphs transports into the exponential space $Y^T$

$$Grph_{Cat}(X \otimes T; Y) \to Grph_{Cat}(X; Y^T) .$$

For a $Cat$-graphs transport $f : X \otimes T \Longrightarrow Y$ it provides the composition of $Cat$-graphs transports

$$f^\# = \lambda_X \circ f^T : X \Longrightarrow Y^T .$$

It is monotone for order relation defined by existing natural transform between primarily parallel $Cat$-graphs transports to the target $Cat$-graph $Y$ and existing natural transforms between secondarily parallel $Cat$-graphs transports to the exponential space of all transforms $Y^T$. It can be checked for edges set mappings which will be functors between correspondent categories.

We again will come down inside the categories $p \in Cat_0$. The natural transform $\alpha_{\langle x,s \rangle, \langle y,t \rangle} : f_{\langle x,s \rangle, \langle y,t \rangle} \to f'_{\langle x,s \rangle, \langle y,t \rangle}$ between two $Cat$-graphs transports $f, f' \in Cat(X \otimes T; Y)$ is defined by collection of arrows from category $Y(f_0(x,s); f_0(y,t))$

$$\alpha_{\langle w,u \rangle} : (f_{\langle x,s \rangle, \langle y,t \rangle})_0(w,u) \to (f'_{\langle x,s \rangle, \langle y,t \rangle})_0(w,u) \uparrow w \in X(x;y)_0, u \in T(s;t)_0 \ .$$

For any vertex $x \in X_0$ appointed section $\Lambda_X(x) : T \Longrightarrow X \otimes T$ is defined with vertexes appointment $\langle e_x, s \rangle \in X(x;x)_0 \times T(s;t)_0$ and edges set mappings $\lambda_X(x)_{s,t} : T(s;t) \to (X(x;x) \otimes T(s;t))^2$. Further the changing of target space by a $Cat$-graphs transport $f : X \otimes T \Longrightarrow Y$ appoints the $Cat$-graphs transport $T \Longrightarrow Y$ which for the vertex $s \in T_0$ appoints the vertex $f_0(x,s) \in Y_0$ and have edges set mappings

$$T(s,t) \to f^{\mathbf{2}}_{\langle x,s \rangle, \langle x,t \rangle}(e_x(*) \otimes T(s;t)) \subset Y(f_0(x,s); f_0(x,t))^{\mathbf{2}} \ .$$

We get a natural transform $\alpha^{\lambda_X(x)_{s,t}(T(s;t))} : f_{\langle x,s \rangle, \langle x,t \rangle} \to f'_{\langle x,s \rangle, \langle x,t \rangle}$ defined with collection of arrows from category $Y(f_0(x,s); f_0(y,t))$ for the source and target points

$$\alpha^1_{(\lambda_X(x)_{s,t})_0(u)} = \alpha_{\langle e^1_x, u \rangle} : (f_{\langle x,s \rangle, \langle x,t \rangle})_0(e^1_x, u) \to (f'_{\langle x,s \rangle, \langle x,t \rangle})_0(e^1_x, u) \ ,$$

$$\alpha^2_{(\lambda_X(x)_{s,t})_0(u)} = \alpha_{\langle e^2_x, u \rangle} : (f_{\langle x,s \rangle, \langle x,t \rangle})_0(e^2_x, u) \to (f'_{\langle x,s \rangle, \langle x,t \rangle})_0(e^2_x, v) \ .$$

For the edges set mappings $X(x,y) \to \prod_{s \in T_0} f_{\langle x,s \rangle, \langle y,s \rangle}$ the natural transform is defined with natural transforms between taken $Cat$-graphs transports

$$\alpha_{\langle x,s \rangle, \langle y,t \rangle} : f_{\langle x,s \rangle, \langle y,t \rangle} \to f'_{\langle x,s \rangle, \langle y,t \rangle} \ .$$

They provide a new natural transform with arrows

$$\alpha^{(X \otimes T)^T} : \prod_{s \in T_0} (f_{\langle x,s \rangle, \langle y,s \rangle})_0(w) \to \prod_{s \in T_0} (f_{\langle x,s \rangle, \langle y,s \rangle})_0(w)$$

equal to Carte product of arrows couples

$$\alpha^1_{\langle x,s \rangle, \langle y,s \rangle}(w, e^1_s) : (f_{\langle x,s \rangle, \langle y,s \rangle})_0(w, e^1_s) \to (f'_{\langle x,s \rangle, \langle y,s \rangle})_0(w, e^1_s) \uparrow u \in T(s;s)_0 \ ,$$

$$\alpha^2_{\langle x,s \rangle, \langle y,s \rangle}(w, e^2_s) : f_{\langle x,s \rangle, \langle y,s \rangle_0})_0(w, e^2_s) \to (f'_{\langle x,s \rangle, \langle y,s \rangle})_0(w, e^2_s) \uparrow u \in T(s;s)_0 \ .$$

Now we shall construct an evaluation functor. In a $Cat$-graph $X$ the contact arrows mappings Will be defined as bifunctors

$$h_{a,b,c} : X(a;b) \otimes X(b;c) \to X(a;c) \ .$$

They defines some $Cat$-graphs transport over *pullbacked product*

$$h : X \otimes_{X_0} X \Longrightarrow X .$$

The pullbacked product of enriched graphs is defined by the vertexes set $X_0 \times X_0$ and the edges sets

$$X \otimes_{X_0} X(p;r) = \bigcup_{q \in X_0} X(p;q) \otimes X(q;r) .$$

The taken bifunctors define edges set mappings

$$\bigcup_{q \in X_0} h_{p,q,r} : \bigcup_{q \in X_0} X(p;q) \otimes X(q;r) \to X(p;r) .$$

We demand that changing of target space with such bifunctors would be monotone

$$Cat(r, X(a,b) \otimes X(b;c)) \to Cat(r; X(a;b) \circ X(b;c)) \subset Cat(r; X(a;c)) ,$$

This means that for the couple of natural transforms $\alpha : f \to f'$ between any two functors $f, f' \in Cat(r; X(a;b))$ and $\beta : g \to g'$ between two functors $g, g' \in Cat(r; X(b;c))$ such changing of target space would provide a new natural transform. In our case the $Cat$-graphs functors always keeps natural transforms. The image of natural transform $\alpha$ by functor $F$ can be denoted $\alpha \times F$.

The contact arrows mappings allows us to define preevaluation

$$\mathbf{ev}_Y^+ : Y^T \otimes T \Longrightarrow Y .$$

Preevaluation for the edges sets $Y^T(\phi; \psi)$ and $T(s;t)$ appoints precontact edges set

$$Y^T(\phi; \psi) \otimes T(s;t) \to Y(\phi_0(s); \psi_0(t)) \circ \psi_{s,t}(T(s;t)) \subset Y(\phi_0(s); \psi_0(t)) .$$

The prerealization appointment

$$Grph_{Cat}(X; Y^T) \to Grph_{Cat}(X \otimes T; Y)$$

is provided by composition of $Cat$-graphs transports. For a $Cat$-graphs transport $g : X \Longrightarrow Y^T$ it appoints the composition

$$g^\flat = (g \otimes \mathtt{Id}_T) \circ \mathbf{ev}_Y^+ : X \otimes T \Longrightarrow Y .$$

We can check that for a natural transforms $\alpha_x : g_0(x) \to g'_0(x)$ between $Cat$-graphs transports $g_0(x), g'_0(x) : T \Longrightarrow Y$ for vertexes $x \in X_0$ and natural transforms $\alpha_{x,y} : g_{x,y} \to g'_{x,y}$ between the edges set mappings

$$X(x;y) \to g_{x,y}(X(x;y)) \geq g'_{x,y}(X(x;y)) \subset Y^T(g'_0(x); g'_0(y))$$

of secondarily parallel $Cat$-graphs transports $g, g' : X \Longrightarrow Y^T$ we get a new natural transform between edges set mappings

$$X(x;y) \otimes T(s;t) \to Y(g_0(x,s); g_0(y,s)) \geq Y(g'_0(x,s); g'_0(y,s)) \ .$$

Let we have a natural transforms $\alpha_x : g_0(x) \to g'_0(x)$ between edges set mappings $g_0(x)_{s,t}, g'_0(x)_{s,t} : T(s;t) \Longrightarrow Y(g_0(x,s); g_0(x,t))$ defined by collection of arrows

$$(\alpha_x)_u : (g_0(x)_{s,t})_0(u) \to (g'_0(x)_{s,t})_0(u) \uparrow u \in T(s;t)_0$$

and natural transforms $\alpha_{x,y} : g_{x,y} \to g'_{x,y}$ between two edges set mappings $X(x;y) \to Y^T(g_0(x); g_0(y)))$. These mappings are defined as a Carte product

$$g_{x,y} = \prod_{s \in T_0} g_{x,y}(s)$$

of functors $g_{x;y}(s) : X(x;y) \to Y(g_0(x;s); g_0(y;s))$, therefore we get the Carte product

$$\alpha_{x,y} = \prod_{s \in T_0} \alpha_{x,y}(s)$$

of natural transforms

$$(\alpha_{x,y}(s))_w : (g_{x,y}(s))_0(w) \to (g'_{x,y}(s))_0(w)) \uparrow w \in X(x;y)$$

defined with arrows in the category $Y(g_0(x,s); g_0(y,s))$.

The realization provided by the preevaluation transport is defined as precontact arrows mappings

$$X(x;y) \otimes T(s;t) \to Y((g_0(x))_0(s); (g_0(y))_0(s)) \circ Y((g_0(y))_0(s); (g_0(y))_0(t)) \subset$$

$$Y((g_0(x))_0(s); (g_0(y))_0(t)) \ .$$

Having the natural transforms between edges set mappings

$$f_{\langle x,s \rangle, \langle y,t \rangle}, f'_{\langle x,s \rangle, \langle y,t \rangle} : X(x;y) \otimes T(s;t) \to Y(g_0(x,s); g_0(y,t)))$$

we get the needed natural transform for the precontact arrows mappings by monotone property of changing target space by contact arrows mappings.

In the case of mappings $g_0(x)_{s,t} : T(s;t) \to Y(g_0(x,s); g_0(x,t))^{\mathbf{2}}$ and

$$g_{x,y} : X(x;y) \to Y(g_0(x,s); g_0(y,s))^{\mathbf{2}}$$

we would get a realization $f : X \otimes T \Longrightarrow Y$ with edges set mappings

$$X(x;y) \otimes T(s;t) \to Y(g_0(x,s); g_0(y,t))^{\mathbf{2}^2} \ .$$

Such realization remains monotone for order relation defined with existent natural transforms between functors.

We begin studies of original $Cat$-graphs. The original $Cat$-enriched graph $X$ needs equalities for chosen unit arrows and contact arrows mappings

$$X(a;b) \to e_a(*) \circ X(a;b) = X(a;b) = X(a;b) \circ e_b(*) \subset X(a;b) \ .$$

At first we shall check such equalities for translation of values category $Cat \implies Set$. The unit arrow $e_x \in X(x;x)$ becomes some point among all points $O_x \in X(x;x)_0$ of category $X(x;x)$. For original $Set$-graph $X_0$ we need the equalities

$$X(a;b)_0 \to O_a(*) \circ X(a;b)_0 = X(a;b)_0 = X(a;b) \circ O_b(*) \ .$$

It means that for arbitrary point $u \in X(a;b)_0$ we have the equalities

$$O_a \circ u = u = u \circ O_b \ .$$

Otherwise for points $u \in X(a;a)$ and the unit point $O_a \in X(a;a)_0$ we get the equalities

$$O_a \circ u = u = u \circ O_a \uparrow u \in X(a;a)_0 \ .$$

We can see that a unit point is unique. If we have another unit point $u \in X(a;a)_0$, then we get equalities

$$O_x = O_x \circ u = u \ .$$

So for the original $Set$-graph $X_0$ all sets $X(x;x)_0$ has unique unit point

$$O_x \in X(x;x)_0 \ .$$

In arbitrary set $X(x;x)_0$ we can add a new unit $O_x \in X(x;x)_0$ and in the set $X(x;x)_0^* = X(x;x)_0 \cup \{O_x\}$ we can extend contact arrows mapping by taking

$$O_x \circ u = u \circ O_x = u \uparrow u \in X(x;x)_0 \ , \quad O_x \circ O_x = O_x \ .$$

The old unit $e \in X(x;x)_0$ looses the property of unit in $X(x;x)_0^*$ as we have

$$e \circ O_x = e \ .$$

Such construction allows us arbitrary $Set$-graph $X$ to change with original $Set$-graph $X^*$. More exactly the category of $Set$-graphs $Grph_{Set}$ is isomorphic to the category of original $Set$-graphs $OGrph_{Set}$ with arrows of all original $Set$-graphs transports. For a $Set$-graphs transport $f_0 : X_0 \implies Y_0$ we appoint its extension $f^* : X^* \implies Y^*$ with

$$(f^*_{x,x})_0(O_x) = O_{f_0(x)} \in Y(f_0(x); f_0(x))_0 \ .$$

Obviously we get a superfunctor $Grph_{Set} \implies OGrph_{Set}$, and inverse superfunctor will be defined by restriction of edges set mappings

$$(f^*_{x,x})_0 : X(x;x)_0^* \to Y^*(f_0(x); f_0(x))_0$$

The same construction can be repeated for a $Cat$-enriched graph $X$. For a category $X(x;x)$ we appoint the set of category's points $X(x;x)_0$. Let we have find one origin point $O_x \in X(x;x)_0$ which is a neutral for the projection of contact arrows mapping to the mappings of $Set$-enriched graph

$$X(x;x)_0 \times X(x;x)_0 \to X(x;x)_0 .$$

In the set of loops with source and target point $O_x$ we shall look for a unit arrow $1_x : O_x \to O_x$ with the equality for the contact arrows mapping

$$1_x * 1_x = 1_x .$$

Otherwise we can add such unit arrow formally. In such way we extend taken category $X(x;x)^*$. Such extension is possible also even having true unit arrow. We again can unanimously extended every functor to the original functor over extended category $X(x;x)^*$. For arbitrary $Cat$-graphs transport $f : X \Longrightarrow Y$ we get original $Cat$-graphs transport

$$f^* : X^* \Longrightarrow Y^* .$$

The restriction of original $Cat$-graphs transport $g : X^* \Longrightarrow Y^*$ provides non original $Cat$-graphs transport $\hat{g} : X \Longrightarrow Y$. We again get isomorphism of categories

$$Grph_{Cat} = \Rightarrow OGrph_{Cat} .$$

However we don't need to apply such construction for all $Cat$-graphs. When we are studying the name appointment

$$Grph_{Cat(}(X \otimes T; Y) \to Grph_{Cat}(X; Y^T) ,$$

we shall apply such construction only for three taken $Cat$-graphs $X^*$, $Y^*$, $T^*$ and we shall probe to construct needed correct extension for biproduct $X \otimes T$ and exponential space $Y^T$.

We shall add a unit loop $1_x : O_x \to O_x$ with additional arrows to other points of category $X(x;x)$. We shall say that we have a category $X(x;x)$ with a *leading unit.* Also the $Cat$-graph $X$ will be called *leading* if it have a leading unit. Let we take leading extensions of source and assistant $Cat$-graphs $X$ and $T$. For the target $Cat$-graph $Y$ we ask existence of correspondent limits, needed to define extension of taken $Cat$-graphs transport $f : X \otimes T \Longrightarrow Y$. In ordered spaces such limits coincide with lower boundaries for the points. In general case we should work with arrow spaces defined as an exponential space $p^{\mathbf{2}}$.

For mapping $f_{\langle x,s \rangle, \langle y,s \rangle} X(x;y) \to Y(f_0(x,s); f_0(y,s))$ we take

$$f^*_{\langle x,s \rangle, \langle y,s \rangle}(w, 1_s) = \lim_{u \to 1_s} f_{\langle x,s \rangle, \langle y,s \rangle}(w, u) \in Y(f_0(x,s); f_0(y,s)) ,$$

$$f^*_{\langle x,s \rangle, \langle x.t \rangle}(1_x, u) = \lim_{w \to 1_x} f_{\langle x,s \rangle, \langle x.t \rangle}(w, u) \in Y(f_0(x,s); f_0(x,t)) ,$$

$$f^*_{\langle x,s \rangle, \langle x,s \rangle}(1_x; 1_s) = 1_{f_0(x,s)} \in Y(f_0(x,s); f_0(x,s)) .$$

These limits can be not unique, then we should get different possible extensions. Any such extension will be useful to explain situation with taken $Cat$-graphs $X$, $Y$, $T$. Also the *partial transport* $g : X \Longrightarrow Y^T$ will be changed with its extension.

We want to remark that we get opportunity to understand a general theory for the name appointment

$$Grph_{Cat}(X \otimes T; Y) \to Grph_{Cat}(X; Y^T) \ .$$

From original graphs $X$, $Y$, $T$ we can remove unit arrows, and to produce all patterns of general case.

Let we have an original $Cat$-graphs transport $f : X \otimes T \Longrightarrow Y$ with upper predecomposition

$$f_{\langle a,s\rangle,\langle b,s\rangle}(X(a;b) \otimes 1_s(*)) \circ f_{\langle b,s\rangle,\langle b,t\rangle}(1_b(*) \otimes T(s;t)) \geq$$

$$f_{\langle a,s\rangle,\langle b,t\rangle}(X(a;b) \otimes T(s;t)) \ .$$

Then appointed prerealization of name $f^+ \geq f$ is the smallest preregular upper bound, i. e. a original $Cat$-graph transport coinciding with the name's prerealization.

For this $Cat$-graphs transport we have the predecomposition equality

$$f^+_{\langle a,s\rangle,\langle b,s\rangle}(X(a;b) \otimes 1_s(*)) \circ f^+_{\langle b,s\rangle,\langle b,t\rangle}(1_b(*) \otimes T(s;t)) =$$

$$f^+_{\langle a,s\rangle,\langle b,t\rangle}(X(a;b) \otimes T(s;t)) \ .$$

Otherwise each original $Cat$-graphs transport $f : X \otimes T \Longrightarrow Y$ with predecomposition equality is preregular.

The name appointment is bijective over the set of preregular $Cat$-graphs transports. This means that for arbitrary $Cat$-graphs transport $g : X \Longrightarrow Y^T$ we have the realization $f : X \otimes T \Longrightarrow Y$ being some preregular $Cat$-graphs transport.

We want to remark that for original $Cat$-graphs transports $f : X \otimes T \Longrightarrow Y$ we use the properties of contact arrows mappings only in the target space $Y$. In the source space $X$ we need only that chosen units will be neutral for taken contact arrows mappings.

Finally we have seen that the name $g : X \Longrightarrow Y^T$ of arbitrary original $Cat$-graphs transport $f : X \times T \Longrightarrow Y$ with upper predecomposition and lower postdecomposition

$$f_{\langle a,s\rangle,\langle b,s\rangle}(X(a;b) \otimes 1_s((*))) \circ f_{\langle b,s\rangle,\langle b,t\rangle}(1_b(*) \otimes T(s;t)) \geq$$

$$f_{\langle a,s\rangle,\langle b,t\rangle}(X(a;b) \otimes T(s;t)) \geq$$

$$f_{\langle a,s\rangle,\langle a,t\rangle}(1_a(*) \otimes T(s;t)) \circ f_{\langle a,t\rangle,\langle b,t\rangle}(X(a;b) \otimes 1_t(*))$$

is corestricted in exponential space of natural transforms

$$g : X \Longrightarrow Y^{(T)} \ .$$

It rests to consider the $Cat$-enriched categories. The $Cat$-enriched category $X$ will be arbitrary original $Cat$-graph with contact arrows mappings providing associative composition, i. e. we can apply for the *associativity mappings*

$$S_{p,q,r} p \otimes (q \otimes r) \to (p \otimes q) \otimes r$$

the *associativity diagram*

$$S_{X(x;y),X(y;z),X(z;u)} \circ (X_{x,y,z} \otimes 1_{X(z,u)}) \circ X_{x,z,u} =$$

$$(1_{X(x,y)} \otimes X_{y,z,u}) \circ X_{x,y,u}$$

$$\begin{array}{ccc}
X(x;y) \otimes (X(y;z) \otimes X(z;u)) & \xrightarrow{S} & (X(x;y) \otimes X(y;z)) \otimes X(z,;u) \\
{\scriptstyle 1_{X(x;y)} \otimes X_{y,z,u}} \downarrow & & \downarrow {\scriptstyle X_{x,y,z} \otimes 1_{X(z,u)}} \\
X(x;y) \otimes X(y;u) & & X(x;z) \otimes X(z;u) \\
{\scriptstyle X_{x,z,u}} \searrow & & \swarrow {\scriptstyle X_{x,z,u}} \\
& X(x;u) &
\end{array}$$

We shall work with various categories of $Cat$-enriched categories. Their transports $f: X \Longrightarrow Y$ will be called functors. The crude functors are transports of original $Cat$-graphs. They maintain the unit arrows

$$f_{x,x}(1_x(*)) = 1_{f_0(x)}(*) \subset Y(f_0(x); f_0(x)) \ .$$

This category I shall denote, in contrast of other users, most simply $Cat_{Cat}$.

Then we take continuous or cocontinuous functors having correspondent diagrams for contact arrows mappings

$$\begin{array}{ccc}
X(a;b) \otimes X(b;c) & \xrightarrow{f_{a,b} \otimes f_{b,c}} & Y(f_0(a); f_0(b)) \otimes Y(f_0(b); f_0(c)) \\
\downarrow & & \downarrow \\
X(a;c) & \xrightarrow[f_{a,c}]{} & Y(f_0(a); f_0(c)) \\
& & \nearrow
\end{array}$$

$$\begin{array}{ccc}
X(a;b) \otimes X(b;c) & \xrightarrow{f_{a,b} \otimes f_{b,c}} & Y(f_0(a); f_0(b)) \otimes Y(f_0(b); f_0(c)) \\
\downarrow & & \downarrow \\
X(a;c) & \xrightarrow[f_{a,c}]{} & Y(f_0(a); f_0(c)) \\
& & \swarrow
\end{array}$$

Their categories will be denoted $Cat^+_{Cat}$ and $Cat^-_{Cat}$. Their intersection also will have own denotation

$$Cat^{\pm}_{Cat} = Cat^+_{Cat} \bigcap Cat^-_{Cat} \ .$$

I shall call the usual functors $f : X \Longrightarrow Y$ as natural functors. They maintain the contact arrows mappings, i. e. we have commuting diagrams

$$\begin{array}{ccc} X(a;b) \otimes X(b;c) & \xrightarrow{f_{a,b} \otimes f_{b,c}} & Y(f_0(a); f_0(b)) \otimes Y(f_0(b); f_0(c)) \\ \downarrow & & \downarrow \\ X(a;c) & \xrightarrow[f_{a,c}]{} & Y(f_0(a); f_0(c)) \end{array}$$

Their category will be denoted $Cat^=_{Cat}$. For it the name appointment provides the true conjunction between adjoint superfunctor of biproduct $X \otimes T$ and coadjoint functor of exponential space $Y^{(T)}$ of natural transforms. Therefore the biproduct is defined by exponential space unanimously up unique isomorphism.

We can see that the values category $Cat$ becomes a $Cat$-graph $Cat^e$ with edges sets $p^s$ between vertexes $s, p \in Cat_0$. The true edges sets $Cat(s; p)$ will be the edges sets for $Set$-graph $Cat$. It can be identified with the set of all functors $Cat^{\mathbf{2}}$, where $\mathbf{2} = \{\cdot \to \cdot\}$

The defined functors $q^p \times r^q \to r^p$ defines edges sets mappings for $Cat$-graphs transport over pullbacked product of $Cat$-graphs

$$Cat^e \times_{Cat_0} Cat^e \Longrightarrow Cat^e \ .$$

The pullbacked product for the $Set$-graphs $Cat$ can be defined as the set of composable functors couples

$$Cat^{\mathbf{2}} \times_{Cat_0} Cat^{\mathbf{2}} = \bigcup_{q \in Cat_0} Cat(p; q) \times Cat(q; r) \ .$$

Its projection to $Set$-graph $Cat^{\mathbf{2}}$

$$Cat^{\mathbf{2}} \times_{Cat_0} Cat^{\mathbf{2}} \Longrightarrow Cat^{\mathbf{2}}$$

is defined with the composition of functors

$$f \circ g \in Cat(p; r) \uparrow f \in Cat(p; q), g \in Cat(q; r) \ .$$

It is interesting to understand is the horizontal composition of natural transform a unique possibility for the $Cat$-graph $Cat^e$ to become $Cat$-enriched category.

In the category $Cat^e$ we have evaluation mappings

$$\mathbf{ev}_q : p \otimes q^p \to q \ .$$

We can compose two evaluation mappings

$$p \otimes (q^p \otimes r^q) = (p \otimes q^p) \otimes r^q \to q \otimes r^q \to r \ ,$$

and get the arrow with the name

$$q^p \otimes r^q \to r^p \ .$$

I am sure that we could find also another useful contact arrows mappings in $Cat$-graph $Cat^e$. But then we need to check that we have got mappings coincided with horizontal composition of natural transforms, This can be done by careful calculation. We present the *graphical instruction* for such calculations.

In **G. Valiukevičius 1992** p. 102 I have proposed to express the transform $\alpha : \phi \to \psi$ between parallel functors $\phi, \psi : p \Longrightarrow q$ with *deformed squares*

$$\begin{array}{ccc} & \phi & \searrow \\ p & \longrightarrow & q \\ \uparrow & & \uparrow \\ p & \longrightarrow & q \\ & \psi & \end{array}$$

We take identity functors $\text{Id}_p : p \Longrightarrow p$ and $\text{Id}_q : q \Longrightarrow q$ as *uparrows* in such diagrams. In this article we have called such squares as continuous squares.

The *prehorizontal composition* $\alpha * \beta : \phi \circ \xi \to \psi \circ \kappa$ of two transforms $\alpha : \phi \to \psi$ between functors $\phi, \psi : p \Longrightarrow q$ and $\beta : \xi \to \kappa$ between functors $\xi, \kappa : q \Longrightarrow r$ are calculated as *horizontal composition of two deformed quadrats*

$$\begin{array}{ccccc} & \phi & \searrow^\alpha & \xi & \searrow^\beta \\ p & \longrightarrow & q & \longrightarrow & r \\ \uparrow & & \uparrow & & \uparrow \\ p & \longrightarrow & q & \longrightarrow & r \\ & \psi & & \kappa & \end{array}$$

$$\alpha *^+ \beta = (\alpha \star \xi) \circ (\psi \star \beta) \ .$$

The associative composition of functors is essential for such application of deformed quadrats.

The edges set mapping of preevaluation functor $\mathbf{ev}_q : p \otimes q^p \Longrightarrow q$ is also expressed as horizontal composition of deformed quadrats. For an arrow $w \in p(u; v)$ we get a precontact arrow

$$\alpha_u \circ \psi(w) : \phi(u) \to \psi(v)$$

$$\begin{array}{ccccc} & u & \searrow^w & \psi & \searrow^\alpha \\ * & \longrightarrow & p & \longrightarrow & q \\ & & \uparrow & & \uparrow \\ * & \longrightarrow & p & \longrightarrow & q \\ & v & & \phi & \end{array}$$

So the *double evaluation* provides the horizontal composition of three deformed quadrats

$$\begin{array}{ccccccc} & u & & \searrow^w_\phi & \searrow^\alpha & \xi & \searrow^\beta \\ * & \longrightarrow & p & \longrightarrow & q & \longrightarrow & r \\ & & \uparrow & & \uparrow & & \uparrow \\ * & \longrightarrow & p & \longrightarrow & q & \longrightarrow & r \\ & v & & \psi & & \kappa & \end{array}$$

and coincides with the preevaluation of earlier defined prehorizontal composition.

Also we can express the same transform $\alpha : \phi \to \psi$ between parallel functors $\phi, \psi : p \Longrightarrow q$ with a *dual deformed quadrat*. Such duality was introduced in **G. Valiukevičius 1992** and thoroughly exploited in **G. Valiukevičius 2009**.

For dual diagram we take the identity arrows $\mathtt{Id}_p : p \Longrightarrow p$ and $\mathtt{Id}_q : q \Longrightarrow q$ as the *downarrows*.

$$\begin{array}{ccc} p & \xrightarrow{\phi} & q \\ \downarrow & & \downarrow \\ p & \xrightarrow{\psi} & q \\ & & \nearrow \end{array}$$

Such deformed quadrats we have called as cocontinuous quadrats.

A *posthorizontal composition* $\alpha * \beta : \phi \circ \xi \to \psi \circ \kappa$ of two transforms $\alpha : \phi \to \psi$ between functors $\phi, \psi : p \Longrightarrow q$ and $\beta : \xi \to \kappa$ between functors $\xi, \kappa : q \Longrightarrow r$ is calculated as *horizontal composition of two deformed quadrats*

$$\begin{array}{ccccc} p & \xrightarrow{\phi} & q & \xrightarrow{\xi} & r \\ \downarrow & & \downarrow & & \downarrow \\ p & \xrightarrow{\psi} & q & \xrightarrow{\kappa} & r \\ & \swarrow_\alpha & & \swarrow_\beta & \end{array}$$

$$\alpha *^- \beta = (\phi \star \beta) \circ (\alpha \star \kappa) \ .$$

The edges set mapping of postevaluation functor $\mathbf{ev}_q : p \otimes q^p \Longrightarrow q$ also is expressed as horizontal composition of deformed quadrats. For an arrow $w \in p(u; v)$ we get a postcontact arrow

$$\phi(w) \circ \alpha_v : \phi(u) \to \psi(v)$$

$$\begin{array}{ccccc} * & \xrightarrow{u} & p & \xrightarrow{\phi} & q \\ & & \downarrow & & \downarrow \\ * & \xrightarrow{v} & p & \xrightarrow{\psi} & q \\ & & \swarrow_w w & & \swarrow_\alpha \end{array}$$

So the *double evaluation* provides the horizontal composition of three deformed quadrats

$$\begin{array}{ccccccc} * & \xrightarrow{u} & p & \xrightarrow{\phi} & q & \xrightarrow{\xi} & r \\ & & \downarrow & & \downarrow & & \downarrow \\ * & \xrightarrow{v} & p & \xrightarrow{\psi} & q & \xrightarrow{\kappa} & r \\ & & \swarrow_w & & \swarrow_\alpha & & \swarrow_\beta \end{array}$$

and coincides with the postevaluation of defined posthorizontal composition.

For a natural transform $\beta : \xi \to \kappa$ the prehorizontal composition coincides with the posthorizontal composition. We need to check only coincidence of preevaluation with postevaluation

$$
\begin{array}{ccccccccc}
* & \xrightarrow{u} & p & \xrightarrow{\phi} & q & \xrightarrow{\xi} & r & & \\
 & \searrow^{w} & \uparrow & \searrow^{\alpha} & \uparrow & \searrow^{\beta} & \uparrow & & \\
* & \xrightarrow{v} & p & \xrightarrow{\psi} & q & \xrightarrow{\kappa} & r & & 
\end{array}
\quad = \quad
\begin{array}{ccccccc}
* & \xrightarrow{u} & p & \xrightarrow{\phi} & q & \xrightarrow{\xi} & r \\
 & & \downarrow & & \downarrow & & \downarrow \\
* & \xrightarrow{v} & p & \xrightarrow{\psi} & q & \xrightarrow{\kappa} & r \\
 & \swarrow_{w} & & \swarrow_{\alpha} & & \swarrow_{\beta} & 
\end{array}
$$

For two natural transforms $\alpha : \phi \to \psi$ and $beta : \xi \to \kappa$ we get a horizontal composition being natural transform

$$\alpha * \beta : \phi \circ \xi \to \psi \circ \kappa \ .$$

First we can check that the horizontal composition of dual deformed quadrats is dual for a horizontal composition of taken deformed quadrats. It is true not only for the pair of opposite identity functors, but also for arbitrary truly joint pair of functors as it was done in **G. Valiukevičius 1992** p. 104.

Then the conclusion for horizontal composition of natural transforms will be obvious.

We can see an instance of our proposition in metric spaces. The scalar space $[0, \infty]$ is compounded of positive real numbers $r \geq 0$ and infinite point $+\infty$. Taking inequality $r \geq s$ as an arrow $r \to s$, we get a graph structure. We need a new sign for it $[0, \infty]^+$. The scalar space will be a set of its vertexes

$$[0, \texttt{inf}]_0^+ = [0, \infty] \ .$$

The contact arrows mapping is provided by transitivity of inequalities

$$r \geq s, s \geq t \Longrightarrow r \geq t \ .$$

Such implication can be understood as a rule for definition of graphs transport over pullbacked product

$$[0, \infty]^+ \times_{[0, \infty]} [0, \infty]^+ \Longrightarrow [0, \infty]^+ \ .$$

Such contact arrow mapping defines associative composition of arrows as we have commuting diagram of graphs transports

$$
\begin{array}{ccc}
(p \geq q, q \geq r), r \geq s & \to & p \geq q, (q \geq r, r \geq s) \\
\downarrow & & \downarrow \\
p \geq r, r \geq s & & p \geq q, q \geq s \\
& \searrow \quad \swarrow & \\
& p \geq s &
\end{array}
$$

Commuting for such diagram is get automatically, as an inequality in ordered space means the unique arrow.

The inequality $r \geq r$ will be a unit arrow for such contact arrows mapping. It also can be expressed as two commuting diagrams of graphs transports

$$\begin{array}{ccc} r \geq s & \to & r \geq r, r \geq s \\ \searrow & & \downarrow \\ & & r \geq s \end{array} \qquad \begin{array}{ccc} r \geq s, s \geq s & \leftarrow & r \geq s \\ \downarrow & & \swarrow \\ r \geq s & & \end{array}$$

So indeed the graph $[0, \infty]^+$ is a category. The point $0 \in [0, \infty]$ is a final point for it as we have exactly one arrow from any point

$$r \geq 0 \uparrow r \in [0, \infty] \ .$$

It is true for arbitrary ordered space $X$ with the smallest member $0 \in X$. Taking inequality $r \geq s$ as an arrow, we get a category $X^+$ with final point $0 \in F$.

Every graphs transport between the categories provided by ordered spaces $f : X^+ \Longrightarrow Y^+$ defines monotone appointment of points

$$p \geq r \Longrightarrow f_0(p) \geq f_0(r) \ .$$

And conversely any monotone appointment of points defines graphs transport for categories provided of taken ordered spaces. Usually we don't distinguish ordered space $X$ from its provided category $X^+$.

We can check that every graphs transport $f : X^+ \Longrightarrow Y^+$ is exact functor between correspondent categories.

Indeed. Every unit arrow $p \geq p$ will be maintained by arbitrary graphs transport, as $f(p) \geq f(p)$. As edge between two points can be only unique, so we get maintained unit arrow.

Also for two *connected inequalities* $p \geq q$ and $q \geq r$ we get connected inequalities in the target space $f(p) \geq f(q)$ and $f(q) \geq f(r)$. Their composition $p \geq r$ is maintained also maintained $f(p) \geq f(q)$.

These properties are sufficient for graphs transport between categories provided by ordered spaces $f : X^+ \Longrightarrow Y^+$ to be an exact functor.

In the category $[0, \infty]^+$ we shall take a biproduct as a graphs transport

$$[0, \infty]^+ \times [0, \infty]^+ \Longrightarrow [0, \infty]^+$$

defined with appointment of vertexes provided by addition of real numbers. For edges we have appointment of inequality provided by monotone addition

$$(p \geq q, r \geq s) \Longrightarrow p + r \geq q + s \ .$$

Obviously such biproduct is commutative and associative. The final vertex $0 \in [0, \infty]$ will be neutral for such biproduct

$$0 + p = p = p + 0 \ .$$

So we have got a monoidal category $[0,\infty]^+$.

For a chosen number $s$ we get a functor $[0,\infty]^+ \Longrightarrow [0,\infty]^+$

$$p \geq r \Longrightarrow p + s \geq r + s \ .$$

For it we define another coadjoint exponential functor

$$p \geq r \Longrightarrow (p - s)^+ \geq (r - s)^+ \ .$$

The equivalence of inequalities

$$r + s \geq t \Longleftrightarrow r \geq (t - s)^+$$

provides the adjunction for true joint pair of functors $r + s$ and $t^s$. So we have got the closed monoidal category $[0,\infty]^+$

This category is complete for arbitrary Carte product

$$\prod_\alpha r_\alpha = \sup_\alpha r_\alpha \ .$$

The $[0,\infty]^+$-enriched graph $X$ is defined by the set of vertexes $X_0$ and the edges sets $X(x;y) \in [0,\infty]$ being positive real numbers. The $[0,\infty]^+$-graphs transport $f : X \Longrightarrow Y$ is defined by appointment of vertexes $f_0 : X_0 \to Y_0$ and edges set mappings

$$f_{x,y} : X(x;y) \to f_{x,y}(X(x;y)) \subset Y(f_0(x); f_0(y)) \ .$$

Obviously this expression is only symbolic. It is only the inequality between real numbers
$$X(x;y) \geq Y(f_0(x); f_0(y)) \ .$$
Such expressions remain mappings in the category of sets $Set$ and they can be useful for constructing some needed commuting diagrams which will be common for all $V$-enriched spaces with different values categories $V$.

The $[0,\infty]^+$-graph $X$ will be called as *distance space*, and $[0,\infty]^+$-graphs transport between $[0,\infty]^+$-enriched graphs $f : X \Longrightarrow Y$ can be called as *non-expanding mapping* of distance spaces.

The original $[0,\infty]^+$-graph needs to have chosen unit arrows

$$* \to 1_x(*) \subset X(x;x) \ .$$

It will be expressed as an inequality

$$0 \geq X(x;x) \uparrow x \in X_0 \ .$$

Also it needs to have chosen contact arrows mappings

$$X_{x,y,z} : X(x;y) \otimes X(y;z) \to X(x;y) \circ X(y;z) \subset X(x;z) \ .$$

It will be expressed as a triangular inequality for distance

$$X(x;y) + X(y;z) \geq X(x;z) \uparrow x,y,z \in X_0 \ .$$

The chosen unit arrows in original $[0,\infty]^+$-graphs automatically will be neutral for contact arrows mappings

$$X(x;y) \to 1_x(*) \circ X(x;y) = X(x;y) = X(x;y) \circ 1_y(*) \uparrow x,y \in X_0 \ ,$$

i. e. we have equalities between different mappings

$$X(x;y) \geq 0 + X(x;y) = X(x;y) = X(x;y) + 0 \ .$$

It must be expressed as commuting diagrams

$$\begin{array}{ccccccc} X(x;y) & \to & 1_x(*) \circ X(x;y) & X(x;y) \circ 1_y(*) & \leftarrow & X(x;y) \\ & \searrow & \downarrow & \downarrow & \swarrow & \\ & & X(x;y) & X(x;y) & & \end{array}$$

So in an original $[0,\infty]^+$-graph $X$ the distance is a true metric. Such graph will be called a metric space $X$. We can check that original $[0,\infty]^+$-graph is an $[0,\infty]^+$-enriched category, We need to check only commuting diagram for associativity of contact arrows mappings

$$X(x;y) + (X(y;z) + X(z;w)) = (X(x;y) + X(y;z)) + X(x;z) \ .$$

Such diagram is constructed for mapping over pullbacked product

$$S : X \times_{X_0} (X \times_{X_0} X) =\to (X \times_{X_0} X) \times_{X_0} X$$

$$\begin{array}{ccc} X(x;y) \otimes (X(y;z) \otimes X(z;u)) & \xrightarrow{s} & (X(x;y) \otimes_X (y;z)) \otimes X(z,;u) \\ 1_{X(x;y)} \otimes X_{y,z,u} \downarrow & & \downarrow X_{x,y,z} \otimes 1_{X(z,u)} \\ X(x;y) \otimes X(y;u) & & X(x;z) \otimes X(z;u) \\ X_{x,z,u} & \searrow \quad \swarrow & X_{x,z,u} \\ & X(x;u) & \end{array}$$

The commuting diagram in ordered space $[0,\infty]^+$ is get automatically, as there between two points can be only unique arrow. The nonexpanding mapping between metric spaces $f : X \implies Y$ will be crude functor of $[0,\infty]^+$-enriched categories. It again will be a natural functor. We need to check the commuting diagram for contact arrows mappings

$$\begin{array}{ccc} X(x;y) \otimes X(y;z) & \xrightarrow{f_{x,y} \otimes f_{y,z}} & Y(f(x);f(y)) \otimes Y(f(y);f(z)) \\ \downarrow & & \downarrow \\ X(x;y) \circ X(y;z) & \xrightarrow[f_{x,z}]{} & Y(f(x);f(y)) \circ Y(f(y);f(z)) \end{array}$$

For metric spaces $X$ and $Y$ we have triangular inequalities

$$X(x;y) + X(y;z) \geq X(x;z) \ ,$$

$$Y(f(x);f(y)) + Y(f(y);f(z)) \geq Y(f(x);f(z)) \ ,$$

and for nonexpanding mapping we have inequalities

$$X(x;y) \geq Y(f(x);f(y)) \ , \quad X(y;z) \geq Y(f(y);f(z)) \ ,$$

$$X(x;z) \geq Y(f(x);f(y)) \ .$$

Therefore we get also the inequality

$$X(x,y) + X(y;z) \geq Y(f(x);f(y)) + Y(f(z);f(z)) \ .$$

So we have all needed arrows for diagram and this diagram in ordered space $Y$ must be commuting

$$\begin{array}{ccc} X(x;y) + X(y;z) & \to & Y(f(x);f(y)) + Y(f(y);f(z)) \\ \downarrow & & \downarrow \\ X(x;z) & \to & Y(f(x);f(z)) \end{array}$$

Also it maintain chosen unit arrows

$$0 \geq X(x;x) \geq Y(f(x);f(x)) \ .$$

We have shown that every nonexpanding mapping $f : X \Longrightarrow Y$ between metric spaces is natural functor between correspondent $[0,\infty]^+$-enriched categories. So we have coinciding categories

$$OGrph_{[0,\infty]^+} = Cat^{=}_{[0,\infty]^+} \ .$$

In these categories for the joint pair of superfunctors we have a natural bijection. We shall check it directly.

**Proposition 10.** *It is equivalent properties for some mapping $\Phi : X \times A \to Y$ to be nonexpanding over taken tensor product of metric spaces $X \otimes A$*

$$X(x;x') + A(a,a') \geq Y(\Phi(x,a);\Phi(x',a'))$$

*and for the partial mapping $\Phi_1 : X \to Y^A$ to be nonexpanding for the metric of uniform boundary in functional space*

$$X(x;x') \geq \sup_{a \in A} Y(\Phi(x,a);\Phi(x',a')) \ .$$

**Proof:** Let we have the first inequality

$$X(x,x') + A(a;a') \geq \Phi(x,a);\Phi(x';a')) \ .$$

Then we get
$$X(x;x') \geq Y(\Phi(x,a); \Phi(x';a)) \uparrow a \in A ,$$
$$X(x;x') \geq \sup_{a \in A} Y(\Phi(x,a); \Phi(x',a')) .$$

Otherwise from the second inequality for arbitrary nonexpanding mapping $\Phi_1(x) \in Y^A$ we get inequality
$$X(x;x') \geq Y(\Phi(x,a); \Phi(x';a)) ,$$
and for the metric in the space $Y$ we can apply triangular inequality
$$Y(\Phi(x,a); \Phi(x';a')) \leq Y(\Phi(x;a); \Phi(x';a)) + Y(\Phi(x';a); \Phi(x';a')) \leq$$
$$X(x;x') + A(a;a') .$$

□

We shall formulate this proposition in the language of $[0,\infty]^+$-graphs. At first we define biproduct $X \otimes Y$ of $[O,\infty]^+$-graphs $X$ and $Y$. We take Carte tensor product of vertexes $(X \otimes Y)_0 = X_0 \times Y_0$ and tensor products of edges sets
$$X \otimes Y(\langle x,y \rangle, \langle x',y' \rangle) = X(x;x;) \otimes Y(y;y') .$$
For distance spaces we get $\ell_1$-distance on product space
$$X \otimes Y(\langle x,y \rangle, \langle x',y' \rangle) = X(x;x') + Y(y;y') .$$

With an assistant distance space $T$ we define the adjoint superfunctor. For a distance space $X$ it appoints the biproduct with the assistant space $X \otimes T$, and for a nonexpanding mapping between distance spaces $f : X \Longrightarrow Y$ it appoints the biproduct with the identity mapping of the assistant space
$$f \otimes \mathtt{Id}_T : X \otimes T \Longrightarrow Y \otimes T .$$

The coadjoint superfunctor will be define with the exponential space $Y^T$. Its vertexes set will be compound by all nonexpanding mappings between distance spaces $\phi : T \Longrightarrow Y$. The edges set between two nonexpanding mappings $\phi, \psi : T \Longrightarrow Y$ will be defined by Carte product in value category $[0,\infty]^+$
$$Y^T(\phi;\psi) = \prod_{s \in T_0} Y(\phi_0(s); \psi_0(s)) .$$

It is provided by upper boundary, therefore
$$Y^T(\phi;\psi) = \sup_{s \in T_0} Y(\phi_0(s); \psi_0(s))) .$$

Between two nonexpanding mappings we have got a distance of uniform convergence.

The exponential space $Y^T$ defines a superfunctor in the category of distance spaces $Grph_{[0,\infty]^+}$. For arbitrary nonexpanding mapping $f : Y \Longrightarrow Y'$ it appoints the changing of target space $f^T : Y^T \Longrightarrow (Y')^T$. For a nonexpanding mapping $\phi : T \Longrightarrow Y$ it appoints the new mapping $\phi \circ f : T \Longrightarrow Y'$ and between changed mappings it provides the new distance

$$(Y')^T(\phi \circ f; \psi \circ f) = \sup_{s \in T} Y'(f_0(\phi_0(s)); f_0(\psi_0(s))) \uparrow \phi, \psi \in Y^T .$$

Easy to check the maintenance of mappings composition

$$(f \circ g)^T = f^T \circ g^T .$$

It is true for general mappings in the category $Set$. We need only to check that a changing of target space with composition of nonexpanding mappings rests nonexpanding mapping. We have inequalities of distances for nonexpanding mappings

$$X(x; y) \geq Y(f_0(x); f_0(y)) , \quad Y(f_0(x); f_0(y)) \geq Y'(g_0(f_0(x)); g_0(f_0(y))) .$$

Therefore

$$\sup_{s \in T} Y'(g_0(f_0(\phi_0(s)))); g_0(f_0(\psi_0(s))))) \leq$$

$$\sup_{s \in T} Y(f_0(\phi_0(s)); f_0(\psi_0(s))) \leq$$

$$\sup_{s \in T} X(\phi_0(s); \psi_0(s)) .$$

Also the identity mapping defines the identity mapping of exponential spaces

$$(\text{Id}_Y)^T = \text{Id}_{Y^T} .$$

The mapping of unit arrow choosing

$$* \to 1_x(*) \subset X(x; x)$$

coincides with inequality for a distance

$$0 \geq X(x; x) .$$

Such inequalities for all points in the source space $x \in X$ and assistant space $s \in T$ allows us to define the sections transport of distance spaces

$$\lambda_X : X \Longrightarrow (X \otimes T)^T .$$

For a point $x \in X$ it appoints the section $\lambda_X(x) : T \Longrightarrow X \otimes T$ which for a point $s \in T$ appoints the couple $\langle x, s \rangle \in X \times T$. The section is nonexpanding mapping as in our case $X(x; x) = 0$ and we have the inequality

$$T(s; t) \geq X(x; x) + T(s; t) .$$

Also the sections transport is nonexpanding mapping as $T(s;s) = 0$

$$X(x;y) \geq \sup_{s \in T}(X(x;y) + T(s;s)) \ .$$

The sections transport provides the unit transform for joint pair of superfunctors. It will be a natural transformation. We have the commuting quadrats for arbitrary nonexpanding mapping of distance spaces $f : X \to X'$

$$\begin{array}{ccc} X & \xrightarrow{f} & X \\ \lambda_X \downarrow & & \downarrow \lambda_{X'} \\ (X \otimes T)^T & \xrightarrow[(f \otimes 1_T)^T]{} & (X' \otimes T)^T \end{array}$$

For an arrows set $X(x, x')$ we appoint the distance between sections mappings

$$(X \otimes T)^T(\lambda_X(x); \lambda_X(x')) = X(x; x')$$

and for a mapping $(f \otimes \text{Id}_T)^T : (X \otimes T)^T \to (X' \otimes T)^T$ we shall get the new sections mappings with distances

$$(X' \otimes T)^T(\lambda_{(X')}(f(x)); \lambda_{(X')}(f(x'))) = X'(f(x); f(x'))$$

the same as it would be a distance for sections mappings defined by the appointed points $f(x) \in X'$ and $f(x') \in X'$.

For a distance space $Y$ the contact arrows mappings

$$Y(x;y) \otimes Y(y;z) \to Y(x;y) \circ Y(y;z) \subset Y(x;z)$$

are expressed as inequalities

$$Y(x;y) + Y(y;z) \geq Y(x;z) \ .$$

In such spaces we define evaluation transport

$$\mathbf{ev}_Y : Y^T \otimes T \Longrightarrow Y \ .$$

It for a nonexpanding mapping $\phi : T \Longrightarrow Y$ and point in assistant distance space $s \in T$ appoints the value of taken mapping $\phi_0(s) \in Y$. It is also nonexpanding. For the pair of couples $\langle \phi, s \rangle \in Y^T \otimes T$ and $\langle \psi, t \rangle \in Y^T \otimes T$ we get the distance between values

$$Y(\phi_0(s); \psi_0(t)) \leq Y(\phi_0(s); \psi_0(s)) + \psi(T(s,t)) \leq \sup_{s \in T} Y(\phi_0(s); \psi_0(s)) + T(s;t) \ .$$

An evaluation transports provides the natural transform between superfunctors. We have the commuting quadrats for arbitrary nonexpanding mapping $g : Y \Longrightarrow Y'$

$$\begin{array}{ccc} (Y)^T \otimes T & \xrightarrow{g^T \otimes 1_T} & Y'^T \otimes T \\ \mathbf{ev}_Y \downarrow & & \downarrow \mathbf{ev}_{Y'} \\ Y & \xrightarrow[g]{} & Y' \end{array}$$

For a couple of nonexpanding mapping $\phi : T \Longrightarrow Y$ and point $s \in T$ we appoint the mapping's value $\phi(s) \in Y$, for a nonexpanding map $g : Y \Longrightarrow Y'$ we get the changed value $g(\phi(s)) \in Y'$.

For a pair of couples $\langle \phi, s \rangle$ and $\langle \psi, t \rangle$ with distance

$$Y^T(\phi; \psi) + T(s, t)$$

the evaluation $\mathbf{ev}_Y : Y^T \otimes T \Longrightarrow Y$ appoints the distance for appointed values $Y(\phi(s); \psi(t))$ and further the distance of changed values

$$Y'(g(\phi(s)); g(\psi(t))) \ .$$

It is the same as the distance calculated between two mappings $\phi : T \Longrightarrow Y$ and $\psi : T \Longrightarrow Y$ changed with nonexpanding mapping $g : Y \Longrightarrow Y'$

$$\langle (Y')^T (\phi \circ g; \psi \circ g) \ .$$

For the truly joint superfunctors $X \otimes A$ and $Y^A$ we need to check both triangular equalities. However we can apply our proposition about bijective name appointment.

For reversing realization appointment we need to check decomposition equality

$$X(a;b) \otimes T(s;t) \to f_{\langle a,s \rangle, \langle b,s \rangle}(X(a;b) \otimes 1_s(*)) \circ f_{\langle y,s \rangle, \langle y,t \rangle}(1_b(*) \otimes T(s;t)) =$$

$$f_{\langle x,s \rangle, \langle b,t \rangle}(X(a;b) \otimes T(s;t)) \subset Y(f(a,s); f(b,t))$$

In our case we have commuting diagram of inequalities

$$X(a;b) + T(s;t) \geq Y(f(a,s); f(b,s)) + Y(f)b,s); f(b,t)) \geq Y(f(a,s); f(b,t))$$

For the reversing the name appointment we need to apply the neutral equalities

$$X(a;b) \to g_{a,b}(X(a;b))_s \circ g(b)_{s,s}(1_s(*)) =$$

$$g_{a,b}(X(a;b))_s \subset Y(g(a)_s; g(b)_s) \ ,$$

$$T(s;t) \to g_{a,a}(1_a(*)) \circ g(a)_{s,t}(T(s;t)) =$$

$$g(a)_{s,t}(T(s;t)) \subset Y(g(a)_s; g(a)_t)$$

In our case we have such arrows only as inequalities. First we can see that nonexpanding mapping maintains the unit arrow mappings. For a nonexpanding mapping $g(b) : T \Longrightarrow X \otimes T$ we get the inequalities

$$0 \geq T(s;s) \geq Y(g(b)_s; g(b)_s) \ ,$$

and the first neutral equality

$$X(a;b) \geq Y(g(a)_s; g(b)_s) + Y(g(b)_s; g(b)_s) = Y(g(a)_s; g(b)_s) \ .$$

Also for a nonexpanding mapping $g : X \Longrightarrow Y^T$ we get the inequalities

$$0 \geq X(a;a) \geq \sup_{s \in T} Y(g(a)_s; g(a)_s) \ ,$$

and the second neutral equality

$$T(s;t) \geq Y(g(a)_s; g(a)_s) + Y(g(a)_s; g(a)_t) = Y(g(a)_s; g(a)_t) \ .$$

We can see that there isn't any need to have the contact arrows mapping in the source space $X$ or in the assistant space $T$.

We have another alternative metric space definition. In scalar space we have Carte product $r \vee s \in [0, \infty]$. It will be another tensor product. It defines a functor which for a scalar $r \in [0, \infty$ appoints another scalar $r \times s = r \vee s$. We have a coadjoint exponential functor which for the scalar $t \in [0, \infty]$ appoints another scalar $\mathtt{Hom}(s;t) \in [0, \infty]$

$$\mathtt{Hom}(a, v) = v \uparrow v > a; 0 \uparrow v \leq a \ .$$

We get a joint pair of functors $r \vee s$ and $\mathtt{Hom}(s;t)$.
For the unit transform we have an arrow $x \geq \mathtt{Hom}(a; x \vee a)$.
In the case $x > a$ we shall have $x \vee a > a$ and

$$x = x \vee a = \mathtt{Hom}(a; x \vee a) \ .$$

In another case $x \leq a$ we shall have $x \vee a \leq a$, therefore

$$\mathtt{Hom}(a; x \vee a) = 0 \leq x \ .$$

Also it is an arrow for the counit transform $\mathtt{Hom}(a, y) \vee a \geq y$.
In the case $y > a$ we shall have $\mathtt{Hom}(a; y) = y$ and

$$\mathtt{Hom}(a; y) \vee a = y \vee a \geq y \ .$$

In another case $y \leq a$ we shall have $\mathtt{Hom}(a; y) = 0$ and

$$\mathtt{Hom}(a; y) \vee a = a \geq y \ .$$

The triangular equalities are get automatically, therefore we have the true joint pair of functors $r \vee s$ and $\mathtt{Hom}(s;t)$ provided with assistant member $s \in [0, \infty]$.

So we have got a closed category again. It will be called *Cartesian closed category* in the scalar space $[0, \infty]$. We shall use a sign $[0, \infty]^+_{Carte}$ for it. Now we get a *translation* of monoidal categories

$$F : [0, \infty]^+ \Longrightarrow [0, \infty]^+_{Carte} \ .$$

Denoting $r \,\square\, s = r \vee s$ we have needed natural transform

$$t_{r,s} : (r \otimes s)^F \to r^F \,\square\, s^F \ ,$$

and identity arrow for coinciding neutral points $\bar{e} = e^F = 0 \in [0, \infty]^+_{Carte}$.

$$t_e : e^F \to \bar{e} \ .$$

We recognize the original $[0, \infty]^+_{Carte}$-graph $X$ as ultrametric space, i. e. we have inequalities for scalars $X(x; y)$

$$0 \geq X(x, y) \ ,$$

$$X(x; y) \vee X(y; z) \geq X(x; z) \ .$$

The former translation provides imbedding for the categories of original graphs

$$OGrph_{[0,\infty]^+_{Carte}} \subset OGrph_{[0,\infty]^+} \ .$$

The new tensor product of ultrametric spaces $X \mathbin{\square} T$ is defined with smaller distance

$$X(x; y) \mathbin{\square} A(a, b) \leq X(x; y) \otimes A(a; b) \ ,$$

but an exponential superfunctor appoints the same functional space. It is compounded of all nonexpanding mappings $\phi : T \to X$. And for two nonexpanding mappings $\phi : A \to Y$ and $\psi : A \to Y$ we have again the distance of uniform convergence

$$\mathtt{Hom}(T; Y)(\phi; \psi) = \sup_{a \in A} Y(\phi(a); \psi(a)) \ .$$

This tensor product $X \mathbin{\square} Y$ can be defined for arbitrary metric spaces $X$ and $Y$ defined with the initial category $Grph_{[0,\infty]}$,

In the case of ultrametric target space $Y$ we again have equivalent properties for a mapping $\Phi : X \times A \to Y$ from arbitrary metric spaces $X$ and $A$ to be nonexpanding over the new tensor product of ultrametric spaces

$$X(x; x') \vee A(a; a') \geq Y(\Phi(x, a); \Phi(x'; a'))$$

and for the partial mapping $\Phi_1 : X \to Y^A$ to be nonexpanding with the metric of uniform boundary

$$X(x; x') \geq \sup_{a \in A} Y(\Phi(x, a); \Phi(x'; a)) \ .$$

**Proposition 11.** *For arbitrary metric spaces $X$ and $A$ and ultrametric space $Y$ we have equivalent properties for mapping $\Phi : X \times A \to Y$ to be nonexpanding over the new tensor product of ultrametric spaces*

$$X(x; x') \vee A(a; a') \geq Y(\Phi(x, a); \Phi(x'; a'))$$

*and for the partial mapping $\Phi_1 : X \to Y^A$ to be nonexpanding with the metric of uniform boundary*

$$X(x; x') \geq \sup_{a \in A} Y(\Phi(x, a); \Phi(x'; a)) \ .$$

**Proof:** Let we have the first inequality

$$X(x, x') \vee A(a; a') \geq Y(\Phi(x, a); \Phi(x'; a')) \ .$$

Then we get

$$X(x; x') \geq Y(\Phi(x, a); \Phi(x'; a)) \uparrow a \in A \ ,$$

$$X(x; x') \geq \sup_{a \in A} Y(\Phi(x, a); \Phi(x', a')) \ .$$

Otherwise from the second inequality for arbitrary nonexpanding mappings $\Phi_1(x) \in Y^A$ we get the inequality

$$X(x; x') \geq Y(\Phi(x, a); \Phi(x'; a)) \ ,$$

and for the ultrametric in space $Y$ we can apply triangular inequality

$$Y(\Phi(x, a); \Phi(x'; a')) \leq Y(\Phi(x; a); \Phi(x'; a)) \vee Y(\Phi(x'; a); \Phi(x'; a')) \leq$$

$$X(x; x') \vee A(a; a') \ .$$

$\square$

Seeking to understand this mathematical fact, I have gone astray for one year in the theory of enriched graphs. For the better illustration we could work with more general values category `Ban` of Banach spaces.

We shall generalize the notion of monoidal category $V$. Let we have defined some biproduct as bifunctor

$$\otimes : V \times V \Longrightarrow V$$

with neutral point $e \in V_0$. For it we have natural transforms for left canceling defined with right projection

$$r_a : e \otimes a \to a$$

and another transform for right canceling defined with left projection

$$\ell_a : a \otimes e \to a \ .$$

Usually we shall take a category's $V$ final object $* \in V_0$ as a neutral point for taken biproduct.

Also we can ask a natural transform for associativity defined with arrows

$$s_{a,b,c} : a \otimes (b \otimes c) \to (a \otimes b) \otimes c \ .$$

So we get the classical monoidal category if such natural transforms are defined with isomorphic arrows.

Demanding the pentagon diagram

$$\begin{array}{ccccc}
a \otimes (b \otimes (c \otimes d)) & \xrightarrow{s} & (a \otimes b) \otimes (c \otimes d) \\
{\scriptstyle 1 \otimes s} \swarrow & & \searrow {\scriptstyle s \otimes 1} \\
a \otimes ((b \otimes c) \otimes d) \xrightarrow{s} & (a \otimes (b \otimes c)) \otimes d & \xrightarrow{s \otimes 1} ((a \otimes b) \otimes c) \otimes d
\end{array}$$

the triangular diagram

$$\begin{array}{ccc} a \otimes (e \otimes c) & \stackrel{s}{\longrightarrow} & (a \otimes e) \otimes c \\ 1 \otimes \ell \searrow & & \swarrow r \otimes 1 \\ & a \otimes c & \end{array}$$

and equal arrows
$$\ell_e = r_e : e \otimes e \to e$$
for classical monoidal category we get a tensor product. By coherence theorem in **S. MacLane 1971** part VII we get a unique canonical isomorphic arrow between tensor products defined for arbitrary two binary words of the same length.

Later for the working with $V$-enriched spaces we shall ask a natural transform of *relator's involution*

$$m_{p,r,s,t} : (p \otimes s) \otimes (r \otimes t) \to (p \otimes r) \otimes (s \otimes t) \ .$$

For associative tensor product this property is weaker than commutativity property.

Also we demand the *upper continuity* for the Carte products

$$\prod_i p_i \otimes \prod_i r_i = \prod_i p_i \otimes r_i \ .$$

Let in the category $V$ we have another biproduct $r \square s$ which is a relator for the first tensor product $r \otimes s$ with a natural transform defined with arrows

$$N_{p,r,s,t} : (p \square s) \otimes (r \square t) \to (p \otimes r) \square (s \otimes t)$$

and the same neutral member with natural cancellation transforms defined by arrows
$$* \square r \to r \leftarrow r \square * \ .$$
We shall demand that these canceling arrows will be isomorphisms
$$* \square r = \to r \leftarrow = r \square * \ .$$
Then we have inclusion arrows
$$r \otimes s \to (r \square *) \otimes (s \square *) \to (r \square s) \otimes (* \square *) \to r \square s \ .$$
which will be *comparison arrows* between two biproducts
$$t_{r,s} : r \otimes s \to r \square s \ .$$
Translation of values category is defined by transport of monoidal categories $F : V \Longrightarrow V'$ with comparison arrows
$$t_{a,b} : a^F \square b^F \to (a \otimes b)^F \ ,$$

$$t_e : \bar{e} \to e^F \ .$$

We get a *continuous translation* of values category $V$ to the new values category $V'$ with the same points $r \in V_0$ and arrows $\alpha \in V$ but with a new biproduct $x \,\square\, s$

$$F : V \Longrightarrow V' \ .$$

For objects and arrows we have identity functor

$$F_0(r) = r \ , \quad \alpha^F = \alpha \uparrow r \in V_0, \alpha \in V(r;s) \ .$$

Such translation will be continuous for biproducts

$$\begin{array}{ccc} V \times V & \xrightarrow{F \times F} & V' \times V' \\ \otimes \downarrow & & \downarrow \square \\ V & \xrightarrow{F} & V' \end{array} \quad \nearrow$$

i. e. we have continuous square with natural transform of comparison defined with arrows

$$t_{r,s} : (r \otimes s)^F \to r^F \,\square\, s^F$$

Also we ask the continuity for chosen neutral member $e^F \to \bar{e}$

$$\begin{array}{ccc} * & = & * \\ e \downarrow & & \downarrow \bar{e} \\ V & \xrightarrow{F} & V' \end{array} \quad \nearrow$$

If the new biproduct is associative with natural transform defined by arrows

$$s'_{a,b,c} : a \,\square\, (b \,\square\, c) \to (a \,\square\, b) \,\square\, c \ ,$$

we should demand that two continuous quadrats would be coherent for associativity transforms. Such continuous quadrats are defined for a functor $F : V \Longrightarrow V'$ with *convergences* defined by bifunctors of different biproducts $V \times V \Longrightarrow V \otimes V$ and $V' \times V' \Longrightarrow V' \,\square\, V'$

$$\begin{array}{ccc} V \times (V \times V) & \xrightarrow{F \times (F \times F)} & V' \times (V' \times V') \\ \downarrow & & \downarrow \\ V \times (V \otimes V) & \xrightarrow{F \times F} & V' \times (V' \,\square\, V') \\ & & \nearrow \\ \downarrow & & \downarrow \\ V \otimes (V \otimes V) & \xrightarrow{F} & V' \,\square\, (V' \,\square\, V') \\ & & \nearrow t \end{array}$$

and

$$\begin{array}{ccc} (V \times V) \times V & \xrightarrow{(F \times F) \times F} & (V' \times V') \times V' \\ \downarrow & & \downarrow \\ (V \otimes V) \times V & \xrightarrow{F \times F} & (V' \,\square\, V') \times V' \\ & & \nearrow \\ \downarrow & & \downarrow \\ (V \otimes V) \otimes V & \xrightarrow{F} & (V' \,\square\, V') \,\square\, V' \\ & & \nearrow t \end{array}$$

The natural transform between two coherent graphs is defined with vertexes in the same spaces and arrows of transform between them. It can be explain by graphs of continuous mapping for two coherent convergences. In our case the convergences are defined by coherent biproducts.

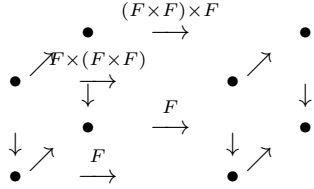

The coherence can be defined by commuting diagram in the target space

$$(V \otimes (V \otimes V))^F \xrightarrow{s^F} ((V \otimes V) \otimes V)^F$$
$$t \downarrow \qquad\qquad \downarrow t$$
$$V' \square (V' \square V') \xrightarrow{s'} (V' \square V') \square V'$$

For objects $a, b, c \in V$ we have equality of arrows compositions

$$s^F_{a,b,c} \circ t_{(a,b),c} = t_{a,(b,c)} \circ s'_{a^F, b^F, c^F} : (a \otimes (b \otimes c))^F \to (a^F \square b^F) \square c^F .$$

We have denoted arrows of horizontal composition for natural transforms

$$t_{a,(b,c)} = (1_a \otimes t_{b,c}) \circ t_{a, b \otimes c} ,$$

$$t_{(a,b),c} = (t_{a,b} \otimes 1_c) \circ t_{a \otimes b, c} .$$

For the inclusion of the sets $a \otimes b \subset a \square b$ it demands that a mapping over the larger set would be restricted to the mapping over the smaller set

$$\begin{array}{ccc} (a \otimes (b \otimes c) & \subset & a \square (b \square c) \\ \downarrow & & \downarrow \\ (a \otimes b) \otimes c & \subset & (a \square b) \square c \end{array}$$

We shall say that translator $F : V \Longrightarrow V'$ is *coherent for associativity transforms*.

In this case the coherence theorem for classical monoidal categories will produce a restricted natural isomorphism for arbitrary binary word with the smaller tensor product.

We want to explain why equality of natural transforms in the *final space* define the coherence of the whole quadrat. We can see that earlier defined equality is a unique choice. We propose the procedure of reduction which can be applied also for more complex coherent *diagrams* .

Let we have the arrows of natural transform for every vertex of diagram

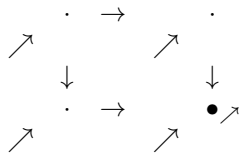

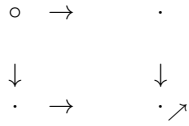

We shall study all possible *circuits* from initial vertex in the first quadrat to the final vertex in the second quadrat. We can "compress" two commuting quadrats in such diagram. This process I have drawn by hand in **G. Valiukevičius 1992** p. 123. Now such drawing must be done by computer, and I haven't time for it at this moment.

As the result we get a problem of coherent transforms for continuous triangle

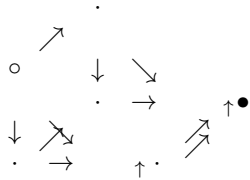

We have drawn also another transported arrow in the final category. It provides a commuting lower side. So in the final category we get a commuting quadrat which is needed and provides coherent transform for the whole diagram.

For more detailed work it is possible to apply another properties of biproducts. We notice that we are working with comparison arrows for continuous translation. The opposite arrows would define a cocontinuous translation, and for it the product must be changed by coproducts.

The biproduct $\otimes : V \times V \Longrightarrow V$ can be considered as a coadjoint functor for the diagonal functor $\triangle : V \Longrightarrow V \times V$. The *diagonal functor* for a object $r \in V_0$ appoints the couple $\langle r, r \rangle \in V_0 \times V_0$ and for an arrow $u : r \to s$ appoints the couple of arrows

$$\langle u, u \rangle : \langle r, r \rangle \to \langle s, s \rangle \ .$$

The unit transform will be defined by *diagonal arrows*

$$i_r : r \to r \otimes r \ .$$

Such biproduct will be called *diagonalized*. The counterexample can be defined with non Hausdorff topology on the product space $X \times Y$.

In **G. Valiukevičius 2009** p.107 the proposition 3c concludes that comparison from diagonalized biproduct provides another contact arrow in the target category $V'$

$$t' = i^F \circ t : r^F \to r^F \otimes r^F \ .$$

This was formulated as changing the initial point of diagram.

$$
\begin{array}{ccc}
① \quad V \times V & \xrightarrow{F \times F} & V' \times V' \\
\quad \triangle \updownarrow \otimes & \xrightarrow{F} & \downarrow \square' \\
② \quad V & & V' \\
& & \nearrow
\end{array}
$$

We shall say that we have got contact arrows *transported by the translation* $F : V \Longrightarrow V'$. If we have diagonal arrows also in the target category $V'$, then translation $F : V \Longrightarrow V'$ will be called *diagonalized* if the diagonal arrows in the target category $V'$ extends the transported diagonal arrow.

The counit transform will be defined with pair of projections

$$\langle p^1_{r,s}, p^2_{r,s} \rangle : \langle r \otimes s, r \otimes s \rangle \to \langle r, s \rangle .$$

Such biproduct will be called projectable.

In **G. Valiukevičius 1992** p. 171 we provide an example of non projectable biproduct for graphs. In **G. Valiukevičius 2009** p. 107 the proposition 3a concludes that comparison from projectable biproduct provides another projections in the target category $V'$

$$t' = \langle t, t \rangle \circ \langle p^1, p^2 \rangle .$$

This was formulated as changing the final point of diagram

$$
\begin{array}{ccc}
① & & \nwarrow \\
V \times V & \xrightarrow{F \times F} & V' \times V' \\
\otimes \downarrow & & \triangle \updownarrow \square' \\
V & \xrightarrow{F} & V' \\
& & \nearrow
\end{array}
$$

We shall call the projections $p^1 : r^F \otimes s^F \to r^F$ and $p^2 : r^F \square s^F \to s^F$ *transported by translation* $F : V \Longrightarrow V'$. If we have projections in the target category $V'$, the translation $F : V \Longrightarrow V'$ will be called *projectable* if the projections in target category $V'$ extends the transported projections.

We shall say that a biproduct $r \otimes s$ in a category $V$ is *almost coadjoint* one, if it is projectable and diagonalized. Accordingly the translation $F : V \Longrightarrow V'$ transports the almost coadjoint biproduct if it is projectable and diagonalized.

The first triangular equality for a biproduct $r \otimes s$

$$\langle 1_r, 1_r \rangle =: \langle r, r \rangle \xrightarrow{\langle i_r, i_r \rangle} \langle r \otimes r, r \otimes r \rangle \xrightarrow{\langle p_1, p_2 \rangle} \langle r, r \rangle$$

provides the injective name appointment. We express that there is only one couple of arrows $p \to r$ and $p \to s$ providing an arrow to this biproduct $p \to r \otimes s$.

The faithful almost coadjoint translation $F : V \Longrightarrow V'$ reflects the property of target category $V'$ to have injective name appointment.

The second triangular equality for a biproduct $r \otimes s$

$$1_{r \otimes s} =: r \otimes s \stackrel{i_{r \otimes s}}{\longrightarrow} \langle r \otimes s, r \otimes s \rangle \stackrel{p^1 \otimes p^2}{\longrightarrow} r \otimes s$$

provides the injective realization appointment. We express that there is only one arrows to the biproduct $p \to r \otimes s$ providing the couple of arrows $p \to r$ and $p \to s$.

The faithful almost coadjoint translation $F : V \Longrightarrow V'$ reflects the property of target category $V'$ to have an injective realization appointment.

The almost coadjoint biproduct with both first and second triangular equalities will be *truly coadjoint* one, therefore it coincides will the Carte product $r \times s$ up unique isomorphism.

These properties were applied in **G. Valiukevičius 1992** part 5.4 to study the commuting diagrams in categories with abstract products. They also will be useful for any detailed work with concrete categories.

The most special translation of monoidal categories is defined with isomorphic comparison arrows

$$t_{r,s} : (r \otimes s)^F =\to r^F \square s^F .$$

We shall call it a *strict translation*. Such property is meaningful even without the property of projectable tensor product.

In monoidal category $V$ the comparison arrow $r \otimes s \to r \square s$ provides the restriction of arbitrary arrow over the new biproduct $r \square s \to t$. It is get as composition

$$r \otimes s \to r \square s \to t .$$

We can take the same exponential functor $t^s$ as coadjoint functor for the functor defined with the new biproduct $r \square s$. We get again a sections arrow

$$\lambda'_r : r \to (r \otimes s)^s \to (r \square s)^s .$$

However we can't get a truly joint pair again. We can't have evaluation arrows $\mathbf{ev}_t : t^s \square s \to t$ for arbitrary member $t \in V_0$. It may be necessary also to use a smaller exponential space to get possibility to define on it the evaluation arrow.

The same situation will be for the $V$-graphs. In the category $Grph_V$ we have a joint pair of superfunctors. An adjoint functor of biproduct $X \otimes A$ is provided with the tensor product $r \otimes s$ in the values category $V$ and a coadjoint functor will be the exponential functor $Y^A$ provided by transforms between $V$-graphs transports $\phi : A \Longrightarrow Y$. The new biproduct of $V$-graphs $X \square A$ is defined by the new biproduct $r \square s$ in the values category $V$. We can work with another values category $V'$ and define the category of $V'$-graphs $Grph_{V'}$. Every $V'$-graphs transport $\phi : A \Longrightarrow Y$ will be also the $V$-graphs transport. So the new category will be imbedded in the former category of $V$-graphs

$$Grph_{V'} \subset Grph_V .$$

However the exponential space $(Y^A)'$ for $V'$-graph $Y$ will be smaller and imbedded in the former exponential space from category $Grph_V$. We expose such remarks more carefully.

Within the category of $V$-graphs $Grph_V$ we can define another biproduct $X \square Y$. For vertexes sets we take the same Carte tensor product $X_0 \times Y_0$ and for edges sets we take the $\square$-biproducts in the values category $V$

$$X \square Y(\langle x,y \rangle; \langle x',y' \rangle) := X(x;x') \square Y(y;y') .$$

For $V$-graphs transports $f : X \Longrightarrow X'$ and $g : Y \Longrightarrow Y'$ we define a new biproduct

$$f \square g : X \square Y \Longrightarrow X' \square Y'$$

having the new edges set mappings

$$(f \square g)_{\langle x,y \rangle, \langle x',y' \rangle} = f_{x,x'} \square g_{y,y'} :$$

$$X(x;x') \square Y(y;y') \to X'(f_0(x); f_0(x')) \square Y'(g_0(y); g_0(y')) .$$

We have got a new bisuperfunctor in the category of $V$-graphs $Grph_V$. We additionally can check that we have got a bisuperfunctor also in smaller categories of original $V$-graphs or in the category of $V$-enriched categories with continuous functors $Cat_V^+$, cocontinuous functors $Cat_V^-$ or exact functors $Cat_V^=$.

Indeed. We take the unit arrow mappings transported with such bisuperfunctor. For the unit arrow mappings $1_x : * \to 1_x(*) \subset X(x;x)$ and $1_y : * \to 1_y(*) \subset Y(y;y)$ it appoints the biproduct

$$1_{\langle x,y \rangle} = 1_x \square 1_y : * \to 1_x(*) \square 1_y(*) \subset X(x;x) \square Y(y;y) .$$

If these $V$-graphs transports have maintained the unit arrow mappings

$$1_x \circ f_{x,x} = 1_{f_0(x)} , \quad 1_y \circ g_{y,y} = 1_{g_0(y)} ,$$

then we get maintained unit arrows mappings also in the biproduct of $V$-graphs $X \square Y$

$$(1_x \square 1_y) \circ (f_{x,x} \square g_{y,y}) = (1_x \circ f_{x,x}) \square (1_y \circ g_{y,y}) = f_{f_0(x)} \square 1_{g_0(y)} .$$

The contact arrows mappings are also transported by taken bisuperfunctor. In the new biproduct of $V$-graphs $X \square Y$ the contact arrows mappings are defined with the $\square$-product of contact arrows mappings for each component

$$(X \square Y(\langle x,y \rangle; \langle x',y' \rangle)) \otimes (X \square Y(\langle x',y' \rangle; \langle x'',y'' \rangle)) =$$

$$(X(x;x') \square Y(y;y')) \otimes (X(x';x'') \square Y(y';y'')) \to$$

$$(X(x;x') \otimes X(x';x'')) \square (Y(y;y') \otimes Y(y';y'')) \to$$

$$(X(x;x') \circ X(x';x'')) \square (Y(y;y') \circ Y(y';y'')) \subset$$

$$X(x;x'') \square Y(y;y'') = X \square Y(\langle x,y \rangle; \langle x'',y'' \rangle) .$$

The equalities for the chosen unit arrow mappings

$$1_x(*) \circ X(x;x') = X(x;x') = X(x;x') \circ 1_{x'}(*) ,$$

$$1_y(*) \circ Y(y;y') = Y(y;y') \circ 1_{y'}(*)$$

provides correspondent equalities in the new biproduct of $V$-graphs $X \square Y$. We shall write only for the first equalities

$$1_{\langle x,y \rangle}(*) \otimes (X \square Y(\langle x,y \rangle; \langle x',y' \rangle)) =$$

$$(1_x(*) \square 1_y(*)) \otimes (X(x;x') \square Y(y;y') \to$$

$$(1_x(*) \otimes X(x;x')) \square (1_y(*) \otimes Y(y;y')) \to$$

$$(1_x(*) \circ X(x;x')) \square (1_y(*) \circ Y(y;y')) =$$

$$X(x;x') \square Y(y;y') =$$

$$X \square Y(\langle x,y \rangle; \langle x',y' \rangle) .$$

So for original $V$-graphs $X$ and $Y$ we get the original $V$-graph $X \square Y$ again. We have got that bisuperfunctor $X \otimes Y$ can be restricted in the category of crude functors between original $V$-graphs $OGrph_V$.

For a natural *relator's transform* defined with arrows

$$n_{r,r',s,s'} : (r \square r') \otimes (s \square s') \to (r \otimes s) \square (r' \otimes s')$$

we can prove that the new bisuperfunctor also maintains natural functors, i. e. we have commuting diagram again for the contact arrows mappings. First we have commuting diagram

$$(X(x;x') \square Y(y;y')) \otimes (X(x';x'') \square Y(y';y'')) \to$$
$$(f_{x;x'} \square g_{y;y'}) \otimes (f_{x';x''} \square g_{y';y''}) \downarrow$$
$$(X'(f_0(x); f_0(x')) \square Y'(g_0(y); g_0(y')) \otimes (X'(f_0(x'); f_0(x'')) \square Y'(g_0(x'); g_0(y''))) \to$$

$$(X(x;x') \otimes X(x';x'')) \square (Y(y;y') \otimes Y(y';y''))$$
$$\downarrow (f_{x,x'} \otimes f_{x',x''}) \square (g_{y,y'} \otimes g_{y',y''})$$
$$(X'(f_0(x); f_0(x')) \otimes X'(f_0(x'); f_0(x''))) \square (Y'(g_0(y); g_0(y')) \otimes Y'(g_0(y'); g_0(y'')))$$

Having commuting diagrams for a natural functor $f : X \Longrightarrow X'$

$$\begin{array}{ccc} X(x;x') \otimes X(x';x'') & \xrightarrow{\circ} & X(x;x'') \\ f_{x,x'} \otimes f_{x';x''} \downarrow & & \downarrow f_{x;x''} \\ X'(f_0(x); f_0(x')) \otimes X'(f_0(x'); f_0(x'')) & \xrightarrow{\circ} & X'(f_0(x); f_0(x'')) \end{array}$$

and for another exact functor $g: Y \Longrightarrow Y'$

$$
\begin{array}{ccc}
Y(y;y') \otimes Y(y';y'') & \xrightarrow{\circ} & Y(y;y'') \\
{\scriptstyle g_{y,y'} \otimes g_{y';y''}} \downarrow & & \downarrow {\scriptstyle g_{y;y''}} \\
Y'(g_0(y); g_0(y')) \otimes Y'(g_0(y'); g_0(y'')) & \xrightarrow{\circ} & Y'(g_0(y); g_0(y''))
\end{array}
$$

we get commuting diagram for the new biproducts of arrows

$$(X(x;x') \otimes X(x';x'')) \square (Y(y;y') \otimes Y(y';y''))$$
$$(f_{x,x'} \otimes f_{x';x''}) \square (g_{y,y'} \otimes g_{y';y''}) \downarrow$$
$$(X'(f_0(x); f_0(x')) \otimes X'(f_0(x'); f_0(x''))) \square (Y'(g_0(y); g_0(y')) \otimes Y'(g_0(y'); g_0(y'')))$$

$$
\begin{array}{ccc}
\rightarrow & X(x;x'') \square Y(y;y'') \\
 & \downarrow f_{x;x''} \square g_{y,y''} \\
\rightarrow & X'(f_0(x); f_0(x'')) \square Y'(g_0(y); g_0(y''))
\end{array}
$$

and finally the commuting diagram for the biproduct of functors

$$(X(x;x') \square Y(y;y')) \otimes (X(x';x'') \square Y(y';y'')) \rightarrow$$
$$(f_{x,x'} \square g_{y,y'}) \otimes (f_{x',x''} \square g_{y',y''}) \downarrow$$
$$(X'(f_0(x); f_0(x')) \square Y'(g_0(y); g_0(y'))) \otimes (X'(f_0(x'); f_0(x'')) \square Y'(g_0(x'); g_0(y''))) \rightarrow$$

$$
\begin{array}{ccc}
\rightarrow & X(x;x'') \square Y(y;y'') \\
 & \downarrow f_{x,x''} \square g_{y,y''} \\
\rightarrow & X'(f_0(x); f_0(x'')) \square Y'(g_0(y); g_0(y''))
\end{array}
$$

So we have shown that taken bisuperfunctor can be restricted in the category of original $V$-graphs with natural functors $OGrph_V^{\overline{=}}$.

For weaker assumption about arrows $n_{r,r',s,s'}$, we can prove restriction of taken bisuperfunctor in larger partial categories. For cocontinuous arrows $n_{r,r',s,s'}$ we get restriction in the category of original $V$-graphs with continuous functors $OGrph_V^+$, and for continuous arrows $n_{r,r',s,s'}$ we get restriction in the category of cocontinuous functors $OGrph_V^-$. The earlier proving for tensor product $X \otimes Y$ will be considered as a partial case.

Finally we shall study biproducts in the category of $V$-enriched categories. Let the new biproduct in values category $V$ is also associative with natural isomorphisms

$$s'_{r,s,t}: r \square (s \square t) =\!\!\Rightarrow (r \square s) \square t$$

and comparison arrows

$$t_{p,r}: p \otimes r \rightarrow p \square r$$

are coherent for associativity isomorphisms, i. e. we have commuting diagram

$$s^F_{a,b,c} \circ t_{(a,b),c} = t_{a,(b,c)} \circ s'_{a^F, b^F, c^F}: (a \otimes (b \otimes c))^F \rightarrow (a^F \square b^F) \square c^F.$$

We shall use relator's involution, which is defined with isomorphic arrows

$$m_{p,r,s,t} : (p \otimes s) \otimes (r \otimes t) =\to (p \otimes r) \otimes (s \otimes t),$$

The contact arrows mappings

$$X_{a,b,c} : X(a,b) \otimes X(b;c) \to X(a;c)$$

defines a functor over pullbacked tensor product

$$\circ : X \otimes_p X \Longrightarrow X.$$

Repeating relator's mapping we get isomorphisms

$$(X \otimes Y) \otimes_p ((X \otimes Y) \otimes_p (X \otimes Y)) =\to$$

$$(X \otimes p(X \otimes_p X) \otimes (Y \otimes_p (Y \otimes_p Y)),$$

$$((X \otimes Y) \otimes_p (X \otimes Y)) \otimes_p (X \otimes Y) =\to$$

$$((X \otimes_p (X \otimes_p X)) \otimes ((Y \otimes_p Y) \otimes_p Y).$$

Therefore the associativity transforms for each component

$$X \otimes_p (X \otimes_p X) \to (Y \otimes_p X) \otimes_p X,$$

$$Y \otimes_p (Y \otimes_p Y) \to (Y \otimes_p Y) \otimes_p Y$$

provides the transform for the associativity of tensor product

$$(X \otimes Y) \otimes_p ((X \otimes Y) \otimes_p (X \otimes Y)) \to ((X \otimes Y) \otimes_p (X \otimes Y)) \otimes_p (X \otimes Y).$$

There isn't any obvious way to define associativity transform for the new biproduct of $V$-categories. We can assume that relator's arrows

$$n_{p,r,s,t} : (p \,\square\, s) \otimes (r \,\square\, t) =\to (p \otimes r) \,\square\, (s \otimes t)$$

are get from involution arrows

$$m_{p,r,s,t} : (p \otimes s) \otimes (r \otimes t) \to (p \otimes r) \otimes (s \otimes t),$$

by some topological extension of continuous mapping. So the associativity transform can be provided also by topological extension of continuous mapping

$$X \,\square\, Y \otimes_p (X \,\square\, Y \otimes_p X \,\square\, Y) \to (X \,\square\, Y \otimes_p X \,\square\, Y) \otimes_p X \,\square\, Y.$$

Then in the new product of $V$-graphs defined contact arrows mappings will again provide an associative composition, These contact arrows mappings will define the functor over pullbacked product

$$(X \,\square\, Y) \otimes_p (X \,\square\, Y) \Longrightarrow X \,\square\, Y$$

as composition of relator's transform with contact arrows mappings for each component $X$ and $Y$

$$(X \square Y) \otimes_p (X \square Y) \xrightarrow{N} (X \otimes_p X) \square (Y \otimes_p Y) \xrightarrow{\circ \square \circ} X \square Y \ .$$

In each $V$-enriched category we have commuting diagram for associative composition

$$\begin{array}{ccc}
X \otimes_p (X \otimes_p X) & \xrightarrow{S} & (X \otimes_p X) \otimes_p X \\
{}_{1_X \otimes \circ} \downarrow & & \downarrow {}_{\circ \otimes 1_X} \\
X \otimes_p X & & X \otimes X \\
{}_{\circ} \searrow & & \swarrow {}_{\circ} \\
& X &
\end{array}$$

and

$$\begin{array}{ccc}
Y \otimes_p (Y \otimes_p Y) & \xrightarrow{S} & (Y \otimes_p Y) \otimes_p Y \\
{}_{1_Y \otimes \circ} \downarrow & & \downarrow {}_{\circ \otimes 1_Y} \\
Y \otimes_p Y & & Y \otimes_p Y \\
{}_{\circ} \searrow & & \swarrow {}_{\circ} \\
& Y &
\end{array}$$

We shall check such commuting diagram for the contact arrows mappings from the biproduct $X \square Y$. It was for the tensor product $X \otimes Y$, i. e. we have commuting *double diagram.* For them we get generated commuting diagrams. For example with commuting triangles

$$\begin{array}{ccc} X & \to & Y \\ & \searrow \downarrow & \\ & Z & \end{array} \qquad \begin{array}{ccc} X' & \to & Y' \\ & \searrow \downarrow & \\ & Z' & \end{array}$$

we get also the commuting triangle of biproducts

$$\begin{array}{ccc} X \otimes X' & \to & Y \otimes Y' \\ & \searrow \downarrow & \\ & Z \otimes Z' & \end{array}$$

So we can get commuting diagrams of arrows composition in the tensor product $X \otimes Y$. This diagram can be changed to the double diagram using relator's involution

$$\begin{array}{ccc}
(X \otimes Y) \otimes_p ((X \otimes Y) \otimes_p (X \otimes Y)) & \xrightarrow{S} & ((X \otimes Y) \otimes_p (X \otimes Y)) \otimes_p (X \otimes Y) \\
{}_{1_{X \otimes Y} \otimes \circ} \downarrow & & \downarrow {}_{\circ \otimes 1_{X \otimes Y}} \\
(X \otimes Y) \otimes_p (X \otimes Y) & & (X \otimes Y) \otimes_p (X \otimes Y) \\
{}_{\circ} \searrow & & \swarrow {}_{\circ} \\
& X \otimes Y &
\end{array}$$

Now we get again a commuting diagram for topological extended continuous mappings

$$\begin{array}{ccc}
(X \square Y) \otimes_p ((X \square Y) \otimes_p (X \square Y)) & \xrightarrow{S} & ((X \square Y) \otimes_p (X \square Y)) \otimes_p (X \square Y) \\
{}_{1_{X \square Y} \otimes \circ} \downarrow & & \downarrow {}_{\circ \otimes 1_{X \square Y}} \\
(X \square Y) \otimes_p (X \square Y) & & (X \square Y) \otimes_p (X \square Y) \\
{}_{\circ} \searrow & & \swarrow {}_{\circ} \\
& X \square Y &
\end{array}$$

We can get restriction of taken bisuperfunctor $X \square Y$ in the categories of $V$-enriched categories with crude functors $Cat_V$, with continuous functors $Cat_V^+$, with cocontinuous functors $Cat_V^-$ or with exact functors $Cat_V^=$.

Let we take an assistant $V$-enriched category $A$. Then we define a superfunctor $F$ for a $V$-enriched category $X$ appointing the $\square$-product $X \square A$ of taken $V$-enriched category $X$ with the assistant $V$-enriched category $A$ and for transports of $V$-enriched categories $f : X \to Y$ appointing the $\square$-product of transports $f \square \mathtt{Id}_A : X \square A \to Y \square A$. We get a superfunctor in the category of $V$-graphs, also in smaller categories of original $V$-graphs or $V$-enriched categories.

We shall try for the new biproduct of $V$-graphs $X \square Y$ to define a coadjoint exponential superfunctor. We take a category of $V'$-graphs $Grph_{V'}$. We have continuous translation of values category $V \Longrightarrow V'$ which change only biproduct with natural arrows of comparison $t_{r,s} : r \otimes s \to r \square s$. Therefore every $V'$-graph $X$ will be a $V$-graph, and every $V'$-graphs transport $f : X \Longrightarrow Y$ will rest a $V$-graphs transport, i. e. we get embedding of categories

$$Grph_{V'} \subset Grph_V .$$

In reality we have coincidence of these two categories.

The choosing of unit arrows $1_x(*) \subset X(x;x)$ is also the same for enriched graphs in both categories.

The difference will appear only when we shall work with contact arrows mappings. Every contact arrows mapping $X(a;b) \square X(b;c) \to X(a;c)$ can be restricted over the tensor product in values category $V$. We get a new mapping with composition

$$X(x;y) \otimes X(y;z) \to X(x;y) \square X(y;z) \to X(y;z) \circ X(y;z) \subset X(x;z) .$$

For a $V'$-graph $X$ with contact arrows mappings we appoint a new $V$-graph $X'$ with restricted contact arrows mappings.

The choice of unit arrow remains neutral for such restricted contact arrows mapping. The neutral property for the source of arrows will be expressed by commuting diagram

$$\begin{array}{ccc}
 & & X(x;y) \\
 & \swarrow & \downarrow \\
* \otimes X(x;y) & = & * \square X(x;y) \\
{\scriptstyle 1_x \otimes 1} \downarrow & & \downarrow {\scriptstyle 1_x \square 1} \\
X(x;x) \otimes X(x;y) & \stackrel{t}{\longrightarrow} & X(x;x) \square X(x;y) \\
 & \searrow & \downarrow \circ \\
 & & X(x;y)
\end{array}$$

The neutral property for the target of arrows will be expressed with the similar commuting diagram.

So we get that every original $V'$-graph remains original $V$-graph with restricted contact arrows mappings i. e. we have faithful imbedding of correspondent category with crude functors

$$OGrph_{V'} \Longrightarrow OGrph_V \ .$$

Any continuous, cocontinuous or natural functor between original $V'$-graphs rests such in the category of original $V$-graphs. We apply the natural property of comparison arrows. There is the diagram for a continuous transport $f : X \Longrightarrow Y$

$$\begin{array}{ccc} X(x;y) \otimes X(y;z) & \to & f(X(x;y)) \otimes f(X(y;z)) \\ t \downarrow & & \downarrow t \\ X(x;y) \square X(y;z) & \to & f(X(x;y)) \square f(X(y;z)) \\ \circ \downarrow & & \downarrow \circ \\ X(x;y) \circ X(y;z) & \to & f(X(x;y)) \circ f(X(y;z)) \\ & & \nearrow \end{array}$$

We have checked faithful embedding of corresponding categories of original graphs with continuous functors $OGrph_{V'}^+ \Longrightarrow OGrph_V^+$, with cocontinuous functors $OGrph_{V'}^- \Longrightarrow OGrph_V^-$, and with natural functors $OGrph_{V'}^{\overline{=}} \Longrightarrow OGrph_V^{\overline{=}}$.

Also we can check that restricted contact arrows mappings maintain the property of associativity.

Let we have commuting diagram for edges set mappings in $V'$-graph $X$

$$\begin{array}{ccc} X(x;y) \square (X(y;z) \square X(z;u)) & \stackrel{s'}{\to} & (X(x;y) \square X(y;z)) \square X(z,;u) \\ 1_{X(x;y)} \square X_{y,z,u} \downarrow & & \downarrow X_{x,y,z} \square 1_{X(z,u)} \\ X(x;y) \square X(y;u) & & X(x;z) \square X(z;u) \\ X_{x,z,u} \searrow & & \swarrow X_{x,z,u} \\ & X(x;u) & \end{array}$$

We must check commuting diagram for deformed restricted diagram with correspondent comparison arrows

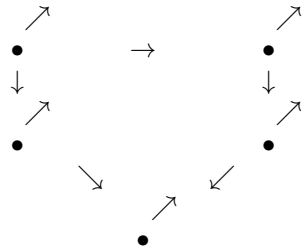

We need only to check commuting quadrats for corespondent comparison arrows. For the upper side we have commuting diagram by an assumption of coherence for associativity transforms by the continuous translation $V \Longrightarrow V''$

$$\begin{array}{ccc} X(x;y) \otimes (X(y;z) \otimes X(z;u)) & \stackrel{s}{\to} & (X(x;y) \otimes X(y;z)) \otimes X(z;u) \\ t \downarrow & & \downarrow t \\ X(x;y) \square (X(y;z) \square X(z;u)) & \stackrel{s'}{\to} & (X(x;y) \square X(y;z)) \square X(z;u) \end{array}$$

Now we shall check commuting diagram for the left side

$$
\begin{array}{ccc}
X(x;y) \otimes (X(y;z) \otimes X(z;u)) & & \\
{\scriptstyle 1\otimes t} \downarrow & \searrow & \\
X(x;y) \otimes (X(y;z) \square X(z;u)) & \xrightarrow{t} & X(x;y) \square (X(y;z) \square X(z;u)) \\
{\scriptstyle 1\otimes \circ} \downarrow & & \downarrow {\scriptstyle 1 \square \circ} \\
X(x;y) \otimes X(y;u) & \xrightarrow{t} & X(x;y) \square X(y;u)
\end{array}
$$

The commuting diagram for the lowest left side is get by definition of restricted contact arrows mapping

$$
\begin{array}{ccc}
X(x;y) \otimes X(y;u) & \xrightarrow{t} & X(x;y) \square X(y;u) \\
& \searrow & \downarrow \circ \\
& & X(x;u)
\end{array}
$$

So we get correspondent faithful embeddings for the different categories of $V$-enriched categories. For the crude functors we get $Cat_{V'} \implies Cat_V$, for the continuous functors we get $Cat_{V'}^+ \implies Cat_V^+$, for the cocontinuous functors we get $Cat_{V'}^- \implies Cat_V^-$ and finally for the natural functor we get $Cat_{V'}^{=} \implies Cat_V^{=}$.

For assistant $V$-graph $A$ we can take the same exponential superfunctor restricted over the image category of faithful embedding of $V'$-graphs

$$Grph_{V'} \implies Grph_V \;,$$

i. e. for $V'$-graph $X$ we take the set of $V$-graphs transports

$$Y_0^A = \{\phi : A \implies Y\} \;,$$

and the edges set between two $V$-graphs transports $\phi, \psi : A \implies Y$ is defined by the product of edges sets

$$Y^A(\phi; \psi) = \prod_{a \in A_0} Y(\phi_0(a); \psi_0(a)) \;.$$

We choose the transported unit arrow

$$1_\phi(*) = \prod_{a \in A_0} 1_{\phi_0(a)}(*) \subset \prod_{a \in A_0} Y(\phi_0(a); \phi_0(a)) \;.$$

This allows to have projections $p_a : Y^A \implies Y$ as $V$-graphs transports maintaining the chosen unit arrows

$$1_\phi^{p_a}(*) = 1_{\phi_0(a)}(*) \subset Y(\phi_0(a); \phi_0(a)) \;.$$

Having the contact arrows mappings

$$\circ : Y(x;y) \square Y(y;z) \to Y(x;z)$$

we take it's restriction and define contact arrows mappings in the exponential $V$-graph

$$Y^A(\phi;\psi) \otimes Y^A(\psi;\xi) = \prod_{a \in A_0} (Y(\phi_0(a);\psi_0(a)) \otimes Y(\psi_0(a);\xi_0(a))) \to$$

$$\prod_{a \in A_0} (Y(\phi_0(a);\psi_0(a)) \,\square\, Y(\psi_0(a);\xi_0(a))) \to$$

$$\prod_{a \in A_0} (Y(\phi_0(a);\psi_0(a)) \circ Y(\psi_0(a);\xi_0(a))) \subset \prod_{a \in A_0} Y(\phi_0(a);\xi_0(a)) = Y^A(\phi;\xi) \ .$$

For an original $V'$-graph $Y$ we get the original exponential $V$-graph again.

For a $V'$-graphs transport $f : X \Longrightarrow Y$ we get the $V$-graphs transport of imbedded $V$-graphs and the changing of target space provides the $V$-graphs transport of exponential $V$-graphs

$$f^A : X^A \Longrightarrow Y^A \ .$$

We get a superfunctor in the category of $V$-graphs $Grph_V$ defined over the partial category of imbedded $V'$-graphs

$$Grph_{V'} \Longrightarrow Grph_V \ .$$

It can also be restricted in the categories of original graphs with natural functors $OGrph_{V'} \Longrightarrow OGrph_V$. with continuous functors $OGrph_{V'}^+ \Longrightarrow OGrph_V^+$, with cocontinuous functors $OGrph_{V'}^- \Longrightarrow OGrph_V^-$, or in the category with natural functors $Grph_{V'}^{=} \Longrightarrow Grph_V^{=}$.

The associativity of contact arrows mappings in $V'$-graph $Y$ provides the associativity of restricted contact arrows mappings in imbedded $V$-graph, and also in exponential $V$-graph $Y^A$. So we get superfunctors in smaller categories of $V$-enriched categories. For the crude functors we get $Cat_{V'} \Longrightarrow Cat_V$, for the continuous functors we get $Cat_{V'}^- \Longrightarrow Cat_V^-$, for the cocontinuous functors we get $Cat_{V'}^+ \Longrightarrow Cat_V^+$, and for the natural functors we get $Cat_{V'}^{=} \Longrightarrow Cat_V^{=}$.

For the $V$-enriched categories we can also define the exponential space of continuous transforms $Y_-^{(A)}$, cocontinuous transforms $Y_+^{(A)}$, or natural transforms $Y_=^{(A)}$. However these spaces defines superfunctors only in the category of natural functors $Cat_V^{=} \Longrightarrow Cat_V^{=}$.

The unit transform for the joint pair of superfunctors in the category of $V$-graphs $Grph_V$ is defined by collection of $V$-graphs transports

$$\lambda_X : X \Longrightarrow (X \,\square\, A)^A \ .$$

They are sections transports defined with chosen unit arrow mappings

$$e_x(*) \subset X(x;x) \ , \quad e_a(*) \subset A(a;a) \ .$$

For a vertex $x \in X_0$ they appoint the section $\Lambda_X(x) : A \Longrightarrow X \,\square\, A$, which for a vertex $a \in A_0$ appoints the couple $\langle x, a \rangle \in X_0 \times A_0$ and have edges set mappings

$$\lambda_X(x)_{a,b} : A(a;b) \to 1_x(*) \,\square\, A(a;b) \subset X(x;x) \,\square\, A(a;b) \ .$$

The edges set mappings for sections transport $\Lambda_X : X \Longrightarrow X \,\square\, A$ is defined by the set of transforms in exponential space $(X \,\square\, A)^A$

$$\lambda_X : X(x;y) \to \prod_{a \in A_0} X(x;y) \,\square\, 1_a(*) \subset \prod_{a \in A_0} X(x;y) \,\square\, A(a;a) \ .$$

The counit transform may be defined by preevaluation $V$-graphs transport defined with precontact arrows mappings

$$\mathbf{ev}_Y^+ : Y^A \,\square\, A \Longrightarrow Y \ .$$

For the former joint pair of superfunctors $X \otimes A$ and $Y^A$ we have preevaluation over tensor product
$$\mathbf{ev}_Y^+ : Y^A \otimes A \Longrightarrow Y \ .$$

It for a $V$-graphs transport $\phi : A \Longrightarrow Y$ and a vertex $a \in A_0$ appoints the value $\phi_0(a) \in Y_0$ and have edges set mappings

$$Y^A(\phi;\psi) \otimes A(a;b) \to Y(\phi_o(a); \psi_0(a)) \otimes Y(\psi_{a,b}(A(a;b))) \to$$

$$Y(\phi_o(a); \psi_0(b)) \circ Y(\phi_{a,b}(A(a;b))) \subset Y(\phi_0(a); \psi_0(b)) \ .$$

However for the new biproduct to have extension of these arrows isn't obvious

$$Y^A(\phi;\psi) \,\square\, A(a;b) \to Y(\phi_o(a); \psi_0(a)) \,\square\, Y(\psi_{a,b}(A(a;b))) \to$$

$$Y(\phi_o(a); \psi_0(b)) \circ Y(\phi_{a,b}(A(a;b))) \subset Y(\phi_0(a); \psi_0(b)) \ .$$

Another counit transform may be defined by postevaluation $V$-graphs transport defined with postcontact arrows mappings

$$\mathbf{ev}_Y^- : Y^A \,\square\, A \Longrightarrow Y \ .$$

It for a $V$-graphs transport $\phi : A \Longrightarrow Y$ and a vertex $a \in A_0$ appoints the value $\phi_0(a) \in Y_0$ and have edges set mappings

$$Y^A(\phi;\psi) \,\square\, A(a;b) \to Y(\phi_{a,b}(A(a;b)) \,\square\, Y(\phi_o(b); \psi_0(b))) \to$$

$$Y(\phi_{a,b}(A(a;b))) \circ Y(\phi_o(b); \psi_0(b)) \subset Y(\phi_0(a); \psi_0(b)) \ .$$

I am sorry that after wasting so many pages of paper I have got so miserable result. We are return to the most general case studying the name appointment in the category of $V$-graphs

$$Grph_V(X \,\square\, A; Y) \to Grph_V(X; Y^A) \ .$$

**Proposition 12.** *Let we have V-graphs X and A with chosen unit arrows*

$$1_x(*) \subset X(x;x) , \quad 1_a(*) \subset A(a;a)$$

*and Y with chosen contact arrows mappings*

$$\circ : Y(x;y) \,\square\, Y(y;z) \to Y(x;z)$$

*allowing preevaluation mappings*

$$Y(\phi_o(a); \psi_0(a)) \,\square\, Y(\psi_{a,b}(A(a;b))) \to$$

$$Y(\phi_o(a); \psi_0(b)) \circ Y(\phi_{a,b}(A(a;b))) \subset Y(\phi_0(a); \psi_0(b)) .$$

*or postevaluation mappings*

$$Y(\phi_{a,b}(A(a;b)) \,\square\, Y(\phi_o(b); \psi_0(b))) \to$$

$$Y(\phi_{a,b}(A(a;b))) \circ Y(\phi_o(b); \psi_0(b)) \subset Y(\phi_0(a); \psi_0(b)) .$$

*For the V-graphs transport* $\Phi : X \,\square\, A \Longrightarrow Y$ *with continuous quadrat of contact arrows mappings*

$$X(x;y) \,\square\, A(a;b) \to \Phi_{\langle x,a\rangle,\langle x,b\rangle}(1_x(*) \,\square\, A(a;b)) \circ \Phi_{\langle x,b\rangle,\langle y,b\rangle}(X(x;y) \,\square\, 1_b(*)) \geq$$

$$\Phi_{\langle x,a\rangle,\langle y,a\rangle}(X(x;y) \otimes 1_a(*)) \circ \Phi_{\langle y,a\rangle,\langle y,b\rangle}(1_y(*) \,\square\, A(a;b)) \subset Y(\Phi_0(x,a); \Phi_0(y;b))$$

*the name* $\Psi : X \Longrightarrow Y^A$ *is V-graphs transport corestricted in the exponential space of continuous transforms* $Y_+^{(A)} \subset Y^A$. *The prerealization is provided with the precontact arrows mapping*

$$\Phi'_{\langle x,a\rangle,\langle y,b\rangle} : X(x;y) \,\square\, A(a;b) =$$

$$\Phi_{\langle x,a\rangle,\langle y,a\rangle}(X(x;y) \otimes 1_a(*)) \circ \Phi_{\langle y,a\rangle,\langle y,b\rangle}(1_y(*) \,\square\, A(a;b)) \subset Y(\Phi_0(x,a); \Phi_0(y;b)) .$$

*The postrealization is provided with the postcontact arrows mappings*

$$\Phi'_{\langle x,a\rangle,\langle y,b\rangle} : X(x;y) \,\square\, A(a;b) =$$

$$\Phi_{\langle x,a\rangle,\langle y,b\rangle}(1_x(*) \,\square\, A(a;b)) \circ \Phi_{\langle x,b\rangle,\langle y,b\rangle}(X(x;y) \,\square\, 1_b(*)) \subset Y(\Phi_0(x,a); \Phi_0(y;b)) .$$

*If the name appointment is corestricted in the exponential space of natural transforms* $Y^{(A)}$, *then we have coinciding both realizations. It coincides with taken V-graphs transport* $\Phi : X \,\square\, A \Longrightarrow Y$ *if we have inequalities of predecomposition and postdecomposition*

$$X(x;y) \,\square\, A(a;b) \to \Phi_{\langle x,a\rangle,\langle y,b\rangle}(1_x(*) \,\square\, A(a;b)) \circ \Phi_{\langle x,b\rangle,\langle y,b\rangle}(X(x;y) \,\square\, 1_b(*)) \geq$$

$$\Phi_{\langle x,a\rangle,\langle y,b\rangle}(X(x;y) \,\square\, A(a;b)) \geq$$

$$\Phi_{\langle x,a\rangle,\langle y,a\rangle}(X(x;y) \otimes 1_a(*)) \circ \Phi_{\langle y,a\rangle,\langle y,b\rangle}(1_y(*) \,\square\, A(a;b)) \subset Y(\Phi_0(x,a); \Phi_0(y;b)) .$$

**Proof:** The proving of such conclusion was provided earlier in various propositions.

□

It can be useful to simplify the notation of contact arrows mappings for $V$-graphs transport $F : X \Longrightarrow Y$

$$F_{x,y} : X(x;y) \to F(X(x;y)) \subset Y(F_0(x); F_0(y)) .$$

In former proposition we can choose $V$-graphs transports $\Phi(x,-) : A \Longrightarrow Y$ with contact arrows mappings

$$A(a;b) \to \Phi_{\langle x,a \rangle, \langle x,b \rangle}(1_x(*) \square A(a;b)) \subset Y(\Phi_0(x,a); \Phi_0(x,b)) ,$$

and $\Phi(-,a) : X \Longrightarrow Y$ with contact arrows mappings

$$X(x;y) \to \Phi_{\langle x,a \rangle, \langle y,a \rangle}(X(x;y) \square 1_a(*)) \subset Y(\Phi_0(x,a); \Phi_0(y,a)) .$$

Then this proposition can be red: If we have coincidence of contact arrows mappings

$$X(x;y) \times A(a;b) \to \Phi(X(x;y),a) \circ \Phi(y,A(a;b)) =$$

$$\Phi(x,A(a;b)) \circ \Phi(X(x;y),b) \subset Y(\Phi_0(x,a); \Phi_0(y;b))$$

then these $V$-graphs transforms can be get as partial ones from some $V$-graphs transport $\Phi : X \times A \Longrightarrow Y$. We can choose such $V$-graphs transport defined by these equal values

$$\Phi(X(x;y), A(a;b)) =$$

$$\Phi(X(x;y),a) \circ \Phi(y,A(a;b)) = \Phi(x,A(a;b)) \circ \Phi(X(x;y),b) .$$

This properties of name appointment remains in smaller categories of $V$-graphs. It may be in the category of original $V$-graphs $OGrph_V$ or the category of $V$-enriched categories with natural functors $Cat_V^{\overline{=}}$.

We remark that having families of functors $\Phi(x,-) : A \Longrightarrow Y$ and $\Phi(-,a) : X \Longrightarrow Y$ with coherent vertexes mappings

$$\Phi(x,-)_0(a) = \Phi(-,a)_0(x)$$

we get two forms with the same vertexes mappings

$$\Phi(x,-)_0(a) = \Phi_0(a,b) = \Phi(-,a)_0(x) .$$

One is with predecomposition property

$$\Phi_{\langle x,a \rangle, \langle y,b \rangle}(X(x;y) \square A(a;b)) = \Phi(X(a;y),a) \circ \Phi(y,A(a;b)) ,$$

another is with postdecomposition property

$$\Phi_{\langle x,a \rangle, \langle y,b \rangle}(X(x;y) \square A(a;b)) = \Phi((x,A(a;b)) \circ \Phi(X(x;y),b) .$$

Both forms has the same name $X \Longrightarrow Y^A$. For vertex $x \in X_0$ it appoints the functor $\Phi(x, -) : A \Longrightarrow Y$ with edges set mappings

$$\Phi_{x,y} : X(x;y) \to \prod_{a \in A_0} \Phi(X(x;y), a)$$

providing the transforms sets between functors $\Phi(x, -)$ and $\Phi(y, -)$.

The situation isn't always so simple. We provide a counterexample with Banach spaces. Arbitrary monoid $X$ can be identified with category $p$ having unique object $* \in p_0$. The edges set coincides with taken monoid

$$X = p(*; *) .$$

The transport of monoids $f : X \to Y$ will be natural functor of correspondent categories $F : p \Longrightarrow r$. The objects mapping is obvious $F_0(*) = *$ and edges set mapping is taken transport

$$f = F_{*,*} : p(*; *) \to r(*; *) .$$

For an assistant monoid $A$ we provide the category $s$. The exponential space $r^s$ will have functors $\phi : s \Longrightarrow r$ with unique objects mapping $\phi_0(*) = *$, only the edges set mapping can be different transports between monoids $f, g : A \to Y$. The transform between such functors will be defined by one member of target monoid $u \in Y$. It is natural exactly when we have commutation

$$f(a) \circ u = u \circ f(a) \uparrow a \in A .$$

We shall take Banach space $X$ as a category with one object. For them we have projective tensor product $X \otimes Y$ which with exponential space of continuous linear mappings $Y^A$ provides truly adjoint pair of functors, i. e. we have bijective name appointment

$$\mathtt{Ban}(X \otimes A; Y) =\to \mathtt{Ban}(X; Y^A) .$$

So for such tensor product we have an evaluation mapping

$$\mathbf{ev}_Y : Y^A \otimes A \to Y .$$

In **G. Valiukevičius 2022** we constructed a relator $X \mathbin{\square} A$ with biequicontinuous tensor product. The contact arrows mapping is defined with extension of sum $+ : X \mathbin{\square} X \to X$ without any difficulty. However extended evaluation mapping
$$\mathbf{ev}_Y : Y^A \mathbin{\square} A \to Y$$

can be defined only in smaller exponential space of integrable linear mappings $f : A \to Y$.

The composition of arrows is the same for both tensor product, but bijective name appointment is get in the smaller category of Banach spaces with integrable linear mappings

$$\texttt{SBan}(X \mathrel{\square} A; Y) =\to \texttt{SBan}(X; Y^A) \ .$$

At the end of this paper I shall fulfill my promise to prove the Yoneda lemma. We wanted to follow the proving from **G. Kelley 1982**. But our way appear more simple. The former complex proving of G. Kelley rests some doubt about its validity. We would wish that a reader would get more beautiful impression from our text.

We remember that a monoidal category $V$ is closed if we have truly joint pair of functors. For chosen an assistant object $s \in V_0$ the adjoint functor is defined by biproduct $s \otimes r$ and for it we demand existence of coadjoint exponential functor $t^s$ with natural bijective name mapping

$$V(s \otimes r; t) =\to V(r; t^s) \ .$$

We have changed the place of assistant object $s$ in biproduct $s \otimes r$. So we shall interchange the prefixes *pre* with *post*. The prefix *pre* is applied when we maintain the order of factors, and the prefix *post* is applied when we change the order of factors to opposite one.

We shall weaken our assumptions. We shall ask that a monoidal category would be *almost closed*. For the joint pair of functors $s \otimes r$ and $t^s$ we ask only existence of unit and counit transforms.

The unit transform is defined by *sections arrows*

$$\lambda_r : r \to (s \otimes r)^s \ ,$$

and the counit transform is defined by *evaluation arrows*

$$\mathbf{ev}_t : s \otimes t^s \to t \ .$$

The sections arrow and the changing of target space by taken form $f : p \otimes s \to p$ provide the *name form*

$$f^\sharp : p \xrightarrow{\lambda_p} (s \otimes p)^s \xrightarrow{f^s} p^s \ .$$

The biproduct of taken name form $g : p \to r^s$ with identity $\texttt{Id}_s : s \to s$ and the evaluation arrow and provide *realization*

$$g^\flat : s \otimes p \xrightarrow{1_s \otimes g} s \otimes r^s \xrightarrow{\mathbf{ev}_r} r \ .$$

In the values category $V$ defined exponential space $t^s \in V_0$ will be called an *inner hom set* $[s; t] = t^s$. It provides the edges sets for an $V$-enriched graph $V^e$

$$V^e(p; r) = [p; r] \ .$$

For the projection to the left factor of biproduct $\sigma : r \otimes * \to r$ we get the *name arrow*
$$i_r : * \xrightarrow{\lambda_*} (r \otimes *)^r \xrightarrow{\sigma^r} r^r .$$

We shall take it as choosing of unit arrow $i_r \in V^e(r; r)$ for an object $r \in V_0$. So we construct an original $V$-graph $V^e$.

Using an associative arrow for biproduct $s \otimes t$
$$s : p \otimes (q \otimes r) \to (p \otimes q) \otimes r$$

together with double evaluation we get *contact arrows mappings* between inner hom sets
$$\cdot : q^p \otimes r^q \to q^p \cdot r^q \subset r^p .$$

First we take the composition of arrows
$$p \otimes (q^p \otimes r^q) \xrightarrow{s} (p \otimes q^p) \otimes r^q \xrightarrow{\mathbf{ev}_q} q \otimes r^q \xrightarrow{\mathbf{ev}_r} r .$$

Then its name arrow will be taken as needed contact arrows mapping.

It is some problem to check under what conditions the chosen unit arrows are neutral for such contact arrows mappings. For the true joint pair of functors $r \otimes s$ and $t^s$ we have bijective name appointment. In our more general case we can demand that chosen special projection $\sigma : r \otimes * \to r$ would be a realization of its name. We could demand this property also for arbitrary *left projection* $\sigma : r \otimes s \to r$. Such property is expressed as commuting diagram

$$\begin{array}{ccccc}
r \otimes s & \xrightarrow{1_r \otimes \lambda_s} & r \otimes (r \otimes s)^r & \xrightarrow{1_r \otimes \sigma^s} & r \otimes r^r \\
& \searrow & \downarrow \mathbf{ev}_{r \otimes s} & & \downarrow \mathbf{ev}_r \\
& & r \otimes s & \xrightarrow{\sigma} & r
\end{array}$$

It will be provided by the first triangular equality

$$\begin{array}{ccc}
r \otimes s & \xrightarrow{1_r \otimes \lambda_s} & r \otimes (r \otimes s)^r \\
& \searrow & \downarrow \mathbf{ev}_{r \otimes s} \\
& & r \otimes s
\end{array}$$

and commuting quadrat

$$\begin{array}{ccc}
r \otimes (r \otimes s)^r & \xrightarrow{1_r \otimes \sigma^r} & r \otimes r^r \\
\mathbf{ev}_{r \otimes s} \downarrow & & \mathbf{ev}_r \\
r \otimes s & \xrightarrow{\sigma} & r
\end{array}$$

First we shall check the neutral property for the arrow's source
$$i_r(*) \cdot t^r = t^r .$$

It is expressed as commuting diagram with right projection $\rho : * \otimes t^r \to t^r$

$$\begin{array}{ccc}
* \otimes t^r & \xrightarrow{i_r \otimes 1_{t^r}} & r^r \otimes t^r \\
\rho \downarrow & \swarrow & \\
t^r & &
\end{array}$$

We shall draw a commuting diagram for the realization of contact arrow mappings over the biproduct $r \otimes (* \otimes t^r) \to (r \otimes *) \otimes t^r$

$$(r \otimes *) \otimes t^r \xrightarrow{(r \otimes \lambda_*) \otimes 1_{t^r}} (r \otimes (r \otimes *)^r) \otimes t^r \xrightarrow{(r \otimes \sigma^r) \otimes 1_{t^r}} (r \otimes r^r) \otimes t^r$$

$$\searrow \qquad \mathbf{ev}_{r \otimes *} \otimes 1_{t^r} \downarrow \qquad\qquad \downarrow \mathbf{ev}_r \otimes 1_{t^r}$$

$$(r \otimes *) \otimes t^r \xrightarrow[\sigma \otimes 1_{t^r}]{} r \otimes t^r$$

which provide needed equality for the name form defining contact arrows mapping for the inner home sets.

In similar way we shall check the neutral property for the arrow's target

$$r^p \cdot i_r(*) = r^p \ .$$

It is expressed as a commuting diagram

$$r^p \otimes * \xrightarrow{1_{r^p} \otimes 1_r} r^p \otimes r^r$$
$$\sigma \downarrow \quad \swarrow$$
$$r^p$$

We shall draw a commuting diagram for the realization of contact arrows mapping over biproduct $p \otimes (r^p \otimes *)$

$$p \otimes (r^p \otimes *) \xrightarrow{1_p \otimes (1_{r^p} \otimes \lambda_*)} p \otimes (r^p \otimes (r^p \otimes *)^{r^p}) \xrightarrow{1_p \otimes (1_{r^p} \otimes \sigma^{r^p})} p \otimes (r^p \otimes (r^p)^{r^p})$$

$$\searrow \qquad 1_p \otimes \mathbf{ev}_{(r^p \otimes *)} \downarrow \qquad\qquad \downarrow 1_p \otimes \mathbf{ev}_{(r^p)}$$

$$p \otimes (r^p \otimes *) \xrightarrow[1_p \otimes \sigma]{} p \otimes r^p$$

which provide needed equality for the name arrow defining contact arrows mapping between inner hom sets.

We shall notice that constructed contact arrows mappings between inner hom sets provides predecomposable form

$$q^p \otimes r^q \to q^p \cdot r^q \subset r^p$$

as we have an equality

$$q^p \cdot r^q = (q^p \cdot i_q) \cdot (1_q \cdot r^q) \ .$$

This property will be changed to another decomposability property, if we take another notation of biproduct. For the biproduct $r^q \otimes q^p$ we get postdecomposable form

$$r^q \otimes q^p \to q^p \cdot r^q \subset r^p \ .$$

For general joint pair of functors $s \otimes r$ and $t^s$ we can have another unit arrows $i_r(*) \subset [r;r]$ and another contact arrows mappings between inner hom sets $q^p \times r^q \to r^p$.

We want to check that the associativity of arrow composition in the values category $V$ provides also the associativity arrow for contact arrows mappings between inner hom sets.

The associativity arrow of biproduct isn't necessary. Later we shall give example of concrete value s category, which isn't banal.

So we probe to define contact arrows mappings between inner hom sets without using such associativity of biproduct $r \otimes s$. We should demand that sections arrows would define natural transform. For arbitrary arrow $f : r \to t$ we used the commuting diagram

$$\begin{array}{ccc} r & \stackrel{\lambda_r}{\longrightarrow} & (s \otimes r)^s \\ f \downarrow & & \downarrow (1_s \otimes f)^s \\ t & \stackrel{}{\underset{\lambda_t}{\longrightarrow}} & (s \otimes t)^s \end{array}$$

In usual cases such property is fulfilled.

Also we demand that evaluation arrow would be natural on its functional argument. For an arrow $g : r \to t$ we get commuting diagram

$$\begin{array}{ccc} s \otimes r^s & \stackrel{\mathbf{ev}_r}{\longrightarrow} & r \\ 1_s \otimes g^s \downarrow & & \downarrow g \\ s \otimes t^s & \underset{\mathbf{ev}_t}{\longrightarrow} & t \end{array}$$

More special property we ask for inner hom sets and biproducts. We need to have isomorphic arrows

$$[p; q] \otimes r = [p; q \otimes r] \ .$$

Then the contact arrows mappings between inner hom sets will be defined as a annulation of middle object

$$A_q : q^p \otimes r^q = (q \otimes r^q)^p \stackrel{(\mathbf{ev}_r)^p}{\longrightarrow} r^p \ .$$

The associativity arrow

$$q^p \cdot (r^q \cdot s^r) \to (q^p \cdot r^q) \cdot s^r$$

will be checked from diagrams with different annulation order.

$$q^p \otimes (r^q \otimes s^r) = (q \otimes (r^q \otimes s^r))^p = (q \otimes (r \otimes s^r)^q)^p \stackrel{A_r}{\longrightarrow} (q \otimes s^q)^p \stackrel{A_q}{\longrightarrow} s^p \ .$$

and

$$(q^p \otimes r^q) \otimes s^r = (q \otimes r^q)^p \otimes s^r \stackrel{A_q}{\longrightarrow} q^p \otimes s^r = (r \otimes s^r)^p \stackrel{A_r}{\longrightarrow} s^p \ .$$

So it rests to check a wanted associativity property by equality or inequality of these two concrete processes.

We shall study a concrete value category produced with recursive schema from the less complex values category. We begin with a category of original graphs with chosen set of original graphs transports

$$V = OGrph \ .$$

The graph $p$ is defined over the set of vertexes $p_0$ with the edges sets $p(u;v) \in Set_0$ for arbitrary couple of vertexes $u, v \in p_0$. An original graph $p$ has chosen unit arrows $1_u \in p(u;u)$. The crude functor $f : p \implies r$ must maintain the unit arrows
$$f_{x,x}(1_x) = 1_{f(x)} \ .$$

In the category of original graphs $OGrph$ we can't expect to have truly joint pair of superfunctors $p \otimes s$ and $r^s$.

Let we have in the target graph $r$ the contact arrows mappings
$$\cdot : r(u;v) \times r(v;w) \to r(u;w) \ .$$

Then for original graphs we can construct the *form* $\Phi : p \otimes s \implies r$ from its *partial forms* $\Phi(p,-) : s \implies r$ and $\Phi(-,s) : p \implies r$ with the help of contact arrows in the target graph $r$. We ask the coherent vertexes mappings. For vertexes $a \in p_0$ and $b \in r_0$ we need the equality for both partial forms
$$\Phi_0(a,b) = \Phi(a,-)_0(b) = \Phi(-,b)_0(a) \ .$$

But for the edges $u \in p(a;a')$ and $v \in s(b;b')$ we have two choices. One is prerealization
$$\Phi^+(u,v) = \Phi(u,b) \circ \Phi(a',v) \ ,$$
another is postrealization
$$\Phi^-(u,v) = \Phi(a,v) \circ \Phi(u,b') \ .$$

The partial forms are defined again from taken form $\Phi$ with the help of chosen unit arrow
$$\Phi(a,-)(v) := \Phi(1_a, u) \ , \quad \Phi(-,b)(u) := \Phi(u, 1_b) \ .$$

The biproduct of two graphs $p \otimes r$ is defined over the Carte tensor product of vertexes sets $p_0 \times r_0$ with the Carte tensor product of edges sets
$$p \otimes r(\langle a,a' \rangle; \langle b,b' \rangle) := p(a;a') \times r(b;b') \ .$$

It will present a Carte product in the category of graphs $Grph$.

Taking a unit arrow
$$1_{\langle a,a' \rangle} = \langle 1_a, 1_{a'} \rangle \in p \otimes r(\langle a,a' \rangle; \langle a,a' \rangle)$$

we get a presentation of Carte product in the category of original graphs $OGrph$ with crude functors.

Additionally taking a product of contact arrows mappings
$$\langle u, u' \rangle \cdot \langle v, v' \rangle = \langle u \cdot v, u' \cdot v' \rangle$$

such biproduct will become a presentation of Carte product in the smaller category of original graphs $OGrph^=$ with natural functors.

Such biproduct defines superfunctor in the category of original graphs. It is true not only for the category of crude functors $OGrph$, but also for smaller categories of continuous functors $OGrph^-$, cocontinuous functors $OGrph^+$, or natural functors $OGrph^=$.

For the functors $f : p \Longrightarrow p'$ and $g : r \Longrightarrow r'$ we take a product functor $f \otimes g : p \otimes r \Longrightarrow p' \otimes r'$. We have the product of vertexes mappings

$$(f \otimes g)_0 = f_0 \times g_0 \;,$$

and for edges set mappings $f \otimes g_{\langle x,y\rangle,\langle x',y'\rangle}$ we take a Cauchy tensor product

$$f_{x,x'} \times g_{y,y'} : p(x;x') \times r(y;y') \to p(f_0(x); f_0(x')) \times r(g_0(y); g_0(y')) \;.$$

If we have continuous quadrats for the graphs transports $f : X \Longrightarrow X'$ and $g : Y \Longrightarrow Y'$

$$\begin{array}{ccc} X \otimes_p X & \xrightarrow{f \otimes f} & X' \otimes_p X' \\ \circ \downarrow & & \downarrow \circ \\ X & \xrightarrow{f} & X' \end{array} \qquad \begin{array}{ccc} Y \otimes_p Y & \xrightarrow{g \otimes g} & Y' \otimes_p Y' \\ \circ \downarrow & & \downarrow \circ \\ Y & \xrightarrow{g} & Y' \end{array}$$

then we get again a continuous quadrats for the product of graphs transports $f \otimes g : X \otimes Y \Longrightarrow X' \otimes Y'$

$$\begin{array}{ccc} (X \otimes Y) \otimes_p (X \otimes Y) & \xrightarrow{(f \otimes g) \otimes_p (f \otimes g)} & (X' \otimes g') \otimes_p (X' \otimes Y') \\ \circ \downarrow & & \downarrow \circ \\ X \otimes Y & \xrightarrow{f \otimes g} & X' \otimes Y' \end{array}$$

We could take some relator as a larger biproduct $p \otimes r \subset p \square r$, but it would be hard task to check when we can get needed biproduct of taken functors $f : p \Longrightarrow p'$ and $g : r \Longrightarrow r'$. This can be imagine as a problem to have some integral of taken functors.

Necessarily we need to have a biproduct

$$f \square g : * \square * \to X \square Y$$

of arbitrary points $f : * \Longrightarrow X$ and $g : * \Longrightarrow Y$.

Also we want to have biproduct for arbitrary functor $f : p \Longrightarrow p'$ with identity functor of assistant graph $\mathtt{Id}_s : s \Longrightarrow s$

$$f \square \mathtt{Id}_s : p \square s \Longrightarrow p' \square s \;.$$

The exponential space $p^s$ is constructed as a graph with vertexes set compounded of all functors

$$(p^s)_0 = \{\phi : s \Longrightarrow p\} \;.$$

The edges set will be compounded by all transforms

$$p^s(\phi; \psi) = \{\alpha : \phi \to \psi\} \;.$$

Such transforms will be defined by collection of edges

$$\alpha = \langle \alpha_a \uparrow a \in s_0 \rangle \in \prod_{a \in s_0} p(\phi_0(a); \psi_0(a)) \ .$$

A unit arrow will be defined by collection of unit arrows

$$1_\phi = \langle 1_{\phi_0(a)} \uparrow a \in A \rangle \in \prod_{a \in A_0} p(\phi_0(a); \phi_0(a)) \ .$$

A contact arrows mapping

$$\cdot : p^s(\phi; \psi) \times p^s(\psi; \xi) \to p^s(\phi; \xi)$$

is also defined as a product of contact arrows mappings in the target space. For transforms $\alpha \in \prod_{a \in A_0} p(\phi_0(a); \psi_0(a))$ and $\beta \in \prod_{a \in A_0} p(\psi_0(a); \xi_0(a))$ it appoints

$$\alpha \cdot \beta = \prod_{a \in s_0} \alpha_a \cdot \beta_a \in \prod_{a \in s_0} p(\phi_0(a); \xi_0(a)) \ .$$

Such exponential space $r^s$ also defines a superfunctor within the category $OGrph$. For a functor between target spaces $f : p \Longrightarrow r$ we appoint the changing of target space $f^s : p^s \Longrightarrow r^s$. Such appointment maintains the identity functor

$$(\mathtt{Id}_p)^s = \mathtt{Id}_{p^s} \ ,$$

also it maintains composition of functors

$$(f \circ g)^s = f^s \circ g^s \ ,$$

as we have associative the composition in the category $OGrph$.

We would want to work also with a smaller exponential space $p^{(s)}$ which has smaller edges sets $p^{(s)}(\phi; \psi)$ compounded by natural transforms $\alpha : \phi \to \psi$, i. e. having commuting diagrams for arbitrary edge $u \in s(a; b)$

$$\begin{array}{ccc} \phi_0(a) & \stackrel{\alpha_a}{\longrightarrow} & \psi_0(a) \\ \phi(u) \downarrow & & \downarrow \psi(u) \\ \phi_0(b) & \stackrel{}{\underset{\alpha_b}{\longrightarrow}} & \psi_0(b) \end{array}$$

However without the property of associativity for contact arrows mappings we can't guaranty that contact arrows mappings rest in the set of natural transforms. Such exponential space would define superfunctor only for natural functors. If functor $f : p \Longrightarrow r$ maintains the contact arrows mappings

$$f(u \cdot v) = f(u) \cdot f(v) \uparrow u \in p(x; y), v \in p(y; z) \ ,$$

then transported transform $\alpha \star f$ rests natural transform as we get commuting diagram again for arbitrary edge $u \in s(a; b)$

$$\begin{array}{ccc} f_0(\phi_0(a)) & \stackrel{f(\alpha_a)}{\longrightarrow} & f_0(\psi_0(a)) \\ f(\phi(u)) \downarrow & \searrow & \downarrow f(\psi(u)) \\ f_0(\phi_0(b)) & \stackrel{}{\underset{f(\alpha_b)}{\longrightarrow}} & f_0(\psi_0(b)) \end{array}$$

For the joint pair of superfunctors $s \square p$ and $r^s$ we shall define unit transform by collection of sections transports $\lambda_p : p \Longrightarrow (s \square p)^s$. The sections $\lambda_p(x) : s \Longrightarrow s \square p$ are defined by biproduct of identity transport $\text{Id}_s : s \Longrightarrow s$ with chosen unit arrow $1_x \in p(x;x)$, i. e. we have vertexes mapping

$$a \to \langle a, x \rangle \in p_0 \times s_0$$

and edges set mappings

$$\lambda_p(x)_{a,b} : s(a;b) \to s(a;b) \square 1_x \subset s(a;b) \square p(x;x) \ .$$

The counit transform is defined by collection of evaluation transports

$$\mathbf{ev}_r : s \square r^s \Longrightarrow r \ .$$

It for a functor $\phi : s \Longrightarrow r$ and vertex $a \in s_0$ appoints the value $\phi_0(a) \in r_0$ and have edges set mappings

$$s(a;b) \square r^s(\phi;\psi) \to r(\phi_0(a);\psi_0(b)) \ .$$

It can be done in two ways. A postevaluation map is provided for an edge $u \in s(a;b)$ and a transform $\alpha : \phi \to \psi$ by extension over the larger new biproduct of postcontact arrows mapping

$$s(a;b) \square r^s(\phi;\psi) \to .\alpha_a \circ \psi_{a,b}(u) \in r(\phi_0(a);\psi_0(b)) \ .$$

May be we shall take a smaller exponential space $(r^s)'$ to achieve this task.

A preevaluation map for an edge $u \in s(a;b)$ and a transform $\alpha : \phi \to \psi$ is provided by extension over the larger new biproduct of precontact arrows mapping

$$s(a;b) \square r^s(\phi;\psi) \to \phi_{a,b}(u) \circ \alpha_b \in p(\phi_0(a);\psi_0(b)) \ .$$

Also may be we shall take a smaller exponential space $(r^s)'$ to achieve this task.

We get a bijective name appointment

$$OGrph(s \square p; r) =\to OGrph(p; r^s)$$

for the forms $\Phi : s \square p \Longrightarrow p$ with a predecomposition property

$$\Phi(s(a;b); p(x;y)) = \Phi(1_x, s(a;b)) \cdot \Phi(p(x;y), 1_b) \ .$$

A postevaluation map for an edge $u \in s(a;b)$ and a transform $\alpha : \phi \to \psi$ is provided by extension over the larger new biproduct of postcontact arrows mapping

$$s(a;b) \square r^s(\phi;\psi) \to \alpha_a \circ \psi_{a,b}(u) \in p(\phi_0(a);\psi_0(b)) \ .$$

Also may be we shall take a smaller exponential space $(r^s)'$ to achieve this task.

We get a bijective name appointment

$$OGrph(s \square p; r) =\to OGrph(p; r^s)$$

for the forms $\Phi : s \square p \Longrightarrow p$ with a postdecomposition property

$$\Phi(s(a;b); p(x;y)) = \Phi(p(x;y), 1_a) \cdot \Phi(1_y; s(a;b)) ,$$

It is hard task to construct a new biproduct by some topological extension. For example in the category of metric spaces $Cat_{[0,\infty]^+}$ we can apply a metric of uniform convergence in the space of nonexpanding mappings $Cat_{[0,\infty]^+}(X;Y)$, but such metric space is already complete. Therefore we haven't got any interesting completion.

There is a more easy task to take some new biproduct in the values category $V = OGrph$ which would be relator to the former biproduct

$$p \otimes r \to p \square r$$

and take as target space only such original graphs $t \in V_0$ which would have contact arrows mappings extended over the new biproduct

$$\circ : t(x;y) \square t(y;z) \to t(x;z) .$$

For the partial forms as functors $\Phi(p,-) : r \Longrightarrow t$ and $\Phi(-,r) : p \Longrightarrow t$ we get again the predecomposable form $\Phi^+ : p \square r \Longrightarrow t$ with equal appointments

$$\Phi^+(p(x;x') \square r(y;y')) := \Phi(p(x;x'), y) \square \Phi(x'; r(y;y'))$$

and postdecomposable form $\Phi^- : p \square r \Longrightarrow t$ with another appointment

$$\Phi^-(p(x;x'), r(y;y')) := \Phi(x; r(y;y')) \square \Phi(p(x;x'), y') .$$

Later we prefer to talk about $V$-enriched graphs. The ordinary graph is also $Set$-enriched graph. A $V$-graph $X$ is defined by the set of vertexes $X_0 \in Set_0$ and edges sets $X(a;b) \in V_0$. The $V$-graphs transport $f : X \Longrightarrow Y$ is defined by the mapping of vertexes $f_0 : X_0 \to Y_0$ as an arrow in the category of sets $Set$ and edges set mappings

$$f_{a,b} : X(a;b) \to Y(f_0(a); f_0(b))$$

as an arrows in the values category $V$. All $V$-graphs transports compound the category $Grph_V$.

Now we can take the category of original $V$-graphs as a new values category. At first we can take the category compounded by all crude functors $W = OGrph_V$. Such functors will be $V$-graphs transports $f : X \Longrightarrow Y$ with edges $V$-object mappings

$$f_{x,y} : X(x;y) \to Y(f_0(x); f_0(y))$$

maintaining unit arrow choosing

$$f_{x,x}(1_x(*)) = 1_{f_0(x)}(*) .$$

We shall apply the words *V-objects* and *V-arrows* for the object and arrows in the category $V$.

For two parallel arrows we can define inequality by existing transform between correspondent functors $f, g : X \Longrightarrow Y$. We have the $V$-object of all transforms equal for the Carte product of edges $V$-objects in the target $V$-graph $Y$

$$Y^X(f;g) = \prod_{a \in X_0} Y(f_0(a); g_0(a))$$

and we shall say $f \geq g$ if exists some point in this product as an arrow from former values category $V$

$$(* \to \prod_{a \in X_0} Y(f_0(a)'g_0(a))) = \prod_{a \in X_0} (* \to Y_0(f_0(a); g_0(a))) \ .$$

Also we can take a stronger order relation provided by the partial $V$-object of continuous, cocontinuous or natural transforms, But for it we shall need contact arrows $V$-mappings in the target space $Y$. An inequality with existing special transform will provide the inequality defined with existing of general transforms.

As usually we take the Carte tensor product of vertexes sets $X_0 \times Y_0$ with biproduct of edges $V$-objects from the former values category $V$. For the monotone biproduct in $p \otimes r$ values category $V$ we get again monotone biproduct of original $V$-graphs $X \otimes Y$.

For any assistant space $T$ we get a superfunctor in the category $OGrph_V$ defined by biproduct with assistant space $T \otimes X$. For a functor $f : X \Longrightarrow Y$ it appoints the biproduct

$$\mathtt{Id}_T \otimes f : T \otimes T \Longrightarrow T \otimes Y \ .$$

.

For cocontinuous relator's involution

$$m_{\langle p,r \rangle, \langle s,t \rangle} : (p \otimes r) \otimes (s \otimes t) \to (p \otimes s) \otimes (r \otimes t)$$

we get such superfunctor restricted in the category $OGrph_V^-$ of continuous transports between original $V$-graphs. Also for continuous relator's involution we get restriction in the category $OGrph_V^-$ of cocontinuous transports between original $V$-graphs,

It is possible to choose a larger biproduct

$$X \otimes Y \Longrightarrow X \square Y$$

with relator's transport to the former biproduct

$$N_{\langle X,Y \rangle, \langle Z,W \rangle} : (X \square Y) \otimes (Z \square W) \Longrightarrow (X \otimes Z) \square (Y \otimes W) \ .$$

We can choose such new biproduct being monotone superfunctor

$$OGrph_V \otimes OGrph_V \Longrightarrow OGrph_V \ ,$$

or having continuous or cocontinuous relator's involution

$$M_{\langle X,Y \rangle, \langle Z,W \rangle} : (X \square Y) \square (Z \square W) \Longrightarrow (X \square Z) \square (Y \square W) \ .$$

Also we are choosing the *exponential space* $Y^T$ defined with the partial set of original $V$-graphs transports as the vertexes set

$$(Y^T)_0 \subset \{f : T \Longrightarrow Y\}$$

and the edges $V$-objects taken as as a partial object from all transforms

$$Y^T(\phi; \psi) \subset \prod_{a \in T_0} Y(\phi_0(a); \psi_0(a)) \ .$$

For a target space $Y$ with contact arrows $V$-mappings

$$\circ : Y(x;y) \otimes Y(y;z) \to Y(x;y) \circ Y(y;z) \subset Y(x;z)$$

we get an *exponential superfunctor* provided by the target changing with an original $V$-graphs transport $f : X \Longrightarrow Y$. For a $V$-graphs transport $\phi : T \Longrightarrow X$ we get the new $V$-graphs transport $\phi \circ f : T \Longrightarrow Y$. For continuous taken transport we get continuous changing of target space. And for cocontinuous one we get cocontinuous changing of target space. It rests to check the maintenance of edges $V$-object.

Such superfunctor will be monotone one. For the inequality between parallel $V$-graphs transports $f, f' : X \Longrightarrow Y$, i. e. we have equal first projection

$$f_0(x) = f'_0(x) \uparrow x \in X_0$$

and inequalities for edges $V$-objects mappings

$$f_{x,y} \geq f'_{x,y} \uparrow x, y \in X_0 \ ,$$

we get $V$-graphs transports

$$f^T, (f')^T : X^T \Longrightarrow Y^T$$

with equality of mappings for the equal secondary projection

$$(f_0)^{T_0} = (f'_0)^{T_0} : (X_0)^{T_0} \to (Y_0)^{T_0} \ .$$

We get the inequalities for the $V$-graphs transports of first projection

$$\phi \cdot f \geq \phi \cdot f' \uparrow \phi \in X^T \ ,$$

i. e. we have for edges $V$-objects mappings we have inequalities

$$T(s;t) \to f_{\phi_0(s), f'_0(t)}(\phi_{s,t}(T(s;t)) \geq$$

$$f'_{\phi_0(s),f'_0(t)}(\phi_{s,t}(T(s;t)) \subset Y(f'_0(\phi_0(s)); f'_0(\phi_0(t)))\ .$$

Also we get inequalities for edges $V$-objects mappings

$$X(a;b) \to f_{a,b}(X(a;b)) \geq$$

$$f'_{a,b}(X(a;b)) \subset Y^T(f'_0(a); f'_0(b)) = \prod_{s \in T_0} Y(f'_0(a)_s; f'_0(b)_s)\ .$$

Also the sections transport becomes monotone

$$\lambda_X : X \Longrightarrow (T \otimes X)^T\ .$$

Additionally we demand that changing of target with contact arrows $V$-mappings would be monotone for arbitrary $V$-object $r \in V_0$

$$V(r; Y(x;y)) \to V(r; Y(x;y) \circ Y(y;z)) \subset Y(r; Y(x;z)) \uparrow r \in V_0\ ,$$

Such assumption provides the monotone preevaluation transports

$$\mathbf{ev}_Y^+ : T \otimes Y^T \Longrightarrow Y$$

defined with precontact arrows mappings

$$T(s;t) \otimes Y^T(\phi; \psi) \to \phi_{s,t}(T(s;t)) \circ Y(\phi_0(t); \psi_0(t)) \subset Y(\phi_0(s); \psi_0(t))\ ,$$

For the monotone target's changing

$$V(r; Y(y;z)) \to V(r; Y(x;y) \circ Y(y;z)) \subset Y(r; Y(x;z)) \uparrow r \in V_0\ .$$

we get monotone

$$\mathbf{ev}_Y^- : T \otimes Y^T \Longrightarrow Y$$

defined with postcontact arrows mappings

$$T(s;t) \otimes Y^T(\phi; \psi) \to Y(\phi_0(s); \psi_0(s)) \circ \psi_{s,t}(T(s;t)) \subset Y(\phi_0(s); \psi_0(t))\ .$$

So the category $OGrph_V$ can be consider as a new values category.

We wish to convert the category of $V$-graphs $Grph_V$ into an enriched graph by itself. . We have constructed the exponential space $Y^X$ in the category of $V$-graphs $Grph_V(X;Y)$ taking all $V$-graphs transports $f : X \Longrightarrow Y$ as vertexes and the sets of transforms between them

$$Y^X(\phi; \psi) = \prod_{a \in X_0} Y(\phi_0(a); \psi_0(a))$$

as edges sets. Usually we denote this exponential space

$$Y^X = [X; Y]\ .$$

We have composition of arrows in this category

$$Grph_V(X;Y) \times Grph_V(Y;Z) \to Grph_V(X;Z) .$$

Now the set of transforms between two $V$-graphs transports $\phi, \psi : X \Longrightarrow Y$ will provide an edges $V$-object $Y^X(\phi; \psi)$ in the exponential space $Y^X$. We need to define contact arrows mappings for such exponential spaces. They will be edges $V$-object mappings for the $V$-graphs transport over pullbacked biproduct

$$\circ : Y^X \otimes_p Z^Y \Longrightarrow Z^X .$$

We shall define the contact arrows mappings $Y^X \otimes Z^Y \Longrightarrow Z^Y$ similarly to the double evaluation.

At first we shall construct the $V$-enriched continuous quadrats as $V$-enriched graph and produce the horizontal contact mappings between them in arbitrary category of original $V$-graphs. In the category of $V$-enriched graphs $Grph_V$ a continuous quadrat $Q$ is a diagram

$$\begin{array}{ccc} \bigcirc & \xrightarrow{M} & \cdot \\ I_1 \downarrow & & \downarrow I_2 \\ \cdot & \xrightarrow[N]{} & \bullet \nearrow \end{array}$$

A *source of diagram* is a $V$-graph $\bigcirc$ and a *target of diagram* is a $V$-graph $\bullet$. A *circuit in diagram* is a chain of composable functors from the source $V$-graph to the target $V$-graph. For taken quadrat $Q$ we have only two circuits. One is $I_1 N : \bigcirc \Longrightarrow \bullet$ and another is $M I_2 : \bigcirc \Longrightarrow \bullet$. The continuous quadrat $Q$ is defined by the partial set of transforms defined in the product of edges sets from the final $V$-graph

$$Q(I_1 N; M I_2) \subset [\bigcirc; \bullet](I_1 N; M I_2) = \prod_{a \in \bigcirc_0} \bullet((I_1 N)_0(a); (M I_2)_0(a)) .$$

For a general diagram we construct a new original graph with all circuits as vertexes, and the sets of transforms between two circuits can be defined by products of edges sets in diagram's target space.

We shall try to define horizontal composition of two continuous quadrats $Q_1 \bullet Q_2$ as decomposable form

$$\begin{array}{ccccccc} \bigcirc & \xrightarrow{M_1} & \cdot & & \xrightarrow{M_2} & \cdot \\ I_1 \downarrow & & \downarrow I_2 & & \downarrow I_3 \\ \cdot & \xrightarrow[N_1]{} & \cdot \nearrow & \xrightarrow[N_2]{} & \bigcirc & \nearrow \end{array}$$

First we shall define partial forms. We shall demand that arrows in value category $V$ provides the images of partial objects, i. e. such category is a *regular* by **I. Freyd, A. Scedrov 1990** §1.5.

The new graphs transport

$$\alpha * 1_{N_2} := \alpha \star N_2$$

will be defined by transport of diagrams target with the functor $N_2$. The edges set mapping of functor $N_2$ will provide the image of partial set

$$N_2(Q(I_1N_1; M_1I_2)) \subset [\circ; \bullet](I_1N_1N_2; M_1N_2N_2) .$$

Another partial transport

$$1_{M_1} * \beta := M_1 \star \beta$$

will be defined by transport of diagram source by functor $M_1$. In the case of $Set$-enriched graphs we take the same arrows of transform $\beta$, only indexes are vertexes from the new source object. For more general values category $V$ we can take the projection of product

$$\prod_x \bullet(((I_2N_2)_0(x); (M_2I_3)_0(x)) \to \prod_{x=(M_1)_0(a)} \bullet((M_1I_2N_2)_0(a); (M_1M_2I_3)_0(a))$$

and the image of partial object $O_2 \subset \prod_x \bullet((I_2M_1)_0(x); (M_2I_3)_0(x))$.

If the target object of new diagram has contact arrows mappings, then we get predecomposable appointment

$$\alpha *^+ \beta := (\alpha \star N_2) \circ (M_1 \star \beta) .$$

It will define the *prehorizontal composition* of continuous quadrats.

Also we can get postdecomposable appointment for cocontinuous quadrats

$$\alpha *^- \beta = (M_1 \star \beta) \circ (\alpha \star N_1) .$$

If the target object has associative contact arrows mappings, then we get again the associative composition of continuous or cocontinuous quadrats.

The deformed quadrats with identic functors $I_1 = \text{Id}_X$, $I_2 = \text{Id}_Y$ will be identified with transforms between functors $M, N : X \Longrightarrow Y$.

The *prehorizontal composition* $\alpha *^+ \beta : \phi \circ \xi \to \psi \circ \kappa$ of two transforms $\alpha : \phi \to \psi$ between functors $\phi, \psi : X \Longrightarrow Y$ and $\beta : \xi \to \kappa$ between functors $\xi, \kappa : Y \Longrightarrow Z$ are calculated as *horizontal composition of two continuous quadrats*

$$\begin{array}{ccccc} X & \xrightarrow{\phi} & Y & \xrightarrow{\xi} & Z \\ \uparrow & & \uparrow & & \uparrow \\ X & \xrightarrow[\psi]{} & Y & \xrightarrow[\kappa]{} & Z \end{array}$$

with $\alpha$ and $\beta$ arrows.

$$\alpha *^+ \beta = (\alpha \star \xi) \circ (\psi \star \beta) .$$

The edges set mapping of preevaluation functor $\mathbf{ev}_Y : X \otimes Y^X \Longrightarrow Y$ is expressed as horizontal composition between continuous quadrats. For an edges set $X(a; b)$ we get a set of precontact arrows

$$\alpha_a \circ \psi(w) : \phi(a) \to \psi(b)$$

$$
\begin{array}{ccccc}
 & & \searrow^{w} & & \searrow^{\alpha} \\
* \xrightarrow{a} & X & \xrightarrow{\psi} & Y & \\
 & \uparrow & & \uparrow & \\
* \xrightarrow{b} & X & \xrightarrow{\phi} & Y &
\end{array}
$$

In the case of associative contact arrows mappings from final $V$-graph $Z$ the *double preevaluation* provides the horizontal composition of three continuous quadrats

$$
\begin{array}{ccccccc}
 & & \searrow^{w} & & \searrow^{\alpha} & & \searrow^{\beta} \\
* \xrightarrow{a} & X & \xrightarrow{\phi} & Y & \xrightarrow{\xi} & Z & \\
 & \uparrow & & \uparrow & & \uparrow & \\
* \xrightarrow{b} & X & \xrightarrow{\psi} & Y & \xrightarrow{\kappa} & Z &
\end{array}
$$

and coincides with the preevaluation of earlier defined prehorizontal composition.

This can be produced also for the horizontal composition of cocontinuous quadrats.

$$
\begin{array}{ccccc}
X & \xrightarrow{\phi} & Y & \xrightarrow{\xi} & Z \\
\uparrow & & \uparrow & & \uparrow \\
X & \xrightarrow{\psi} & Y & \xrightarrow{\kappa} & Z \\
 & & \swarrow_{\alpha} & & \swarrow_{\beta}
\end{array}
$$

$$\alpha *^{-} \beta = (\phi * \beta) \circ (\alpha \star \kappa) \ .$$

For G. Kelley it was important to explain his notation by constructing a functor $\text{Hom}_X : X^{op} \otimes X \Longrightarrow V$ which would have the appointment of vertexes

$$(\text{Hom}_X)_0(a;b) = X(a;b) \in V_0$$

and the edges sets mappings

$$\text{Hom}_X : X^{op}(a;b) \times X(c;d) \to [X(a;c); X(b;d)]$$

extending target's changing mappings $X(1;\beta) : X(a;b) \to X(a;c)$ and source' changing mappings $X(\alpha;1) : X(b;d) \to X(a;d)$. For us it will be an interesting exercise.

For the form's $\Phi : X_1 \otimes X_2 \Longrightarrow Y$ name $X_1 \Longrightarrow [X_2;Y]$ we shall apply more special notation. We take the exponential space $[X_1 \otimes X_2; Y]$ and construct partial form over the first argument applying identic mapping $X_1 \Longrightarrow X_1$ composed with taken exponential space

$$X_1 \Longrightarrow [\bigcirc \otimes X_2; Y] \ .$$

Also we can construct another name with another composition over the second argument

$$X_2 \Longrightarrow [X_1 \otimes \bigcirc; Y] \ .$$

So we can express different partial forms more graphically.

We shall use a diagram to explain the notation in our problem of functor $Hom_X$

$$\begin{array}{ccc} b & \to & d \\ \alpha \downarrow & & \uparrow \beta \\ a & \to & c \end{array}$$

Such functor can be get in two ways. One is defined as an predecomposable form, another is defined as an postdecomposable form. We shall apply the names of predecomposable form defined by contact arrows mappings in taken original $V$-graph $X$.

$$X(a,b) \otimes X(b;c) \to X(a;b) \circ X(b;c) \subset X(a;c) \ .$$

For it we can check equality

$$X(a;b) \circ X(b;c) = (X(a;b) \circ 1_b(*)) \circ (1_b(*) \circ A(b;c)) \ .$$

One name will be

$$X(b;c) \to [X(a;c) \circ \bigcirc; X(a;c)] = [X(a;b); X(a;c)] \ .$$

Another *dual name* will be

$$X(b;c) \to [\bigcirc \circ X(b;c); X(a;c)] = [X(b;c); X(a;c)] \ .$$

With such names we define partial forms for desired predecomposable form

$$\text{Hom}_X^+(X^{op}(a;b) \times X(c;d)) = \text{Hom}_X(X^{op}(a;b) \otimes 1_c(*)) \cdot \text{Hom}_X(1_b(*) \times X(c;d))$$

One partial form will be

$$\text{Hom}_X(X^{op}(a;b) \times 1_c(*)) = [\bigcirc \circ X(a;c); X(b;c)] \subset [X(a;c); X(b;c)]$$

$$\begin{array}{cc} b & \\ \downarrow & \searrow \\ a & \to c \end{array}$$

and another partial form

$$\text{Hom}_X(1_b(*) \times X(c;d)) = [X(b;c) \circ \bigcirc; X(b;d)] \subset [X(b;c); X(c;d)]$$

$$\begin{array}{ccc} b & \to & d \\ & \searrow & \uparrow \\ & & c \end{array}$$

The defined predecomposable form can be explain by diagram

$$\begin{array}{ccc} b & \to & d \\ \downarrow & \searrow & \uparrow \\ a & \to & c \end{array}$$

We shall calculate such mappings in more details for the prehorizontal composition between exponential spaces in our values category $V = OGrph$. Later such category will be called a *pointed* one.

Let we have original graphs $p = X^{op}(a;b)$, $q = X(a;c)$, $r = X(c;d)$. We denote the vertexes $u \in p_0$, $v \in q_0$, $w \in r$ and edges $\alpha \in p(u, u')$, $\beta \in q(v; v')$, $\gamma \in r(w; w')$.

The functor between original graphs

$$X^{op}(a;b) \times X(c;d) \Longrightarrow [\bigcirc \circ X(a;c); X(b;c)] \cdot [\bigcirc \circ X(c;d); X(b;d)]$$

is composed as functor to the *first factor*

$$X^{op}(a;b) \Longrightarrow [\bigcirc \circ X(a;c); X(a;d)] \ ,$$

then the functor over two arguments to the *second factor*

$$X^{op}(a;b) \times [\bigcirc \circ X(a;c); X(a;d)] \times X(c;d) \Longrightarrow [\bigcirc \circ X(c;d); X(b;d)] \ ,$$

and finally we calculate the form over such two mappings

$$[\bigcirc \circ X(a;c); X(b;c)] \cdot [\bigcirc \circ X(c;d); X(b;d)] \ .$$

We shall calculate the realization of such composition. For the first functor we get a realization defined with contact arrows mappings in $V$-graph $X$

$$X^{op}(a;b) \times X(a;c) \Longrightarrow X(b;a) \circ X(a;c) \subset X(b;c) \ .$$

For the second factor we get a realization

$$X^{op}(a;b) \times X(b;c) \times X(c;d) \Longrightarrow (X(b;a) \circ X(a;c)) \circ X(c;d) \subset X(b;d) \ .$$

Finally we need to calculate the prehorizontal composition of transforms between provided original graphs transports.

$$\begin{array}{ccccccc} & & & \searrow^{A(\alpha)} & & & \searrow^{A(\gamma)} \\ X^{op}(a;b) & \xrightarrow{T_u} & X(b;c) & & \xrightarrow{T_w} & X(b;d) & \\ \uparrow & & \uparrow & & & \uparrow & \\ X^{op}(a;b) & \xrightarrow[T_{v'}]{} & X(b;c) & & \xrightarrow[T_{w'}]{} & X(b;d) & \end{array}$$

$$A(\alpha) *^+ A(\gamma) = (A(\alpha) \star T_w) \cdot (T_{v'} \star A(\gamma)) \ .$$

For a vertex $u \in p_0$ we get a functor $T_u : X(a;c) \Longrightarrow X(b;c)$ in the first factor defined with contact arrow mappings from taken $V$-graph $X$. For a vertex $v \in q_0$ it appoints $u \circ v \in X(b;c)$ and for an edge $\beta : v \to v'$ it appoints the new edge $1_u \circ \beta : u \circ v \to u \circ v'$.

For an edge $\alpha : u \to u'$ it appoints the transform between former functors $A(\alpha) : T_u \to T_{u'}$ defined with collection of edges $\alpha \circ 1_v : u \circ v \to u' \circ v$

For a vertex $w \in r_0$ we get a functor $T_w : X(b;c) \Longrightarrow X(b;d)$ in the second factor defined with contact arrows mappings from taken $V$-graph $X$. For a vertex $u \circ v \in X(b;c)_0$ it appoints the new vertex $(u \circ v) \circ w$ and for an edge $(\alpha \circ \beta) : u \circ v \to u' \circ v'$ it appoints the new edge $(\alpha \circ \beta) \circ 1_w : (u \circ v) \circ w \to (u' \circ v') \circ w$. For an edge $\gamma : w \to w'$ it appoints the transform between such functors $A(\gamma) : T_w \to T_{w'}$ defined with collection of edges

$$(1_u \circ 1_v) \circ \gamma : (u \circ v) \circ w \to (u \circ v) \circ w; \ .$$

The prehorizontal composition of such transforms

$$A(\alpha) = \alpha \circ 1_v : u \circ v \to u' \circ v \ ,$$

$$A(\gamma) = (1_u \circ 1_v) \circ \gamma : (u \circ v) \circ w \to (u \circ v) \circ w'$$

will be

$$(\alpha \circ 1_v) \circ 1_w : (u \circ v) \circ w \to (u' \circ v) \circ w$$

composed with

$$(1_{u'} \circ 1_v) \circ \gamma : (u' \circ v) \circ w \to (u' \circ v) \circ w' \ ,$$

i.e. we get transform defined with collection of arrows

$$(\alpha \circ 1_v) \circ \gamma : (u \circ v) \circ w \to (u' \circ v) \circ w'$$

between functors defined by the vertexes $u \in p_0$ and $w \in r_0$ with appointment for a vertex $v \in q_0$ the vertex $(u \circ v) \circ w$ and for an edge $\beta : v \to v'$ appointing the new edge

$$(1_u \circ \beta) \circ 1_w : (u \circ v) \circ w \to (u \circ v') \circ w \ .$$

For an edges $\alpha : u' \to u$ and $\gamma : w \to w'$ we appoint the transform between such functors

$$(\alpha \circ 1_v) \circ \gamma : (u \circ v) \circ w \to (u' \circ v) \circ w' \ .$$

This can be expressed as a new functor to the new exponential space

$$X^{op}(a;b) \times X(c;d) \Longrightarrow [(\bigcirc \circ X(a;c)) \circ \bigcirc; X(b;d)] \subset [X(a;c); X(b;d)] \ .$$

It will be a partial form got from the form

$$X^{op}(a;b) \times X(c;d) \times X(a;c) \Longrightarrow (X^{op}(a;b) \circ X(a;c)) \circ X(c;d) \subset X(b;d) \ .$$

There is no surprise that the same result will be got also for posthorizontal composition of edges in value category $V$.

$$A(\alpha) *^- A(\gamma) = (T_v \star A(\gamma)) \cdot (A(\alpha) \star T_{w'}) .$$

For the postdecomposable form

$$\text{Hom}_X^+(X^{op}(a;b) \times X(c;d)) = \text{Hom}_X(1_a(*) \otimes X(c;d)) \cdot \text{Hom}_X(X^{op}(a;b) \times 1_d(*))$$

we take one partial form

$$\text{Hom}_X(1_a(*) \times X(c;d)) = [\bigcirc \circ X(a;c); X(a;d)] \subset [X(a;c); X(a;d)]$$

$$\begin{array}{ccc} & & d \\ & \nearrow & \uparrow \\ a & \to & c \end{array}$$

and another partial form

$$\text{Hom}_X(X^{op}(a;b) \times 1_d(*)) = [\bigcirc \circ X(c;d); X(b;d)]$$

$$\begin{array}{ccc} b & \to & d \\ \downarrow & \nearrow & \\ a & & \end{array}$$

The defined postdecomposable form can be explain by diagram

$$\begin{array}{ccc} b & \to & d \\ \downarrow & \nearrow & \uparrow \\ a & \to & c \end{array}$$

For the values category $V = OGrph$ the postdecomposable form can be expressed again as a functor to the new exponential space

$$X^{op}(a;b) \times X(c;d) \Longrightarrow [\bigcirc \circ (X(a;c) \circ \bigcirc); X(b;d)] \subset [X(a;c); X(b;d)] .$$

It will be a partial form got from the form

$$X^{op}(a;b) \times X(c;d) \times X(a;c) \Longrightarrow X^{op}(a;b) \circ (X(a;c) \circ X(c;d)) \subset X(b;d) .$$

Having associativity arrows in taken original $V$-graph $X$

$$S : X(a;b) \circ (X(b;c) \circ X(c;d)) \to (X(a;b) \circ (X(b;c)) \circ X(c;d) ,$$

we get an arrow also between predecomposable and postdecomposable forms

$$\text{Hom}_X^+ =\!\!\Rightarrow \text{Hom}_X^-$$

which are partial forms for the correspondent forms. The same is true for an opposite arrow, or for isomorphic arrows.

Now we are ready to prove the Yoneda lemma. It will be proved for arbitrary original $V$-graph $X$ with contact arrows mappings

$$X(a;b) \otimes X(b;c) \to X(a;b) \circ X(b;c) \subset X(a;c) \ .$$

Let we have an $V$-graphs transport from the opposite $V$-graph $X^{op}$ to the values category $V$ considered as $V$-graph

$$F : X^{op} \Longrightarrow V^e \ .$$

It will be called a *character* over taken original $V$-graph $X$. The representable character is defined by any vertex $x \in X_0$

$$X(-;x) : X^{op} \Longrightarrow V^e \ .$$

Such appointment will be understood as a Yoneda mapping

$$(Y_X)_0 : X_0 \to (V^{X^{op}})_0 \ .$$

It can be extended as graphs transport

$$Y_X : X \Longrightarrow V^{X^{op}}$$

and this extension will be *fullfaithful imbedding*. It will be called a *Yoneda imbedding*. Usually mathematicians were interesting when such embedding can be corestricted in the set of natural transforms $V^{(X^{op})}$. We extend such problematic with the set of continuous transforms $V_-^{(X^{op})}$.

The Yoneda lemma will describe the set of transforms between representable character $X(-;x)$ and other character $F : X^{op} \Longrightarrow V$ which is an original $V$-graphs transport, i. e. we get identity mapping

$$F(i_x) = \mathtt{Id}_{F_0(x)} : F_0(x) \to F_0(x) \ .$$

**Proposition 13.**  *For an original $V$-graphs transport $F : X^{op} \Longrightarrow V^e$ we get a fullfaithful mapping from the value object $F_0(x)$ to the exponential space of transforms between representable character $X(-;x)$ and taken $V$-graphs transport $F$*

$$F_0(x) \to [X(-;x); F] \ .$$

**Proof:**  For a $V$-graphs transport $F : X^{op} \Longrightarrow V^e$ we have a mapping

$$X(a;x) \to [F_0(x); F_0(x)] \ .$$

It is a name mapping for predecomposable form

$$X(a;x) \otimes F_0(x) \to F(x) \ .$$

The dual name form will be
$$F_0(x) \to [X(a;x); F_0(a)] .$$
So we have got a mapping to the product
$$\rho : F_0(x) \to \prod_{a \in X_0} [X(a;x); F_0(a)] .$$
The inverse mapping is provided by projection to the factor
$$\prod_{x \in X_0} [X(a;x); F_0(a)] \to [X(x;x); F_0(x)]$$
and taking image of unit arrow $i_x(*) \subset X(x;x)$
$$[i_x(*); F_0(x)] \subset F_0(x) .$$

By our assumption taken functor $F : X^{op} \Longrightarrow V^e$ is an original $V$-graphs transport, therefore it maintains the choosing of unit arrows, therefore we have got the composition as identity mapping
$$F_0(x) \to \prod_{a \in X} [X(a;x); F_0(a)] \to [i_x(x); F_0(x)] \subset F_0(x) .$$
Otherwise we have also identic mapping
$$\prod_{a \in X_0} [X(a;x); F_0(a)] \to [i_x(*); F_(x)] \xrightarrow{\rho} \prod_{a \in X_0} [X(a;x); F_0(a)] .$$
It is enough to notice that every factor is defined by the value object $F_0(x)$. $\square$

Let we have an original $V$-graph $X$ with upper associative contact arrows mappings, i. e. we have associative arrow between mappings
$$X(a;b) \circ (X(b;c)) \circ X(c;x)) \to (X(a;b) \circ X(b;c)) \circ X(c;x) .$$
Every representable character now will be continuous one
$$\begin{array}{cc} X(a;b) \otimes X(b;c) & \xrightarrow{\circ} \\ X(-;x) \downarrow \otimes X(-;x) & \\ [\bigcirc \circ X(b;c); X(a;c)] \otimes [\bigcirc \circ X(c:x); X(b;x)] & \to \end{array}$$

$$\begin{array}{cc} \to & X(a;c) \\ & \downarrow X(-;x) \\ \to & [\bigcirc \circ X(c;x); X(b;x)] \cdot [\bigcirc \circ X(b;x); X(a;x)] \\ & \nearrow \end{array}$$

We need only to calculate the double evaluation for composition of contact arrows mappings in $V$-graph $V^e$
$$X(b;c) \otimes ([\bigcirc \circ X(c;x); X(b;x)] \cdot [\bigcirc \circ X(b;x); X(a;x)] \xrightarrow{\mathbf{ev}}$$
$$(X(b;c) \circ X(c;x)) \otimes [\bigcirc \circ X(b;x); X(a;x)] \xrightarrow{\mathbf{ev}}$$
$$X(a;b) \circ (X(b;c) \circ X(c;x) \subset X(a;x) .$$
For it we have associativity arrow to another predecomposable form
$$X(a;b) \circ (X(b;c) \circ X(c;x) \to (X(a;b) \circ X(b;c)) \circ X(c;x) \subset X(a;x) .$$

We shall improve the former proposition for continuous functor between original $V$-graphs $F : X^{op} \Longrightarrow V^e$.

**Proposition 14.** *In original $V$-graph $X$ with upper associative contact arrows mappings*

$$X(a;b) \circ (X(b;c)) \circ X(c;x)) \to (X(a;b) \circ X(b;c)) \circ X(c;x)$$

*the set of transforms between a representable functor $X(-;x) : X^{op} \Longrightarrow V^e$ and a continuous original $V$-graphs transport $F : Xop \Longrightarrow V^e$ coincides with the set of continuous transforms $V_{-}^{(X^{op})}(X(-;x);F)$.*

**Proof:** We need to check continuous inequality between two forms One is get with preevaluation transport

$$X^{op}(a;b) \otimes \prod_{a\in X_0} [X(a;x); F_0(a)] \to [\bigcirc \circ X(a;x); X(b;x)] \cdot [X(b;x); F_0(b)]$$

and another is get with postevaluation transport

$$X^{op}(a;b) \otimes \prod_{a\in X_0} [X(a;x); F_0(a)] \to [X(a;x); F_0(a)] \cdot [F_0(a); F_0(b)] \ .$$

The predecomposable form will be calculated with the double evaluation

$$X^{op}(a;b) \otimes [\bigcirc \circ X(a;x); X(b;x)] \cdot [X(b;x); F_0(b)] \xrightarrow{\mathbf{ev}}$$
$$(X(b;a) \circ X(a;x)) \otimes [X(b;x); F_0(b)] \xrightarrow{[(X(b;a)\circ X(a;x); F_0(b)]} \ .$$

It was calculated as

$$F_{x,b}(X(b;a) \circ X(a;x)) \subset [F(x); F(b)] \ .$$

The postdecomposable form will be calculated also with double evaluation

$$X^{op}(a;b) \otimes [X(a;x); F_0(a)] \cdot [F_0(a); F_0(b)] \to$$
$$[X(a;x); F_0(a)] \cdot X^{op}(a;b) \otimes [F_0(a); F_0(b)]_)$$

It was calculated as

$$F_{x,a}(X(a;x)) \cdot F_{a;b}(X(b;a)) \subset F_{x,b}(X(b;x) \ .$$

For the continuous original $V$-graphs transport $F$ we have an arrow

$$F_{x,a)(X(a;x))\cdot F_{a;b}(X(b;a))=\to F(b;x}(X(b;a) \circ X(a;x)) \ .$$

$\square$

As corollary we get fullfaithful imbedding of original $V$-graphs transport with upper associative contact arrows mappings into the set of continuous transforms between representable characters, i. e. we have isomorphic objects in values category $V$

$$X(x;y) \to V_{-}^{(X^{op})}(X(-;x); X(-;y)) \ .$$

As partial case we get this result also for the $V$-enriched categories, as they are original $V$-graphs with associative contact arrows mappings.

The isometry of metric spaces will be example of such isomorphism for $V$-enriched categories. We shall remark that for morphism of $V$-enriched categories such property of isometry is more weak than to be isomorphic.

We append some theoretical constructions which will provide more insight about the set of natural or continuous transforms. We shall use a translation $V \Longrightarrow Set$ of values category to the category of sets. The image will be called a *pointed category*. We shall denote it as an exponential space over final object $V^*$. For the object $r \in V_0$ we take the set of its *points*

$$r^* = \{* \to r\} \ .$$

A point is understood as a mapping from category's final object. For the mapping $u : a \to b$ between the objects $a, b \in V_0$ we get a mapping between the sets of points $u^* : a^* \to b^*$.

$$\begin{array}{ccc} * & & \\ \downarrow & \searrow & \\ a & \xrightarrow[u]{} & b \end{array}$$

For a point $* \to a$ it appoints another point $* \to b$. Obviously such appointing of mappings maintains an identity arrow and composition of arrows, therefore we get indeed a functor $V \Longrightarrow Set$.

For the biproduct $r \otimes s$ the points will be defined by biproduct of correspondent mappings

$$* =\to * \otimes * \to r \otimes s \ .$$

So we get a mapping of points

$$r^* \times s^* \to (r \otimes s)^* \ .$$

For two mappings between the objects $u : a \to a'$ and $v : b \to b'$ we get a biproduct $u \otimes v : a \otimes b \to a' \otimes b'$. Correspondingly for points $* \to a$ and $* \to b$ we get the point in the biproduct $* \to * \otimes * \to a \otimes b$. So we defined the mapping $a^* \otimes b^* \to (a \otimes b)^*$ from the product of points sets $a^* \times b^*$ to the set of points in the biproduct of objects $(a \otimes b)^*$.

We have got a *natural translation* with isomorphic comparison arrows

$$t_{a,b} : r^* \otimes b^* =\to (r \otimes t)^* \ .$$

In the final object $* \in V_0$ we get the unique point $* \to *$, The onepoint set is a final in the category of sets $Set$, therefore the translation of pointed category maintains the final object.

For the partial category of *Birkhoff algebras* in the category of sets $Set$ the final object will be onepoint set and for each such algebra $X \in Set_0$ the points coincide with members of the set carrying this algebra. Every transport of algebras is defined by some mapping between the carrying sets. So we get that functor of pointed category is identic for this category.

In the category of Banach spaces Ban the objects are Banach spaces with arrows being nonexpanding linear mappings $u : X \to Y$

$$|u(x)| \leq |x| .$$

The final object will be any 1-dimensional space $x \cdot R$ with norm $|x| \leq 1$. So the points in Banach space $X$ will be any 1-dimensional subspace $x \cdot R$ with the points $x$ of unit ball $|x| \leq 1$. The set of points of Banach space for the category of nonexpanding linear mappings can be identified with the unit ball

$$X^* = \{x \in X : |x| \leq 1\} \subset X .$$

However in bigger category of arbitrary bounded linear mappings $u : X \to Y$ the points of Banach space will coincide with the whole space $X^* = X$.

For us will be interesting to note that in the scalar category $[0, \infty]^+$ we have the final object $0 \in [0; \infty]$ and only one point $0 \geq 0$.

In the category of graphs $Grph$ the final object will be the graph with the unique vertex and the unique edge. For another graph $X$ all vertexes $x \in X_0$ with some loop $x \to x$ defines some points. If we have different loops, we shall have different points with the same vertex.

The situation will be simpler in the category $OGrph$ of original graphs with their transports. The original graphs transport must maintain chosen unit arrows, so the points coincide with the vertexes of taken graph

$$(X^*)_0 = X_0 .$$

The transform $\alpha : S \to T$ between two parallel graphs transports $S, T : X \Longrightarrow Y$ is defined by collection of edges

$$\alpha_a \in Y(S_0(a); T_0(a)) .$$

The transform between two points $a, b : * \Longrightarrow X$ will be defined by one edge $\alpha \in (X(a; b)$. And all edges are possible. So we get the same graph

$$X^* = X .$$

For $V$-enriched graphs $X$ we haven't any particular transform between two $V$-graphs transforms. We have defined only the set of transforms between two graphs transports $\phi, \psi : X \Longrightarrow Y$

$$Y^X(\phi; \psi) = \prod_{a \in X_0} Y(\phi_0(a); \psi_0(a)) .$$

The compromise is achieved by defining **2**-category over $Grph_V$. Between parallel graphs transports $\phi, \psi : X \Longrightarrow Y$ we are chosen the particular transforms. We understand the 2-category only as possibility to define order relation in the set of parallel graphs transports $Grph_V(X;Y)$. The transform will be defined by arbitrary collection of edges in the target space

$$\alpha = \prod_{a \in X_0} \alpha_a(*) \subset \prod_{a \in X_0} Y(\phi_0(a); \psi_0(a)) \ .$$

The final object in the category of $V$-graphs $Grph_V$ will be arbitrary one-point set with unique edge $* \in V_0$. For original $V$-graph $X$ we again have the same set of vertexes

$$(X^*)_0 = X_0 \ .$$

The edge between two points $a \in V_0$ and $b \in V_0$ will be a transforms between corespondent original $V$-graphs transports $a : * \Longrightarrow X$ and $b : * \Longrightarrow X$. It is defined by some edge $\alpha \in X(a;b)$. We get a new graph. For an original $V$-graphs transport between two original $V$-graphs $f : X \Longrightarrow Y$ we get mapping of particular edge between two points $\alpha \in X(a;b)$ to the new edge between new points

$$f(\alpha) \in Y(f_0(a); f_0(b)) \ .$$

So we get superfunctor $OGrph_V \Longrightarrow OGrph$ defined by translation of pointed category $V \Longrightarrow V^*$. Appointed new graph will be called an *underlying graph*.

With translation of underlying category $V \Longrightarrow V^*$ we get a new $V^*$-graph $X^*$ with the same vertexes $X_0^* = X_0$ but remaining only points in edges sets

$$X^*(a;b) = \{* \to X(a;b)\} \ .$$

For $V$-graphs transport $f : X \Longrightarrow Y$ we get mappings in the category of sets

$$(f^*)_{a,b} : X(a;b)^* \to Y(f_0(a); f_0(b))^* \ .$$

The using of such notation could be dangerous, but it will be useful for us. Therefore we used two different names for these different situations.

The set of all transform $\alpha : f \to g$ between two parallel $V$-graphs transports $f : X \Longrightarrow Y$ and $g : X \Longrightarrow Y$ will be defined by collection of chosen edges

$$\alpha_a(*) \subset Y(f_0(a); g_0(a)) \uparrow a \in X_0 \ .$$

Let we have $V$-graph $Y$ with contact arrows mappings defining functor over pullbacked biproduct

$$\circ : Y \otimes_p Y \Longrightarrow Y \ .$$

Then we get the postevaluation functor

$$\mathbf{ev}_Y^- : X \otimes Y^X \Longrightarrow Y$$

defined with postcontact arrows mappings

$$\mathbf{ev}_Y^- : X(a;b) \otimes Y^X(\phi;\psi) \to Y(\phi_0(a); \psi_0(a)) \circ \psi(X(a;b)) \subset Y(\phi_0(a); \psi_0(b)) \ .$$

Also we have the preevaluation functor

$$\mathbf{ev}_Y^+ : X \otimes Y^X \Longrightarrow Y$$

defined with precontact arrows mappings

$$\mathbf{ev}_Y^+ : X(a;b) \otimes Y^X(\phi;\psi) \to \phi(X(a;b)) \circ Y(\phi_0(b);\psi_0(b)) \subset Y(\phi_0(a);\psi_0(b)) \ .$$

In the values category $V$ for arbitrary objects $p \in V_0$ and $r \in V_0$ we have defined an exponential space $r^p \in V_0$ Using contact arrows mappings in any object $r \in V_0$ we have defined also contact arrows mappings between exponential spaces

$$q^p \times r^q \to r^p \ .$$

The underlying graph coincides with taken values category

$$[p;r]^* = V(p;r) \ .$$

In **G. Kelly 1982** p. 14 author uses also the equality

$$p = [*;p] \ ,$$

i. e. every object can be get as exponential space over the set of its points. It is fulfilled in the scalar category $[0,\infty]^+$

$$p = (p - 0)^+ \ .$$

But it isn't universal even in monoidal closed categories. So we shall try to apply such property more carefully.

We demanded isomorphic arrow for the right projection $\rho : * \otimes p \to p$. For almost coadjoint biproduct the inverse arrow will be provided as a biproduct with the uniquely defined arrow to the final object $* \in V_0$

$$p \xrightarrow{!} * \ .$$

Then we can get an arrow begun with diagonal arrow

$$p \to p \otimes p \xrightarrow{! \otimes 1_p} * \otimes p \ .$$

In values category $V$ we have associative composition of arrows, therefore the isomorphic arrow $\pi : * \otimes r \to r$ provides a natural transform of representable characters defined with collection of arrows set mappings

$$\alpha_a : V(a; * \otimes p) \to V(a;p) \ .$$

We need to check that such natural transform is bijective. The inverse natural transform will be defined with inverse arrows

$$\beta_a : V(a;p) \to V(a; * \otimes p) \ .$$

In values category we can apply associative property for arrows composition. So we get commuting diagrams for arbitrary arrow $f \in V(a; x)$

$$\begin{array}{ccccc} x & \xrightarrow{\theta_x} & * \otimes x & \xrightarrow{\pi_x} & x \\ f \uparrow & \nearrow & \uparrow 1_* \otimes f & \nearrow & \uparrow f \\ a & \xrightarrow{\theta_a} & * \otimes x & \xrightarrow{\pi_a} & a \end{array}$$

Also we have commuting diagram for arbitrary arrow $g : a \to * \otimes x$.

We continue with introduction some notion which were essential for G. Kelley way of proving Yoneda lemma for $V$-enriched categories.

A natural transform $\alpha : F \to g$ between two graphs transports $f, g X \Longrightarrow Y$ is defined with collection of edges $\alpha_a \in Y(f_0(a); g_0(a)) \uparrow a \in X_0$ and demands commuting diagram for each edge $u \in X(a; b)$

$$\begin{array}{ccc} f_0(a) & \xrightarrow{\alpha_a} & g_0(a) \\ f(u) \downarrow & & \downarrow g(u) \\ f_0(b) & \xrightarrow{\alpha_b} & g_0(b) \end{array}$$

A natural transform between two $V$-graphs transports needs more sophisticate definition. Such transform is defined by taking some edges as mappings from final object of values category $* \in V_0$

$$* \xrightarrow{\alpha_a} Y(f_0(a); g_0(a)) \ .$$

Now we demand the commuting diagrams

$$\begin{array}{ccc} * \otimes X(a; b) & \xrightarrow{\alpha_a \otimes g} & Y(f_0(a); g_0(a)) \otimes Y(g_0(a); g_0(b)) \\ \nearrow & & \searrow \circ \\ X(a; b) & & Y(f_0(a); g_0(b)) \\ \searrow & & \nearrow \circ \\ X(a; b) \otimes * & \xrightarrow[f \otimes \alpha_b]{} & Y(f_0(a); f_0(b)) \otimes Y(f_0(b); g_0(b)) \end{array}$$

We shall rewrite such diagram more compactly.

$$\begin{array}{ccc} X(a; b) & \xrightarrow{g} & Y(g_0(a); g_0(b)) \\ f \downarrow & & \downarrow Y(\alpha_a; 1) \\ Y(f_0(a); f_0(b)) & \xrightarrow[Y(1; \alpha_b)]{} & Y(f_0(a); g_0(B)) \end{array}$$

The mapping $Y(1; \beta) : Y(x; y) \to Y(x; z)$ is defined as changing of target space with an edge $\beta : * \to Y(y; z)$

$$\beta_* : Y(x; y) \to Y(x; y) \otimes * \xrightarrow{1 \otimes \beta} Y(x; y) \otimes Y(y; z) \xrightarrow{\circ} Y(x; z) \ .$$

The mapping $Y(\alpha; 1) : Y(y; z) \to Y(x; z)$ is defined as changing of source space with an edge $\alpha : * \to Y(x; y)$

$$\alpha^* : Y(z; z) \to * \otimes Y(y; z) \xrightarrow{\alpha \otimes 1} Y(x; y) \otimes Y(y; z) \xrightarrow{\circ} Y(x; z) \ .$$

Between two $V$-graphs transports $X(-;x), F : X^{op} \Longrightarrow V^e$ the set of natural transforms is defined as an exponential space of natural transforms $K = F^{(X(-;a))}$ identified with the equalizer as partial object within the space of all transforms between taken functors

$$\to K \to \prod_{a \in X_0} [X(a;x); F_0(x)] \Rightarrow \prod_{a,b \in X_0} [X(a;x); F_0(b)]^{X^{op}(a;b)} \ .$$

The two mappings were defined as names for postevaluation and preevaluation mappings with the set of edges $X^{op}(a;b)$

$$X^{op}(a;b) \otimes \prod_{a \in X_0} [X(a;x); F_0(a)] \to [X(a;x); F_0(a)] \cdot F_{a,b}(X^{op}(a;b)) \ ,$$

$$X^{op}(a;b) \otimes \prod_{a \in X_0} [X(a;x); F_0(a)] \to [X(a;x); X(b;x)] \cdot [X(b;x); F_0(b)] \ .$$

Such equalizer categorically is expressed as an end for the *bifunctor*

$$T(a,b) = [X(-;x); K] : X^{op} \times X \Longrightarrow V^e \ ,$$

clr. **S. MacLane 1971** part IX. Special limits. It is a partial object embedded in the product of all values $T(a,a)$

$$\to \int_{a \in X_0} T(a,a) \to \prod_{a \in X_0} T(a,a) \Rightarrow \prod_{a,b \in X_0} T(a,b)^{X^{op}(a;b)} \ .$$

A bifunctor $T : X^{op} \times X \Longrightarrow V^e$ is defined by mapping over the set of vertexes

$$T(a,b) = V(X(a;x); F_0(b))$$

and is a $V$-graphs functor for each argument. In our case one functor is

$$T(a,-) : X^{op} \Longrightarrow V^e$$

with edges set mappings

$$X^{op}(a;b) \to [X(a;x); F_0(a)] \cdot F_{a,b}(X^{op}(a;b)) \subset [X(a;x); F_0(b)]$$

We express such mappings with a diagram

$$\begin{array}{ccc} X(a;x) & \to & F_0(a) \\ & \searrow & \downarrow X^{op}(a;b) \\ & & F_0(b) \end{array}$$

The contact arrows mappings in original $V$-graph $V^e$ was defined as pre-decomposable form

$$q^p \otimes r^q \to q^p \cdot r^q \subset r^p$$

so we get a postdecomposable form

$$X^{op}(a;b) \otimes T(a,a) \to T(a,a) \cdot F_{a,b}(X^{op}(a;b)) \subset T(a,b)$$

and its name form
$$\sigma_{a;b} : T(a,a) \to T(a;b)^{X^{op}(a;b)} .$$

Another functor is
$$T(-,b) : X \Longrightarrow V^e$$

with edges set mappings

$$X^{op}(a;b) \to [\bigcirc \circ X(a;x); X(b;x)] \cdot [X(b;x); F_0(b)] \subset [X(a;x); K_0(b)]$$

We express such mappings with a diagram

$$\begin{array}{ccc} & X(a;x) & \\ X^{op}(a;b) \downarrow & \searrow & \\ & X(b;x) & \to & F(b) \end{array}$$

It defines another predecomposable form

$$X^{op}(a;b) \otimes T(b,b) \to [\bigcirc \circ X(a;x); X(b;x)] \cdot T(b,b)$$

with its name form
$$\rho_{a,b} : T(b,b) \to T(a,b)^{X^{op}(a;b)} .$$

So the end $K$ will be defined as initial object for all commuting diagrams

$$\begin{array}{ccc} \searrow & & \\ K & \xrightarrow{\lambda_a} & T(a;a) \\ \lambda_b \downarrow & & \downarrow \rho_{a,b} \\ T(b,b) & \xrightarrow[\sigma_{a,b}]{} & [X^{op}(a;b); T(a,b)] \end{array}$$

For the postdecomposable form

$$X^{op}(a;b) \otimes T(a,a) \to T(a,b)$$

we have bijection between *dual name forms* defined as mappings

$$T(a,a) \to [X^{op}(a;b); T(a;b)] , \quad X^{op}(a;b) \to [T(a,a); T(a,b)] .$$

Such bijection exists even without existence of predecomposable or postdecomposable form as a functor over biproduct $X^{op}(a;b) \otimes T(a,a)$. For name forms we need to have only partial forms of arbitrary bifunctor

$$T : X^{op}(a;b) \times T(a,a) \Longrightarrow T(a,b) .$$

However for the next property existence of predecomposable or postdecomposable form is essential.

We shall check that a triangular diagram for source changing between the name mappings

$$\begin{array}{ccc} K & \xrightarrow{\lambda_a} & T(a,a) \\ {}_{\phi_1}\searrow & & \downarrow \psi_1 \\ & [X^{op}(a;b); T(a,b)] & \end{array}$$

provides the triangular diagram for target changing between the dual name mappings

$$\begin{array}{ccc} X^{op}(a;b) & \xrightarrow{\psi_2} & [T(a,a); T(a,b)] \\ {}_{\phi_2}\searrow & & \downarrow [\lambda_a; 1_{T(a,b)}] \\ & [K; T(a,b)] & \end{array}$$

and vice versa.

**Proposition 15.** *Let we have a pair of dual name forms provided by a pair of predecomposable forms. Then an existing source changing between name forms provides the target changing for the dual name forms. And vice versa an existing target changing between name forms provides the source changing for the dual name forms.*

**Proof:** Let we have a pair of predecomposable forms

$$\phi : X(a;b) \otimes K \to T(a;b) \ , \quad \psi : X(a;b) \otimes T(a,a) \to T(a,b) \ .$$

If we have source changing transform $\theta : K \to T(a,a)$ between name forms

$$\begin{array}{ccc} K & \xrightarrow{\theta} & T(a,a) \\ {}_{\phi_1}\searrow & & \downarrow \psi_1 \\ & [X(a;b); T(a,b)] & \end{array}$$

then we get the source changing transform between predecomposable forms

$$\begin{array}{ccccc} X(a;b) \otimes K & \xrightarrow{1_{X(a;b)}\otimes\phi_1} & X(a;b) \otimes [X(a;b); T(a,b)] & \xrightarrow{\mathbf{ev}_{T(a,b)}} & T(a,b) \\ {}_{1_{X(a;b)}\otimes\theta}\downarrow & \nearrow & {}_{1_{X(a;b)}\otimes\psi_1} & & \\ X(a;b) \otimes T(a,a) & & & & \end{array}$$

For the dual name forms we get commuting diagram

$$\begin{array}{ccccc} & & [T(a,a); T(a,a) \otimes X(a;b)] & \xrightarrow{[1_{T(a,b)};\psi]} & T(a,a); T(a,b)] \\ & {}^{\lambda_{X(a;b)}}\nearrow & {}_{[\theta;1_{T(a,a)\otimes X(a;b)}]}\downarrow & & \downarrow {}^{[\theta;1_{T(a,b)}]} \\ X(a;b) & & [K; T(a,a) \otimes X(a;b)] & \xrightarrow{[1_K;\psi]} & [K; T(a;b)] \\ & {}_{\lambda_{X(a;b)}}\searrow & {}_{[1_K;\theta\otimes 1_{T(a,b)}]}\uparrow & \nearrow {}_{[K;\phi]} & \\ & & [K; K \otimes X(a,b)] & & \end{array}$$

Otherwise for the target changing of dual name mappings

$$\begin{array}{ccc}
X(a;b) & \xrightarrow{\psi_2} & [T(a,a);T(a,b)] \\
& \phi_2 \searrow & \downarrow {[\theta;1_{T(a;b)}]} \\
& & [K;T(a,b)]
\end{array}$$

we get the source changing of predecomposable forms

$$\begin{array}{ccccc}
T(a;a) \otimes X(a;b) & \xrightarrow{1_{T(a,a)} \otimes \psi_2} & T(a,a) \otimes [T(a,a);T(a,b)] & & \\
\theta \otimes 1_{X(a;b)} \uparrow & & \theta \otimes 1_{[T(a,a);T(a,b)]} \uparrow & \searrow \mathbf{ev}_{T(a,b)} & \\
K \otimes X(a;b) & \xrightarrow{1_K \otimes \psi_2} & K \otimes [T(a,a);T(a,b)] & & T(a,b) \\
& \phi_2 \searrow & 1_k \otimes [\theta;1_{T(a,b)}] \downarrow & \nearrow \mathbf{ev}_{T(a,b)} & \\
& & K \otimes [K;T(a,b)] & &
\end{array}$$

This produces the source changing transform between dual name forms

$$\begin{array}{ccccc}
K & \xrightarrow{\lambda_K} & [X(a;b); K \otimes X(a;b)] & \xrightarrow{\phi} & [X(a;b);T(a,b)] \\
\downarrow \theta & & \downarrow {[1_{X(a;b)}; \theta \otimes 1_{X(a;b)}]} & \nearrow \psi & \\
T(a,a) & \xrightarrow{\lambda_{T(a,a)}} & [X(a;b); T(a,a) \otimes X(a;b)] & &
\end{array}$$

We must remark that for the dual name forms we used the dual sections arrows
$$\lambda_K : K \to [X(a;b); K \otimes X(a;b)] \ ,$$
$$\lambda_{T(a,a)} : T(a,a) \to [X(a;b); T(a,a) \otimes X(a;b)] \ .$$

Also we used the dual evaluation arrows
$$\mathbf{ev}_{T(a;b)} : [X(a;b); T(a;b) \otimes X(a;b)] \to T(a;b) \ .$$

$\square$

As a consequence we get a correspondence between commuting diagrams. The commuting diagram of name mappings

$$\begin{array}{ccc}
K & \xrightarrow{\lambda_a} & T(a,a) \\
\lambda_b \downarrow & \searrow & \downarrow \sigma_{a,b} \\
T(b,b) & \xrightarrow{\rho_{a,b}} & [X^{op}(a;b); T(a,b)]
\end{array}$$

provides also the commuting diagram between dual name mappings

$$\begin{array}{ccc}
X^{op}(a;b) & \xrightarrow{\bar{\sigma}_{a,b}} & [T(a,a);T(a,b)] \\
\bar{\rho}_{a,b} \downarrow & \searrow & \downarrow {[\lambda_a; 1_{T(a;b)}]} \\
[T(b,b);T(a,b)] & \xrightarrow{[\lambda_b; 1_{T(a,b)}]} & [K;T(a,b)]
\end{array}$$

and vice versa.

2023.05.28